\documentclass[journal]{IEEEtran}

\usepackage[utf8]{inputenc}
\usepackage[english]{babel}
\usepackage{amsmath}

\usepackage{subfig}
\usepackage{tikz}
\tikzset{mynode/.style={draw,solid,circle,inner sep=1pt}}

\usepackage{booktabs}
\usepackage{rotating,tabularx}
\usepackage{multirow}
\usepackage{multicol}
\usepackage{color, colortbl}

\usepackage{graphicx,caption}
\graphicspath{{./figure/}}
\usepackage{bm}
\usepackage{amsfonts}
\usepackage{amsmath}
\usepackage{amssymb}
\usepackage{amsbsy}
\usepackage{mathrsfs} 
\usepackage{mathabx}
\DeclareMathSymbol{\shortminus}{\mathbin}{AMSa}{"39}

\usepackage{verbatim}
\usepackage{blindtext}
\usepackage{tcolorbox}
\usepackage{scalerel}
\usepackage{lipsum}
\usepackage{url}
\usepackage{adjustbox}
\usepackage{soul}
\usepackage{cite}
\usepackage{tikz}
\usepackage{algpseudocode}

\usepackage[wby]{callouts}
\usetikzlibrary{shapes.multipart}

\usepackage{array}
\usepackage{makecell}

\allowdisplaybreaks
\usepackage{amsthm}

\usepackage[linesnumbered, ruled]{algorithm2e}
\SetKw{Continue}{continue}
\SetKw{And}{and}
\SetKw{Elif}{elif}
\SetKw{True}{True}
\SetKw{False}{False}
\SetKw{To}{to}
\SetKw{Def}{def}
\SetKw{Not}{not}

\SetKwComment{Comment}{\color{green!50!black}// }{}

\SetKwProg{Function}{function}{}{}

\theoremstyle{definition}

\theoremstyle{remark}

\newcommand*{\mrn}{\textcolor{black}}

\usepackage{xargs}                      
\newcommandx{\unsure}[2][1=]{\todo[linecolor=red,backgroundcolor=red!25,bordercolor=red,#1]{#2}}
\newcommandx{\change}[2][1=]{\todo[linecolor=blue,backgroundcolor=blue!25,bordercolor=blue,#1]{#2}}
\newcommandx{\info}[2][1=]{\todo[linecolor=OliveGreen,backgroundcolor=OliveGreen!25,bordercolor=OliveGreen,#1]{#2}}
\newcommandx{\improvement}[2][1=]{\todo[linecolor=Plum,backgroundcolor=Plum!25,bordercolor=Plum,#1]{#2}}
\newcommandx{\thiswillnotshow}[2][1=]{\todo[disable,#1]{#2}}
\presetkeys{todonotes}{inline}{}

\renewcommand{\baselinestretch}{.955}
\newif\ifarxiv
\arxivtrue
\title{Tightening QC Relaxations of AC Optimal Power Flow through Improved Linear Convex Envelopes}

\author{Mohammad Rasoul Narimani,$^{\ast}$ Daniel K. Molzahn,$^{\dagger}$ Katherine R. Davis,$^{**}$ and Mariesa L. Crow$^{\ddagger}$
\thanks{${\ast}$: Department of Electrical and Computer Engineering, California State University Northridge (CSUN). Rasoul.narimani@csun.edu. Support from NSF contract \#2308498.}% 
\thanks{$^{\dagger}$: School of Electrical and Computer Engineering, Georgia Institute of Technology. molzahn@gatech.edu. Support from NSF contract \#2023140.}% 
\thanks{${**}$: Electrical and Computer Engineering Department, Texas A\&M University. kdavis@tamu.edu.}
\thanks{${\ddagger}$: Electrical and Computer Engineering Department, Missouri University of Science and Technology. crow@mst.edu.}
}%

\begin{document}

\maketitle
\begin{abstract}
AC optimal power flow (AC OPF) is a fundamental problem in power system operations. Accurately modeling the network physics via the AC power flow equations makes AC OPF a challenging nonconvex problem. To search for global optima, recent research has developed various convex relaxations that bound the optimal objective values of AC OPF problems. The QC relaxation convexifies the AC OPF problem by enclosing the non-convex terms within convex envelopes. The QC relaxation's accuracy strongly depends on the tightness of these envelopes. This paper proposes two improvements for tightening QC relaxations of OPF problems. We first consider a particular nonlinear function whose projections are the nonlinear expressions appearing in the polar representation of the power flow equations. We construct a polytope-shaped convex envelope around this nonlinear function and derive convex expressions for the nonlinear terms using its projections. Second, we use sine and cosine expression properties, along with changes in their curvature, to tighten this convex envelope. We also propose a coordinate transformation to tighten the envelope by rotating power flow equations based on individual bus-specific angles. We compare these enhancements to a state-of-the-art QC relaxation method using PGLib-OPF test cases, revealing improved optimality gaps in 68\% of the cases.    
\end{abstract}


\IEEEpeerreviewmaketitle

\section{Introduction}
\IEEEPARstart{T}{he} optimal power flow (OPF) problem seeks an operating point that optimizes a specified objective function (often generation cost minimization) subject to constraints from the network physics and engineering limits.
Using the nonlinear AC power flow model to accurately represent the power flow physics results in the AC~OPF problem, which is non-convex, may have multiple local optima~\cite{bukhsh_tps}, and is generally NP-Hard~\cite{bienstock2015nphard}. 

{\color{black}Since first being formulated by Carpentier in 1962~\cite{carpentier}, a wide variety of optimization algorithms have been applied to the OPF problem~\cite{opf_litreview1993IandII,ferc4,abdi2017}. Much of this research has focused on algorithms for obtaining locally optimal or approximate OPF solutions. Recently, many convex relaxation techniques have been applied to OPF problems to obtain bounds on the optimal objective values, certify infeasibility, and in some cases, achieve globally optimal solutions~\cite{molzahn2019fnt}.

As a key metric for solution quality, the objective value bounds obtained via convex relaxations characterize how close a local solution is to being globally optimal. Thus, local algorithms and relaxations are often used together in spatial branch-and-bound algorithms to solve nonlinear programs (NLPs) and mixed-integer nonlinear programs (MINLPs)~\cite{burer2012}. With nonlinear constraints modeling AC power flows, many problems relevant to power systems take the form of NLPs and MINLPs (e.g., OPF~\cite{gopinath2022}, unit commitment~\cite{castillo2016}, and topology reconfiguration~\cite{barrows2014,coffrin2014pscc} problems with AC power flow models as well as a variety of emerging problems related to power systems resilience and restoration~\cite{bhusal2020,austgen2023comparisons,coffrin2018,rhodes2023long}). Similar to the integral relaxations in branch-and-bound solvers for mixed-integer linear programs (MILPs), spatial branch-and-bound algorithms relax and then iteratively tighten nonconvex expressions associated with the power flow equations. Thus, the tightness of power flow relaxations and the quality of their associated objective value bounds are of key importance in such algorithms. The wide range of power system optimization problems formulated as NLPs and MINLPs for which spatial branch-and-bound algorithms are applicable motivates the development of tighter power flow relaxations. Notable recent developments include new commercial and open-source solvers (e.g., a spatial branch-and-bound algorithm in recent versions of Gurobi~\cite{gurobi_tech_talk} as well as  Alpine~\cite{alpine} and Minotaur~\cite{mahajan2021}, among others) along with related applications to various power systems optimization problems, e.g.,~\cite{phan2012,gopalakrishnan2012,kocuk2018,harsha2018pscc,liu2019,mehrtash2022,bynum2019}.

Beyond spatial branch-and-bound, we emphasize that power flow relaxations are also key to algorithms developed for a range of other applications, including solving robust OPF problems~\cite{pscc2018robust,lorca2017robust}, calculating voltage stability margins~\cite{molzahn_lesieutre_demarco-pfcondition,hicss2016}, exploring feasible operating ranges~\cite{molzahn-opf_spaces, NarimaniACC}, designing resilient networks~\cite{nagarajan2017}, assessing severe contingencies~\cite{coffrin2018}, mitigating wildfire ignition risk~\cite{rhodes2023long}, protecting against geomagnetic storms~\cite{lu2018gic}, computing operating envelopes for aggregators of distributed energy resources~\cite{verbic2017irep,hicss2019}, etc. With the need to repeatedly bound the objective values of certain subproblems, only convex relaxations provide the rigorous guarantees needed for many of these algorithms. Stronger relaxations that provide tighter objective value bounds are thus a key enabling methodology for many important applications. For applications where good estimates of the optimal decision variables are also important, we note recent work that enables high-accuracy recovery of AC power flow solutions from the outputs of power flow relaxations~\cite{taheri2024}.}

% Convex relaxation methods are under active research with ongoing efforts targeting to improve their tightness. See~\cite{molzahn2019fnt} for a survey of power flow relaxations.

This paper focuses on improving a particular formulation known as the ``Quadratic Convex'' (QC) relaxation. The QC relaxation encloses the trigonometric, squared, and product terms in a polar representation of power flow equations within convex envelopes~\cite{coffrin2015qc}. Since the quality of these envelopes determines the tightness of the QC relaxation, a number of research efforts have focused on improving these envelopes. These include tighter trigonometric envelopes that leverage sign-definite angle difference bounds~\cite{coffrin2016strengthen_tps,coffrin2016quadtrig};
Lifted Nonlinear Cuts that exploit voltage magnitude and angle difference bounds~\cite{coffrin2016strengthen_tps,chen2015cuts}; cuts based on the voltage magnitude differences between connected  buses~\cite{NarimaniPSCC2018}; tighter envelopes for the product terms~\cite{GlobalSIP2018, harsha2018pscc}; and other valid inequalities, convex envelopes, and cutting planes~\cite{StrongSOCPRelaxations,arctan2}. Most recently, we developed a ``rotated QC'' relaxation~\cite{NarimaniTPS1} which applies a coordinate transformation via a complex per unit base power normalization to tighten envelopes for the trigonometric terms. 

This paper proposes two additional improvements for tightening the QC relaxation. The first improvement considers a particular nonlinear function which has projections that are the nonlinear expressions appearing in a polar representation of the power flow equations. We construct a convex envelope around this nonlinear function that takes the form of a polytope and then use projections of this envelope to obtain convex expressions for the nonlinear terms in the OPF problem. The second improvement uses certain characteristics of the sine and cosine expressions along with the changes in their curvature to tighten the first improvement's convex envelope. 
We also extend our previous work on the coordinate transformation~\ifarxiv\cite{NarimaniTPS1, narimani2020strengthening} 
\else 
\cite{NarimaniTPS1}\fi~via rotating the power flow equations by an angle specific to each bus in order to obtain a tighter envelope. A heuristic approach is proposed for choosing reasonable values for these rotation angles. The proposed relaxation improves the optimality gaps for 68\% of the PGLib-OPF test cases compared to a state-of-the-art QC relaxation~\cite{Sundar}.  

% \mrn{To check proximity to AC feasibility, we conduct an AC power flow (AC-PF) utilizing set points derived from different convex relaxation approaches including the proposed LRQC relaxation and assess the ensuing constraint violations. In addition, to evaluate the proximity to local optimality, we employed a metric defined as the average normalized distance of the variable values obtained using the inexact convex relaxation to those of a locally optimal solution derived from a non-convex solver~\cite{venzke2020inexact}.}

\mrn{%Note that the proposed LRQC relaxation approach is versatile enough to address different types of nonconvex OPF problems, including those for radial and toroidal systems. Additionally, it is adaptable to both single-phase and three-phase configurations, accommodating balance and unbalanced three-phase OPF problems. This adaptability is achieved by preparing envelopes for each phase and incorporating corresponding constraints for accurate modeling.
We develop and demonstrate our proposed relaxation on balanced single-phase equivalent networks without requirements on the network topology (no restriction to radial systems). Such network representations are most appropriate for typical transmission systems. However, we note that the underlying convex envelopes upon which we build our algorithm are suitable for general trigonometric functions without restriction on the range of the input arguments. Thus, extensions to unbalanced three-phase network models are conceptually straightforward by simply constructing envelopes for each phase along with corresponding constraints and variables~\cite{fobes2020}.}

This paper is organized as follows. Sections~\ref{OPF overview} and~\ref{Overview of QC Relaxation} review the OPF formulation and the previous QC relaxation, respectively. Section~\ref{Rotated_OPF} presents a rotated OPF problem and associated QC relaxation with multiple rotation angles (one per bus). 
Section~\ref{Higher_dimansional_space} describes a nonlinear function which has projections that are the nonlinear expressions appearing in the polar representation of the power flow equations. This section also presents a convex envelope that encloses this function. Section~\ref{Envelopes_trig_terms} exploits characteristics of the trigonometric terms to tighten this envelope. Bringing this all together, Section~\ref{TLRQC_relaxation} formulates our proposed tightened QC relaxation.
Section~\ref{Best_Rotation_Angle} presents a method for selecting the rotation angles at each bus to tighten the relaxation. \mrn{Section~\ref{intersection_volume} presents a method for selecting the number of extreme points for the polytopes that formulate our envelopes in order to balance tradeoffs in tractability and tightness of the relaxation.}  Section~\ref{Numerical_results} provides an empirical evaluation, and Section~\ref{conclusion} concludes the paper.
\ifarxiv
\else
An extended version of this paper in~\cite{narimani2023tightening} further describes
and analyzes the envelopes proposed for the trigonometric terms. 
\fi

\section{Overview of the Optimal Power Flow Problem}
\label{OPF overview}
This section formulates the OPF problem using a polar voltage phasor representation. 
The sets of buses, generators, and lines are $\mathcal{N}$, $\mathcal{G}$, and $\mathcal{L}$, respectively. The set $\mathcal{R}$ contains the index of the bus that sets the angle reference.
Let $S_i^d = {P}_i^d + j{Q}_i^d$ and $S_i^g = {P}_i^g + j {Q}_i^g$ represent the complex load demand and generation, respectively, at bus~$i\in\mathcal{N}$, where $j = \sqrt{-1}$. Let $V_i$ and $\theta_i$ represent the voltage magnitude and angle at bus~$i\in\mathcal{N}$. Let $g_{sh,i} + j b_{sh,i}$ denote the shunt admittance at bus~$i\in\mathcal{N}$. For each generator, define a quadratic cost function with coefficients $c_{2,i} \geq 0$, $c_{1,i}$, and $c_{0,i}$. Upper and lower bounds for all variables are indicated by $\left(\overline{\,\cdot\,}\right)$ and $\left(\underline{\,\cdot\,}\right)$, respectively.
For ease of exposition, each line $\left(l,m\right)\in\mathcal{L}$ is modeled as a $\Pi$ circuit with mutual admittance $g_{lm}+j b_{lm}$ and shunt admittance $j b_{c,lm}$. The voltage angle difference between buses $l$ and $m$ for $(l,m)\in\mathcal{L}$ is denoted as $\theta_{lm}=\theta_{l}-\theta_{m}$. 
The complex power flow into each line terminal $(l,m)\in\mathcal{L}$ is denoted by ${P}_{lm}+j{Q}_{lm}$, and the apparent power flow limit is $\overline{S}_{lm}$. 
The OPF problem is
\begin{subequations}
\label{OPF formulation}
\begin{align}
\label{eq:objective}
& \!\!\!\!\!\!\!\min\quad \sum_{{i}\in \mathcal{G}}
c_{2,i}\left(P_i^g\right)^2+c_{1,i}\,P_i^g+c_{0,i}\\ 
&\!\!\!\!\!\!\!\nonumber \text{subject to} \quad \left(\forall i\in\mathcal{N},\; \forall \left(l,m\right) \in\mathcal{L}\right) \\
\label{eq:pf1}
&\!\!\!\!\!\!\! P_i^g-P_i^d = g_{sh,i}\, V_i^2+\sum_{\substack{(l,m)\in \mathcal{L},\\\text{s.t.} \hspace{3pt} l=i}} \!P_{lm}+\!\!\sum_{\substack{(l,m)\in \mathcal{L},\\\text{s.t.} \hspace{3pt} m=i}} \!\!P_{ml}, \\
\label{eq:pf2}
&\!\!\!\!\!\!\! Q_i^g-Q_i^d = -b_{sh,i}\, V_i^2+\!\!\!\!\!\!\sum_{\substack{(l,m)\in \mathcal{L},\\ \text{s.t.} \hspace{3pt} l=i}} \!\!Q_{lm}+\!\!\!\sum_{\substack{(l,m)\in \mathcal{L},\\ \text{s.t.} \hspace{3pt} m=i}} \!\!\!\!Q_{ml},\\
\label{eq:angref}
&\!\!\!\!\!\!\! \theta_r=0, \quad r\in\mathcal{R}, \\
\label{eq:active_gen_limits}
&\!\!\!\!\!\!\! \underline{P}_i^g\leq P_i^g\leq \overline{P}_i^g, \quad %\label{eq:reactive_gen_limits}
\underline{Q}_i^g\leq Q_i^g\leq \overline{Q}_i^g,\\
\label{eq:bus_magnitude_limits}
&\!\!\!\!\!\!\! \underline{V}_i\leq V_{i} \leq \overline{V}_i,\\
\label{eq:bus_angle_limits}
&\!\!\!\!\!\!\!\underline{\theta}_{lm}\leq \theta_{lm}\leq \overline{\theta}_{lm},\\
\label{eq:pik}
&\!\!\!\!\!\!\! P_{lm} \!=\! g_{lm} V_l^2\! -\! g_{lm} V_l V_m\cos\left(\theta_{lm}\right)\! -\! b_{lm} V_l V_m\sin\left(\theta_{lm}\right),\\
\label{eq:qik}
&\!\!\!\!\!\!\! \nonumber Q_{lm} = -\left(b_{lm}+b_{c,lm}/2\right) V_l^2 + b_{lm} V_l V_m\cos\left(\theta_{lm}\right)\\ &\qquad\qquad  - g_{lm} V_l V_m\sin\left(\theta_{lm}\right),\\
\label{eq:pki}
&\!\!\!\!\!\!\!P_{ml}\! =\! g_{lm} V_m^2\! -\! g_{lm} V_l V_m\cos\left(\theta_{lm}\right)\! +\! b_{lm} V_l V_m\sin\left(\theta_{lm}\right),\\
\label{eq:qki}
&\!\!\!\!\!\!\!\nonumber Q_{ml} = -\left(b_{lm}+b_{c,lm}/2\right) V_m^2 + b_{lm} V_l V_m\cos\left(\theta_{lm}\right)\\ &\qquad\qquad  + g_{lm} V_l V_m\sin\left(\theta_{lm}\right),\\
\label{eq: apparent power limit sending}
&\!\!\!\!\!\!\! \left(P_{lm}\right)^2+\left(Q_{lm}\right)^2 \leq \left(\overline{S}_{lm}\right)^2, \quad %\label{eq: apparent power limit receiving}
\left(P_{ml}\right)^2+\left(Q_{ml}\right)^2 \leq \left(\overline{S}_{lm}\right)^2.
\end{align}
\end{subequations}

The objective~\eqref{eq:objective} minimizes the generation cost. Constraints~\eqref{eq:pf1} and~\eqref{eq:pf2} enforce power balance at each bus. Constraint~\eqref{eq:angref} sets the reference bus angle. The constraints in~\eqref{eq:active_gen_limits} bound the active and reactive power generation at each bus. Constraints~\eqref{eq:bus_magnitude_limits}--\eqref{eq:bus_angle_limits}, respectively, bound the voltage magnitudes and voltage angle differences. Constraints~\eqref{eq:pik}--\eqref{eq:qki} relate the active and reactive power flows with the voltage phasors at the terminal buses. The constraints in~\eqref{eq: apparent power limit sending} limit the apparent power flows into both terminals of each line.
\section{Traditional QC Relaxation}
\label{Overview of QC Relaxation}

As typically formulated, the QC relaxation convexifies the OPF problem~\eqref{OPF formulation} by enclosing the nonconvex expressions ($V_i^2$, $\forall i\in\mathcal{N}$, $V_l V_m \cos(\theta_{lm})$ and $V_l V_m \sin(\theta_{lm})$, $\forall (l,m)\in\mathcal{L}$) in convex envelopes~\cite{coffrin2015qc,Sundar}. The envelope for the generic squared function $x^2$ is $\langle x^2\rangle^T$:
\begin{align}
\label{eq:squareenvelopes}
\langle x^2\rangle^T =
\begin{cases}
\widecheck{x}: \begin{cases}\check{x} \geq x^2,\\
\widecheck{x} \leq \left({\overline{x}+\underline{x}}\right) x-{\overline{x}\underline{x}},\\
\end{cases}
\end{cases}
\end{align}
where $\widecheck{x}$ is a lifted variable representing the squared term. 
Envelopes for the generic trigonometric functions $\sin(x)$ and $\cos(x)$ are $\left\langle \sin(x)\right\rangle^S$ and $\left\langle \cos(x)\right\rangle^C$:
\label{eq:convex_envelopes_sin&cos}
\begin{small}
\begin{align}
\label{eq:sine envelope}
 &\left\langle \sin(x)\right\rangle^S\!\!\! =\!
\begin{cases}
\!\widecheck{S}\!:\!\begin{cases}
\widecheck{S}\!\leq\!\cos\left(\frac{x^m}{2}\right)\!\!\left(x\!-\!\frac{x^m}{2}\right)\!\!+\!\sin \left(\frac{x^m}{2}\right) \text{if~} \underline{x}\!\leq0\!\le\overline{x},\\
\widecheck{S}\!\geq\!\cos\left(\frac{x^m}{2}\right)\!\!\left(x\!+\!\frac{x^m}{2}\right)\!\!-\!\sin\left(\frac{x^m}{2}\right)\text{if~} \underline{x}\!\leq0\!\le\overline{x},\\
\widecheck{S}\geq\frac{\sin\left({\underline{x}}\right)-\sin\left(\overline{x}\right)}{{\underline{x}-\overline{x}}}\left(x-{\underline{x}}\right)+\sin\left({\underline{x}}\right) \text{if~} \underline{x}\geq0,\\
\widecheck{S}\leq\frac{\sin\left({\underline{x}}\right)-\sin\left({\overline{x}}\right)}{{\underline{x}-\overline{x}}}\left(x-{\underline{x}}\right)+\sin\left({\underline{x}}\right) \text{if~} {\overline{x}}\leq0,
\end{cases}
\end{cases}\\[-0.3em]
\raisetag{18pt} \label{eq:cosine envelope}
&\left\langle\cos(x)\right\rangle^C =
\begin{cases}
\widecheck{C}:\!\begin{cases}
\widecheck{C}\leq 1-\frac{1-\cos\left({x^m}\right)}{\left(x^m\right)^2}x^2,\\
\widecheck{C}\geq\frac{\cos\left(\underline{x}\right)-\cos\left({\overline{x}}\right)}{{\underline{x}-\overline{x}}}\left(x-{\underline{x}}\right)+\cos\left({\underline{x}}\right),
\end{cases}
\end{cases}
\end{align}
\end{small}%
where $x^m = \max(\left|\underline{x}\right|,\left|\overline{x}\right|)$. The envelopes $\left\langle \sin(x)\right\rangle^S$ and $\left\langle \cos(x)\right\rangle^C$ in~\eqref{eq:sine envelope} and~\eqref{eq:cosine envelope} are valid for $-\frac{\pi}{2} \leq x\leq \frac{\pi}{2}$. 

 The lifted variables $\widecheck{S}$ and $\widecheck{C}$ are associated with the envelopes for the functions $\sin(\theta_{lm})$ and $\cos(\theta_{lm})$. 
The QC relaxation of the OPF problem in~\eqref{OPF formulation} is:% formulated 
\begin{subequations}
\label{eq:QC relaxation}
\begin{align}
& \min \quad \sum_{\footnotesize{{i}\in \mathcal{N}}} c_{2,i}\left(P_i^g\right)^2+c_{1,i}\,P_i^g+c_{0,i}\\ 
&\nonumber \text{subject to} \quad \left(\forall i\in\mathcal{N},\; \forall \left(l,m\right) \in\mathcal{L}\right) \\
\label{eq:qc_p}
& P_i^g-P_i^d = g_{sh,i}\, w_{ii}+\!\!\!\!\!\sum_{\footnotesize{\substack{(l,m)\in \mathcal{L},\\ \text{s.t.} \hspace{3pt} l=i}}} \!\!\!P_{lm}+\sum_{\footnotesize{\substack{(l,m)\in \mathcal{L},\\ \text{s.t.} \hspace{3pt} m=i}}} \!\!P_{ml}, \\
\label{eq:qc_q}
& Q_i^g-Q_i^d = -b_{sh,i}\, w_{ii}+\!\!\!\sum_{\footnotesize{\substack{(l,m)\in \mathcal{L},\\ \text{s.t.} \hspace{3pt} l=i}}} \!\!Q_{lm}+\sum_{\footnotesize{\substack{(l,m)\in \mathcal{L},\\ \text{s.t.} \hspace{3pt} m=i}}} \!\!Q_{ml},\\
\label{eq:qc_V}
&  (\underline{V}_i)^2\leq w_{ii} \leq (\overline{V}_i)^2, \qquad w_{ii} \in\left\langle V_i^2 \right\rangle^T,\\
\label{eq:qc_pik}
& P_{lm} = g_{lm} w_{ll} - g_{lm} c_{lm} - b_{lm} s_{lm},\\
\label{eq:qc_qik}
& Q_{lm} = -\left(b_{lm}+b_{c,lm}/2\right) w_{ll} + b_{lm} c_{lm}- g_{lm} s_{lm}, \\
\label{eq:qc_pki}
& P_{ml} = g_{lm} w_{mm} - g_{lm} c_{lm} + b_{lm} s_{lm}, \\
\label{eq:qc_qki}
& Q_{ml} = -\left(b_{lm}+b_{c,lm}/2\right) w_{mm} + b_{lm} c_{lm}+ g_{lm} s_{lm},\\
\label{eq: apparent power limit sending_QC}
& \left(P_{lm}\right)^2+\left(Q_{lm}\right)^2 \leq \left(\overline{S}_{lm}\right)^2, \quad 
\left(P_{ml}\right)^2+\left(Q_{ml}\right)^2 \leq \left(\overline{S}_{lm}\right)^2,\\
    \label{eq:qc_cik}
    \nonumber & {c}_{lm}  = \sum_{\footnotesize{k=1,\ldots,8}}\!\! \mu_{lm,k}\, \rho^{(k)}_{lm,1} \rho^{(k)}_{lm,2} \rho^{(k)}_{lm,3},  \quad \widecheck{C}_{lm} \in \left\langle\cos(\theta_{lm})\right\rangle^C,\\
    \nonumber & V_l =\!\!\!\! \sum_{\footnotesize{k=1,\ldots,8}}\!\! \mu_{lm,k}\rho^{(k)}_{lm,1},\; V_m =\!\!\!\! \sum_{\footnotesize{k=1,\ldots,8}}\!\! \mu_{lm,k}\rho^{(k)}_{lm,2},\\
    \nonumber&\widecheck{C}_{lm} =\!\!\!\! \sum_{\footnotesize{k=1,\ldots,8}}\!\! \mu_{lm,k}\rho^{(k)}_{lm,3},\\
     & \sum_{k=1,\ldots,8} \mu_{lm,k} = 1, \quad \mu_{lm,k} \geqslant 0, \quad k=1,\ldots,8,\\
        \label{eq:qc_sik}
          \nonumber & {s}_{lm}  = \sum_{\footnotesize{k=1,\ldots,8}}\!\!\! \gamma_{lm,k}\, \zeta^{(k)}_{lm,1} \zeta^{(k)}_{lm,2} \zeta^{(k)}_{lm,3}, \quad \widecheck{S}_{lm} \in \left\langle\sin(\theta_{lm})\right\rangle^S,\\
    \nonumber & V_l =\!\!\!\! \sum_{\footnotesize{k=1,\ldots,8}}\!\! \gamma_{lm,k}\zeta^{(k)}_{lm,1},\; V_m =\!\!\!\! \sum_{\footnotesize{k=1,\ldots,8}}\!\! \gamma_{lm,k}\zeta^{(k)}_{lm,2},\\ \nonumber&\widecheck{S}_{lm} =\!\!\!\! \sum_{\footnotesize{k=1,\ldots,8}}\!\!\ \gamma_{lm,k}\zeta^{(k)}_{lm,3},\\
     & \sum_{\footnotesize{k=1,\ldots,8}} \gamma_{lm,k} = 1, \quad \gamma_{lm,k} \geqslant 0, \quad k=1,\ldots,8,  \\
\label{eq:ell_apparent}
& P_{lm}^2 + Q_{lm}^2 \leq w_{ll}\, \ell_{lm},\\
\label{eq:ell_expression}
&\ell_{lm}\! =\! \left(Y_{lm}^2- \frac{b_{c,lm}^2}{4}\right)\!\ w_{ll} + Y_{lm}^2 w_{mm} - 2Y_{lm}^2c_{lm} - b_{c,lm}Q_{lm},\\
\label{eq:qc_linking}
& \begin{bmatrix}
\mu_{lm,1}+\mu_{lm,2}-\gamma_{lm,1}-\gamma_{lm,2} \\
\mu_{lm,3}+\mu_{lm,4}-\gamma_{lm,3}-\gamma_{lm,4} \\
\mu_{lm,5}+\mu_{lm,6}-\gamma_{lm,5}-\gamma_{lm,6} \\
\mu_{lm,7}+\mu_{lm,8}-\gamma_{lm,7}-\gamma_{lm,8}
\end{bmatrix}^\intercal
\begin{bmatrix}
\underline{V}_l \underline{V}_m\\
\underline{V}_l \overline{V}_m \\
\overline{V}_l \underline{V}_m \\
\overline{V}_l \overline{V}_m
\end{bmatrix}=0,\\
\label{eq:qc_others}
&\text{Equations~}\eqref{eq:angref}\text{--}\eqref{eq:bus_angle_limits},\eqref{eq: apparent power limit sending},
\end{align}
\end{subequations}
where $\ell_{lm}$ represents the squared magnitude of the current flow into terminal $l$ of line $(l,m)\in\mathcal{L}$ and $\left(\cdot\right)^\intercal$ is the transpose operator. The relationship between $\ell_{lm}$ and the power flows $P_{lm}$ and $Q_{lm}$ in~\eqref{eq:ell_apparent} tightens the QC relaxation~\cite{farivar2011,coffrin2015qc}.
Also, as shown in~\eqref{eq:qc_V}, $w_{ii}$ is associated with the squared voltage magnitude at bus~$i$. Note that~\eqref{eq: apparent power limit sending_QC} and~\eqref{eq:ell_apparent} are convex quadratic constraints, while all other constraints are linear.

The lifted variables $c_{lm}$ and $s_{lm}$ represent relaxations of the product terms $V_lV_m\cos(\theta_{lm})$ and $V_lV_m\sin(\theta_{lm})$, respectively, with~\eqref{eq:qc_cik} and~\eqref{eq:qc_sik} formulating an ``extreme point'' representation of the convex hulls for the product terms $V_l V_m \widecheck{C}_{lm}$ and $V_l V_m \widecheck{S}_{lm}$\cite{rikun1997,harsha2018pscc, GlobalSIP2018}.\footnote{An extreme point representation formulates a polytope as a convex combination of its vertices~\cite{rikun1997}.} 
The auxiliary variables $\mu_{lm,k},\; \gamma_{lm,k}\in[0,1]$, $k=1,\ldots,8$, $(l,m)\in\mathcal{L}$, are used in the formulations of these convex hulls.
The extreme points of $V_lV_m\widecheck{C}_{lm}$ are $\rho^{(k)}_{lm} \in [\underline{V_l},\overline{V_l}] \times [\underline{V_m},\overline{V_m}]\times[\underline{\widecheck{C}}_{lm},\overline{\widecheck{C}}_{lm}]$, $k=1,\ldots,8$, and the extreme points of $V_l V_m\widecheck{S}_{lm}$ are $\zeta^{(k)} \in [\underline{V_l},\overline{V_l}] \times [\underline{V_m},\overline{V_m}]\times[\underline{\widecheck{S}}_{lm},\overline{\widecheck{S}}_{lm}]$, $k=1,\ldots,8$. Since sine and cosine are odd and even functions, respectively, $c_{lm} = c_{ml}$ and $s_{lm} = -s_{ml}$.

``Linking constraints'' \eqref{eq:qc_linking} associated with the $V_lV_m$ terms that are shared in $V_lV_m\cos(\theta_{lm})$ and $V_lV_m\sin(\theta_{lm})$ are also enforced to tighten the QC relaxation~\cite{Sundar}.%~\eqref{eq:qc_cik} and~\eqref{eq:qc_sik}.

\section{Exploiting Rotational Degrees of Freedom}
\label{Rotated_OPF}
To provide tighter envelopes for the nonlinear terms in the OPF problem, our previous work in~\cite{NarimaniTPS1} considered a polar representation of the branch admittances, $Y_{lm} \mathrm{e}^{j \delta_{lm}}$, as opposed to the rectangular admittance representation $g_{lm} + j b_{lm}$ used in~\eqref{eq:QC relaxation}. We also used a per unit normalization with a complex base power, i.e., $S_{base}e^{j\psi}$, to improve the QC relaxation’s trigonometric envelopes. The angle of the base power, $\psi$, affects the arguments of the trigonometric functions~\cite{NarimaniTPS1}:
\begin{subequations}
\label{rotated power flow complex}
\begin{align}
\label{rotated power flow sending}
&\nonumber \tilde{S}_{lm}=S_{lm}/e^{j\psi}= \left(Y_{lm}e^{-j(\delta_{lm}+\psi)}+(b_{c,lm}/2)e^{-j(\frac{\pi}{2}+\psi)}\right)V_l^2\\&\qquad\qquad\qquad\qquad - Y_{lm}V_l V_m e^{j(-\delta_{lm}+\theta_{lm}-\psi)},\\
  \label{rotated power flow receiving}
&\nonumber \tilde{S}_{ml}=S_{ml}/e^{j\psi} = \left(Y_{lm}  e^{-j(\delta_{lm}+\psi)}+(b_{c,lm}/2)  e^{-j(\frac{\pi}{2}+\psi)}\right)V_m^2\\&\qquad\qquad\qquad\qquad - Y_{lm}V_m V_l e^{-j(\delta_{lm}+ \theta_{lm}+\psi)}.
\end{align}
\end{subequations}
The angle of the complex base power, $\psi$, linearly enters the arguments of the trigonometric terms, thus providing a rotational degree of freedom in the power flow equations~\cite{NarimaniTPS1}. In~\cite{NarimaniTPS1}, we exploited this rotational degree of freedom to improve the QC relaxation's envelopes. In this section, we extend this prior work by considering multiple rotation angles (one per bus) as opposed to the single rotation angle in~\cite{NarimaniTPS1}. We first describe the new rotated OPF formulation and then formulate its convex relaxation.

\subsection{Rotated OPF Formulation}
Permitting each bus to have a different rotation angle extends our previous work~\cite{NarimaniTPS1}. We define an angle $\psi_l$ for each bus~$l$ via a unit-length complex parameter $e^{j\psi_l}$. To ensure that the power balance constraints at each bus consider quantities that have been rotated consistently, the power flow equations for each line connected to bus~$l$ must use the same angle~$\psi_l$. 
 %In the formulation of power flow equations for each bus $l$, it is crucial to address the assignment of rotation angles to individual lines.
Thus, when formulating the power balance equations for a specific bus, e.g., bus $l$, the power flow equations for every line connected to bus $l$ are rotated by a consistent angle, denoted as $\phi_l$. To achieve this, we form rotated versions of the line flow equations for all lines connected to bus $l$ as follows: %Specifically, considering a line connecting bus $l$ and bus $m$, we introduce the rotated complex power flow variables as follows:
\begin{equation*}
\tilde{S}_{lm} = \frac{S_{lm}}{e^{j\psi_l}}, \quad \tilde{S}_{ml} = \frac{S_{ml}}{e^{j\psi_l}}.
\end{equation*}
Rotated quantities are accented with a tilde, $\left(\tilde{\,\cdot\,}\right)$.
% This notation elucidates the differentiation between the initial and rotated power flows.}
The power generation and load demands are adapted by the rotation angles as formulated in~\eqref{eq:gen transformation} and~\eqref{eq:load transformation}:
\begin{align}
\label{eq:gen transformation}
& \begin{bmatrix}
\tilde{P}_l^g \\
\tilde{Q}_l^g
\end{bmatrix} = \begin{bmatrix}
\cos(\psi_l) & \sin(\psi_l) \\
-\sin(\psi_l) & \cos(\psi_l)
\end{bmatrix}
\begin{bmatrix}
P_l^g \\
Q_l^g
\end{bmatrix},\\
\label{eq:load transformation}
& \begin{bmatrix}
\tilde{P}_l^d \\
\tilde{Q}_l^d
\end{bmatrix} = \begin{bmatrix}
\cos(\psi_l) & \sin(\psi_l) \\
-\sin(\psi_l) & \cos(\psi_l)
\end{bmatrix}
\begin{bmatrix}
P_l^d \\
Q_l^d
\end{bmatrix}.
\end{align}
The rotation angles, $\psi_l$, linearly enter the arguments of the trigonometric terms in the power flow equations in the rotated OPF problem, as shown in~\eqref{rotated power flow}, where $\Re(\,\cdot\,)$ and $\Im(\,\cdot\,)$ are the real and imaginary parts of a quantity:  
\begin{subequations}
\interdisplaylinepenalty=10000
\label{rotated power flow}
\begin{align}
\nonumber & \tilde{P}_{lm}\! =\! \Re(\tilde{S}_{lm})\! =\! \left(Y_{lm}\cos(\delta_{lm}+\psi_l)-(b_{c,lm}/2)\sin(\psi_l)\right)V_l^2\\ \label{eq:shifted OPF active power flow lm} &\qquad\quad\quad\qquad\qquad -
  Y_{lm}V_l V_m\cos(\theta_{lm}-\delta_{lm}-\psi_l),\\
  \label{eq:shifted OPF reactive power flow lm}
\nonumber & \tilde{Q}_{lm} \!=\! \Im(\tilde{S}_{lm}) \!=\! -\left(Y_{lm}\sin(\delta_{lm}+\psi_l)\!+\!(b_{c,lm}/2)\cos(\psi_l)\right)\!V_l^2\\ &\qquad\quad\quad\qquad\qquad -
  Y_{lm}V_l V_m\sin(\theta_{lm}-\delta_{lm}-\psi_l),\\
  \label{eq:shifted OPF active power flow ml}
\nonumber & \tilde{P}_{ml} \!=\! \Re(\tilde{S}_{ml}) \!=\! \left(Y_{lm} \cos(\delta_{lm}+\psi_l)\!-\!(b_{c,lm}/2)  \sin(\psi_l)\right)\!V_m^2
  \\ &\qquad\quad\quad\qquad\qquad - Y_{lm}V_m V_l \cos(\theta_{lm}+\delta_{lm}+\psi_l),\\
\nonumber & \tilde{Q}_{ml} \!=\! \Im(\tilde{S}_{ml}) \!=\! - \left(Y_{lm}  \sin(\delta_{lm}+\psi_l)\!+\!(b_{c,lm}/2) \cos(\psi_l)\right)\!V_m^2 \\ &\qquad\quad\quad\qquad\qquad + Y_{lm}V_m V_l \sin(\theta_{lm}+\delta_{lm}+\psi_l).
  \label{eq:shifted OPF reactive power flow ml}
\end{align}
\end{subequations}

Applying~\eqref{eq:gen transformation}--\eqref{rotated power flow} to~\eqref{OPF formulation} yields a ``rotated'' OPF problem. 
The rotation angles, $\psi_l$, add degrees of freedom to the arguments of the trigonometric terms in~\eqref{rotated power flow}. As we will discuss in Section~\ref{Best_Rotation_Angle}, appropriately chosen values for $\psi_l$ can yield tighter envelopes for these terms. 
\subsection{Rotated QC Relaxation}
\label{Rotated_QC}
Enclosing the product and trigonometric terms in the rotated OPF problem yield a “Rotated QC” (RQC) relaxation:
\begin{subequations}
\label{eq:RQC OPF}
\begin{align}
&\nonumber \min\quad \textstyle\sum_{\footnotesize{{k}\in \mathcal{G}}}
c_{2,k}\left( {\tilde P_k^g\cos (\psi_l ) - \tilde Q_k^g\sin (\psi_l)} \label{eq:RQC obj} \right)^2\\[-0.5em]
&\quad\qquad\qquad+c_{1,k}\left( {\tilde P_k^g\cos (\psi_l ) - \tilde Q_k^g\sin (\psi_l )} \right)+c_{0,k}\\
&\nonumber \text{subject to} \quad \left(\forall i\in\mathcal{N}, \forall   \left(l,m\right) \in\mathcal{L}\right)\\
& \nonumber \tilde{P}_i^g-\tilde{P}_i^d =\left(g_{sh,i}  \cos(\psi_l)-b_{sh,i} \sin(\psi_l)\right)w_{ii}\\ \label{eq:RQC active power injection}&\quad\quad\quad\quad +\sum_{\footnotesize{\substack{(l,m)\in \mathcal{L},\\ \text{s.t.} \hspace{3pt} l=i}}} \tilde{P}_{lm}+\sum_{\footnotesize{\substack{(l,m)\in \mathcal{L},\\ \text{s.t.} \hspace{3pt} m=i}}} \tilde{P}_{ml},\\
& \nonumber \tilde{Q}_i^g-\tilde{Q}_i^d =-\left(g_{sh,i}\sin(\psi_l)+b_{sh,i}  \cos(\psi_l)\right)w_{ii}\\ \label{eq:RQC reactive power injection} &\quad\quad\quad\quad +\sum_{\footnotesize{\substack{(l,m)\in \mathcal{L},\\ \text{s.t.} \hspace{3pt} l=i}}} \tilde{Q}_{lm}+\sum_{\footnotesize{\substack{(l,m)\in \mathcal{L},\\ \text{s.t.} \hspace{3pt} m=i}}} \tilde{Q}_{ml},\\
\label{eq:rqc_V}
&  (\underline{V}_i)^2\leq w_{ii} \leq (\overline{V}_i)^2, \qquad w_{ii} \in\left\langle V_i^2 \right\rangle^T,\\
\label{eq:shifted OPF slack bus}
&  \theta_{ref}=0,\\
\label{eq:shifted OPF active power limits}
& \underline{P}_i^{g}\leq \tilde P_i^g\cos (\psi_l ) - \tilde Q_i^g\sin (\psi_l )\leq \overline{P}_i^{g},\\
 \label{eq:shifted OPF reactive power limits}
& \underline{Q}_i^{g}\leq \tilde Q_i^g\cos (\psi_l ) + \tilde P_i^g\sin (\psi_l )\leq \overline{Q}_i^{g},\\
 \label{eq:shifted OPF voltage limits}
& \underline{V}_i\leq V_{i} \leq \overline{V}_i, \quad\qquad
\underline{\theta}_{lm}\leq \theta_{lm}\leq \overline{\theta}_{lm},\\
  \label{eq:shifted OPF power flow limit lm}
& (\tilde{P}_{lm})^2+(\tilde{Q}_{lm})^2 \leq (\overline{S}_{lm})^2,\quad
 (\tilde{P}_{ml})^2+(\tilde{Q}_{ml})^2 \leq (\overline{S}_{lm})^2,\\
\label{eq:RQC active power flow lm}
& \nonumber \tilde{P}_{lm} \!=\! \left(Y_{lm}\cos(\delta_{lm} \!+\! \psi_l) - b_{c,lm}/2 \sin(\psi_l)\right)w_{ll}\\  &\qquad - Y_{lm}\tilde{c}_{lm},\\[-0.35em]
  \label{eq:RQC reactive power flow lm}
& \nonumber \tilde{Q}_{lm}\!=\! -\left(Y_{lm}\sin(\delta_{lm}\!+\!\psi_l)+b_{c,lm}/2\cos(\psi_l)\right)w_{ll}\\  
  &\qquad -Y_{lm}\tilde{s}_{lm},\\[-0.35em]
  \label{eq:RQC active power flow ml}
& \tilde{P}_{ml}\!=\!
  -Y_{lm}\tilde{c}_{lm} \!+\! \left( Y_{lm} \cos(\delta_{lm}\!+\! \psi_l) - b_{c,lm}/2 \sin(\psi_l)\right)w_{mm},\\[-0.35em]
  \label{eq:RQC reactive power flow ml}
& \tilde{Q}_{ml} \!=\!
  Y_{lm}\tilde{s}_{lm} - \left(Y_{lm}  \sin(\delta_{lm}\!+\!\psi_l) + b_{c,lm}/2 \cos(\psi_l)\right)w_{mm},\\
 \label{eq:RQC squared current}
& \tilde{P}_{lm}^2 + \tilde{Q}_{lm}^2 \leq w_{ll}\, \tilde{\ell}_{lm},\\
 \label{eq:RQC squared current expression}
&\nonumber \tilde{\ell}_{lm} =\bigg(b_{c,lm}^2/4+Y_{lm}^2- Y_{lm}b_{c,lm}\cos(\delta_{lm}+\psi_l)\sin(\psi_l)\\[-0.35em]
&\nonumber\quad\quad+Y_{lm}b_{c,lm}\sin(\delta_{lm}+\psi_l)\cos(\psi_l)\bigg)V_l^2+Y_{lm}^2 V_m^2\\&\quad\quad\nonumber+\left(-2Y_{lm}^2\cos(\delta_{lm}+\psi_l)+Y_{lm}b_{c,lm}\sin(\psi_l)\right)\tilde{c}_{lm}\\&\quad\quad
+\left(2Y_{lm}^2\sin(\delta_{lm}+\psi_l)+Y_{lm}b_{c,lm}\cos(\psi_l)\right)\tilde{s}_{lm},\\
\nonumber & \tilde{c}_{lm} = \!\!\!\! \sum_{k=1,\ldots,4\tilde{N}}\!\! \lambda_k\, \eta^{(k)}_1 \eta^{(k)}_2 \eta^{(k)}_3\!\!, \;\;\;\; \tilde{s}_{lm} = \!\!\!\! \sum_{k=1,\ldots,4\tilde{N}}\!\!\! \lambda_k\, \eta^{(k)}_1 \eta^{(k)}_2 \eta^{(k)}_4\!\!\!\!, \\
    \nonumber & V_l =\!\!\!\!\! \sum_{k=1,\ldots,4\tilde{N}} \!\!\!\!\!\lambda_k\eta^{(k)}_1,\;\;
    V_m =\!\!\!\!\! \sum_{k=1,\ldots,4\tilde{N}}\!\!\!\!\! \lambda_k\eta^{(k)}_2\!\!,\;\; 
    \tilde{S}_{lm} =\!\!\!\!\! \sum_{k=1,\ldots,4\tilde{N}}\!\!\!\!\! \lambda_k\eta^{(k)}_4\!\! ,\\
    \nonumber & \tilde{C}_{lm} =\!\!\!\!\!\!\!\!\!\! \sum_{k=1,\ldots,4\tilde{N}} \!\!\!\!\!\!\!\lambda_k\eta^{(k)}_3\, \!\!\!,
     \!\!\! \sum_{k=1,\ldots,4\tilde{N}} \!\!\!\!\lambda_k = 1, \qquad \!\!\!\!\!\!\lambda_k \geqslant 0, \quad \!\!\!k=1,\ldots,4\tilde{N},\\
     \label{eq:convex hull trilinear terms} & \tilde{C}_{lm} \in \left\langle\cos(\theta_{lm}\!-\!\delta_{lm}\!-\!\psi_l)\right\rangle^{C}\!\!, \; \tilde{S}_{lm}\in\left\langle\sin(\theta_{lm}\!-\!\delta_{lm}\!-\!\psi_l)\right\rangle^{S}\!\!\\ &\text{Equation~}\text{\eqref{eq:qc_linking}}.
\end{align}
\end{subequations}

\section{Convexifications and Projections of an Alternative Nonlinear Function}
\label{Higher_dimansional_space}
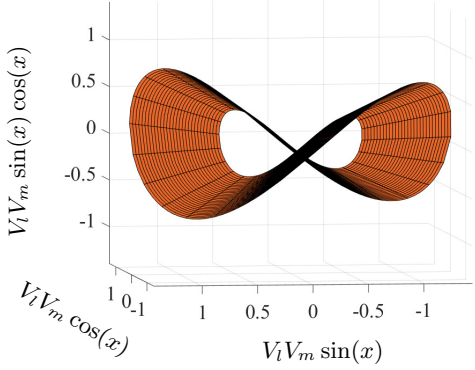
\begin{figure}[tb]
 %\vspace{-0.5cm}
    \centering
\begin{tikzpicture}
\node at (0,0) {\includegraphics[scale=0.42]{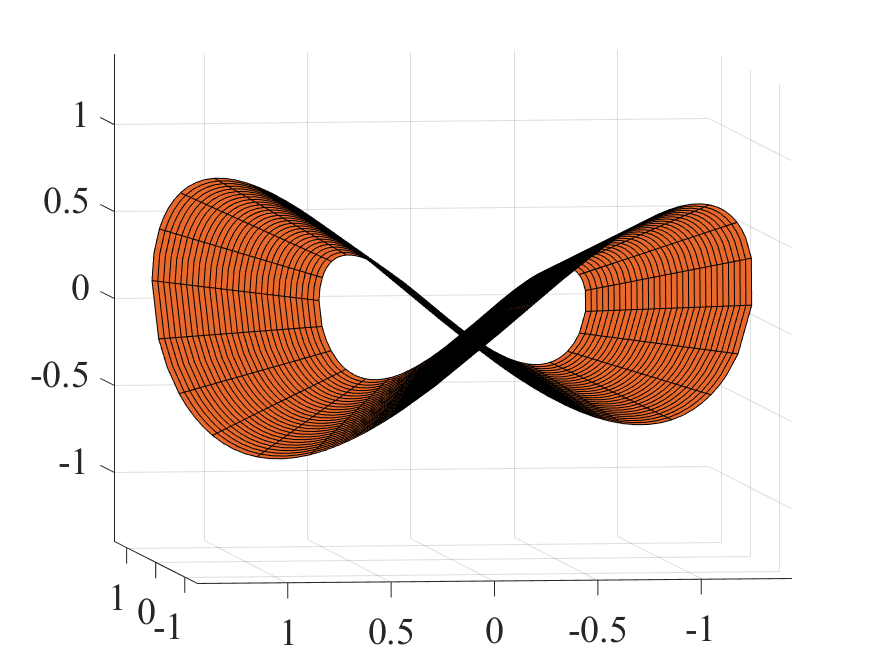}};
\node [anchor=west, rotate=-30, font=\small] at (-3.6,-1.8) {$V_lV_m\cos(x)$};
\node [anchor=west, rotate=0, font=\small] at (-.4,-2.7) { $V_lV_m\sin(x)$};
\node [anchor=west, rotate=90, font=\small] at (-3.5,-1.3) { $V_lV_m\sin(x)\cos(x)$ };
\end{tikzpicture}%
	\caption{A projection of $V_lV_m\cos(x)\sin(x)$. The argument $x$ corresponds to $\theta_{lm} - \delta_{lm} - \psi_l$.
 }
	\label{fig_manifold}
	    %\vspace{-0.3cm}
\end{figure}

In this section, we propose tighter convex envelopes for the nonlinear terms in the power flow equations. Most previous research convexifies the $V_lV_m\cos(\theta_{lm})$ and $V_lV_m\sin(\theta_{lm})$ terms in~\eqref{eq:pik}--\eqref{eq:qki} independently. Instead of individually convexifying these terms, we focus on a different function, $V_lV_m\cos(\theta_{lm} - \delta_{lm} - \psi_l)\sin(\theta_{lm} - \delta_{lm} - \psi_l)$. Fig.~\ref{fig_manifold} shows a projection of this function. 
As we will describe, projections of a convexified form of this function provide tighter envelopes for the product terms in the rotated power flow equations~\eqref{rotated power flow}. 
We first summarize prior QC formulations for comparison purposes and then discuss our proposed formulation. 

% \vspace{-0.5cm}

\begin{figure*}[t]
\centering
\def\tabularxcolumn#1{m{#1}}
\begin{tabularx}{\linewidth}{@{}cXX@{}}
\begin{tabular}{ccc}
\subfloat[The original QC formulation from~\cite{coffrin2015qc}.]{\label{Fig1.a}
\centering
\begin{tikzpicture} \hspace{-.4cm}
    \node[] (image) at (0,0) {\includegraphics[width=0.47\columnwidth, trim={0.5cm 0.3cm 1.3cm 3.5},clip]{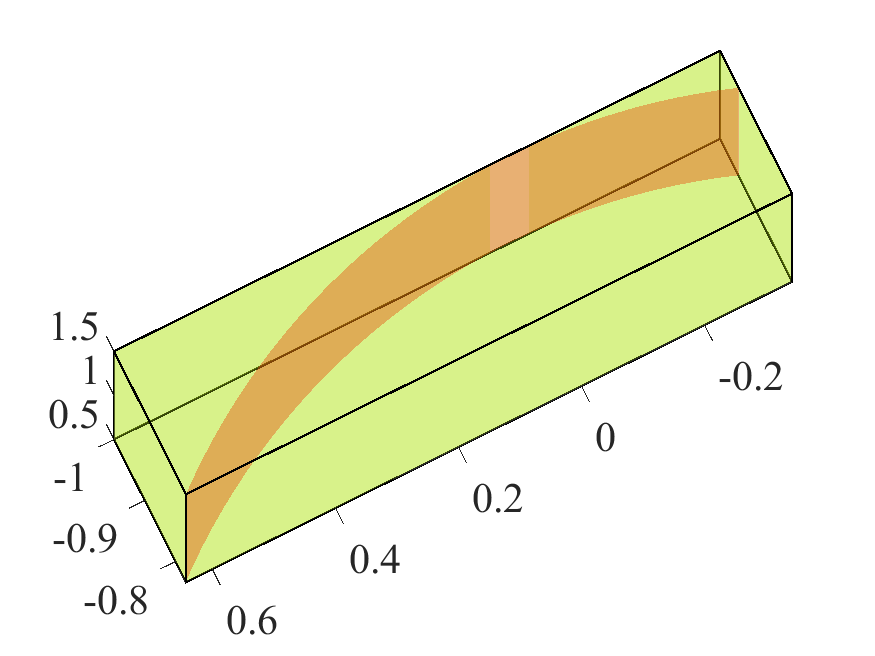}};
    %\begin{scope}[]
 \node [anchor=west, font=\normalsize] at (-0.2,-1.6) {$\sin(x)$};
       \node [anchor=west, font=\normalsize] at (-3.5,-1.4) {$\cos(x)$};
       \node [anchor=west, font=\normalsize] at (-2.8,-0.4) {$V_{l}$};
\end{tikzpicture}
}
\subfloat[The rotated QC relaxation from~\cite{NarimaniTPS1}.]{\label{Fig1.b}
\centering
\centering
\begin{tikzpicture} \hspace{-.4cm}
    \node[] (image) at (0,0) {\includegraphics[width=0.47\columnwidth, trim={0.5cm 0.3cm 1.5cm 3.5},clip]{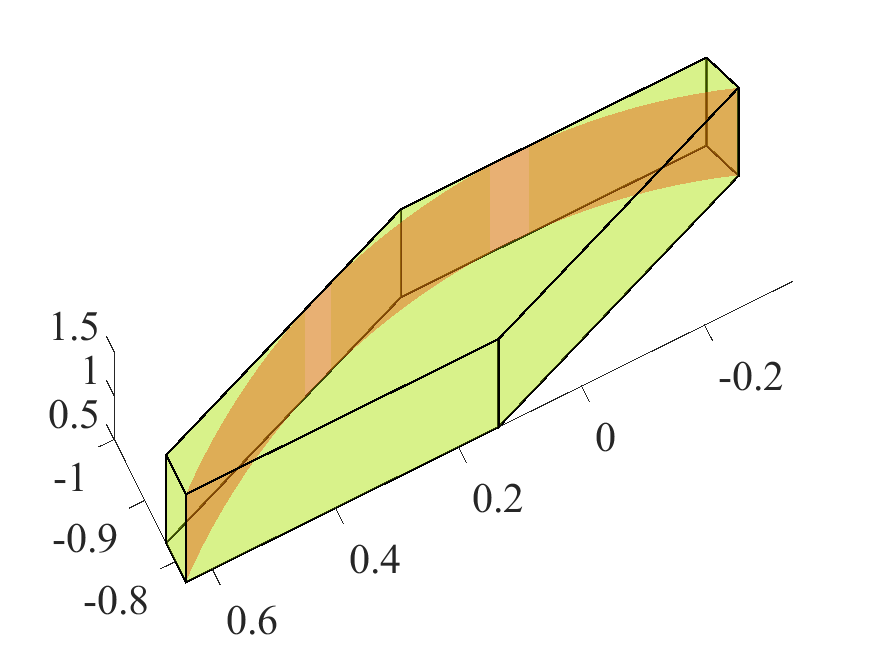}};
 \node [anchor=west, font=\normalsize] at (-0.2,-1.6) {$\sin(x)$};
       \node [anchor=west, font=\normalsize] at (-3.5,-1.4) {$\cos(x)$};
       \node [anchor=west, font=\normalsize] at (-2.8,-0.4) {$V_{l}$};
\end{tikzpicture}
}
\subfloat[The proposed QC relaxation from Section~\ref{subsec:proposed_envelopes}. The envelope (light green region) appears nearly coincident with the function itself (orange region).]{\label{Fig1.c}
\centering
\centering
\begin{tikzpicture} \hspace{-.4cm}
    \node[] (image) at (0,0) {\includegraphics[width=0.47\columnwidth, trim={0.5cm 0.3cm 1.5cm 3.5},clip]{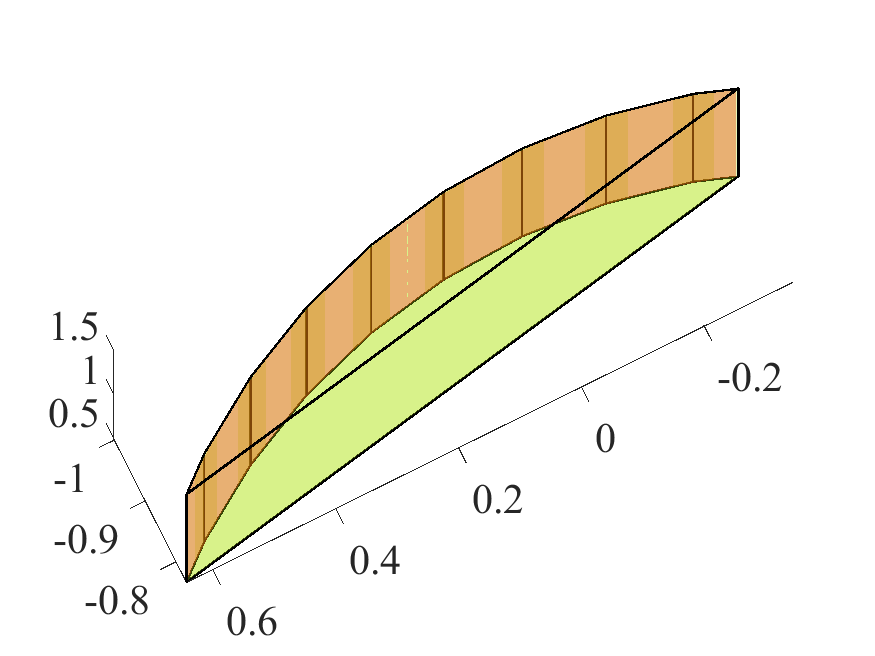}};
 \node [anchor=west, font=\normalsize] at (-0.2,-1.6) {$\sin(x)$};
       \node [anchor=west, font=\normalsize] at (-3.5,-1.4) {$\cos(x)$};
       \node [anchor=west, font=\normalsize] at (-2.8,-0.4) {$V_{l}$};
\end{tikzpicture}
}
 \end{tabular}
\end{tabularx}
\caption{Projections of various envelopes for the function $V_lV_m\cos(x)\sin(x)$ in terms of $V_l$, $\cos(x)$, and $\sin(x)$. The argument $x$ indicates the angle difference $\theta_{lm}$ for the original QC relaxation in Fig.~\ref{Fig1.a} and the rotated argument from the polar admittance representation, $\theta_{lm}-\delta_{lm}-\psi_l$, for the rotated QC relaxations in Figs.~\ref{Fig1.b} and~\ref{Fig1.c}. The pink region common to Figs.~\ref{Fig1.a}--\ref{Fig1.c} is the function $V_lV_m\cos(x)\sin(x)$ that we seek to enclose in a convex envelope. The light green regions in Figs.~\ref{Fig1.a}, \ref{Fig1.b}, and~\ref{Fig1.c} are the surfaces of the convex envelopes proposed in the original QC relaxation~\cite{coffrin2015qc}, the rotated QC relaxation from~\cite{NarimaniTPS1}, and our proposed formulation from Section~\ref{subsec:proposed_envelopes}, respectively.}
\end{figure*}
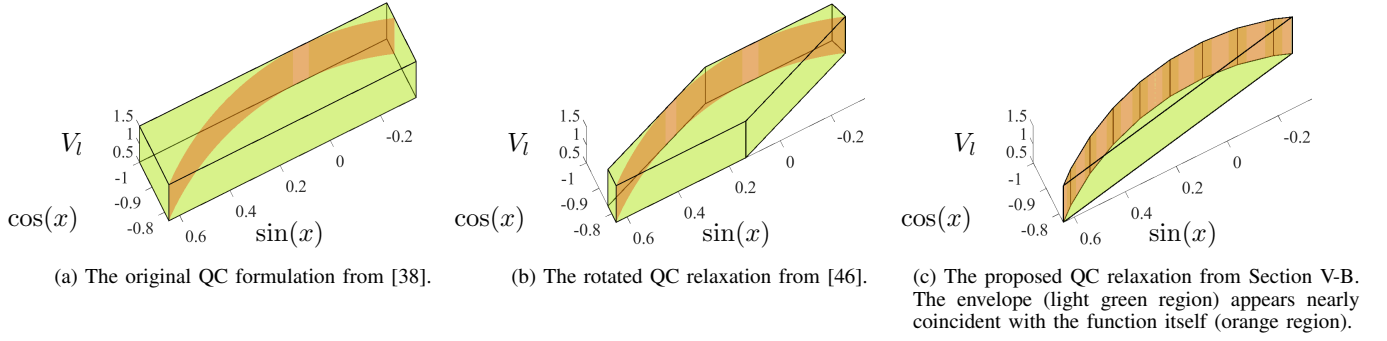

\begin{figure}[t]
\vspace{-0.6cm}
    \centering
\begin{tikzpicture}
\node at (0,0) {\includegraphics[scale=0.45]{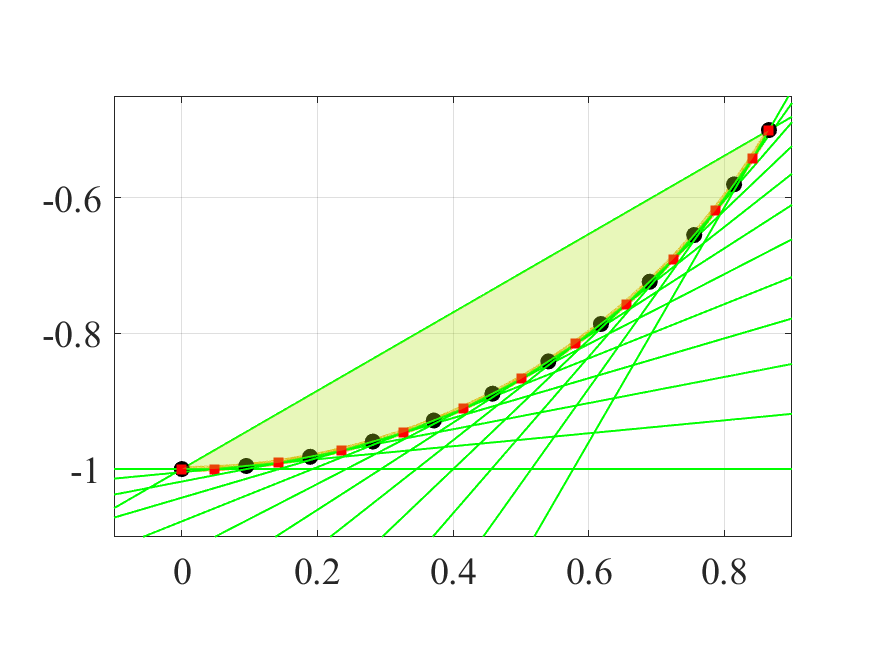}};
\node [anchor=west, font=\normalsize] at (-4.0,0.0) {\begin{tabular}{c} \begin{turn}{90} $\sin(\theta_{lm}-\delta_{lm}-\psi_l)$ \end{turn}\end{tabular}};
\node [anchor=west, font=\normalsize] at (-1.8,-2.3) {\begin{tabular}{c} \begin{turn}{0} $\cos(\theta_{lm}-\delta_{lm}-\psi_l)$ \end{turn}\end{tabular}};
\end{tikzpicture}%
	\caption{Projection of the function $V_lV_m\cos(\theta_{lm}-\delta_{lm}-\psi_l)\sin(\theta_{lm}-\delta_{lm}-\psi_l)$ in terms of $\cos(\theta_{lm}-\delta_{lm}-\psi_l)$ and $\sin(\theta_{lm}-\delta_{lm}-\psi_l)$. The black line is the function $\cos(\theta_{lm}-\delta_{lm}-\psi_l)\sin(\theta_{lm}-\delta_{lm}-\psi_l)$ that we seek to enclose in a convex envelope. The light green lines are tangent to the function at the equally spaced black points. The convex region enclosed by these lines is depicted in light green, encompassing the nonconvex trigonometric function shown in orange. The extreme points of the convex envelope for this function are shown by the red squares and are at the intersections of the green lines. See Algorithm~\ref{alg:algorithm} for details.}
	\label{fig:sin_cos}
\end{figure}

\subsection{Previous Envelopes for Product Terms}
 
By independently convexifying the terms in the products $V_l V_m \cos(\theta_{lm})$ and $V_l V_m \sin(\theta_{lm})$, the original QC relaxation proposed in~\cite{coffrin2015qc} effectively encloses these terms in a rectangle defined by the bounds on $\cos(\theta_{lm})$, $\sin(\theta_{lm})$, $V_l$, and $V_m$. Fig.~\ref{Fig1.a} shows a projection of this envelope.

The approach in~\cite{NarimaniTPS1} also uses the bounds on $\cos(\theta_{lm}-\delta_{lm}-\psi_l)$, $\sin(\theta_{lm}-\delta_{lm}-\psi_l)$, $V_l$, and $V_m$ to create a rectangle enclosing the expressions in~\eqref{eq:shifted OPF active power flow lm} and~\eqref{eq:shifted OPF reactive power flow lm}. Another rectangle is similarly constructed using the bounds on $\cos(\theta_{lm}+\delta_{lm}+\psi_l)$, $\sin(\theta_{lm}+\delta_{lm}+\psi_l)$, $V_l$, and $V_m$ in~\eqref{eq:shifted OPF active power flow ml} and~\eqref{eq:shifted OPF reactive power flow ml}. Considering the intersection of these rectangles yields a convex envelope in the form of a polytope. As shown in Fig.~\ref{Fig1.b}, the envelopes from the rotated QC relaxation~\cite{NarimaniTPS1} can be tighter than those from the original QC relaxation~\cite{coffrin2015qc}.

\subsection{Proposed Envelopes for Product Terms} \label{subsec:proposed_envelopes}

This paper tightens the QC relaxation by constructing an envelope tailored to the function $V_l V_m \cos(\theta_{lm}-\delta_{lm}-\psi_l) \sin(\theta_{lm}-\delta_{lm}-\psi_l)$. To accomplish this, we consider the projection of this function in terms of $\cos(\theta_{lm}-\delta_{lm}-\psi_l)$ and $\sin(\theta_{lm}-\delta_{lm}-\psi_l)$, as shown by the solid black line in Fig.~\ref{fig:sin_cos}. In this projection, the function is an arc of the unit circle defined using the angle difference bounds $\underline{\theta}_{lm}$ and $\overline{\theta}_{lm}$, i.e., $\underline{\theta}_{lm}-\delta_{lm}-\psi_l$ and $\overline{\theta}_{lm}-\delta_{lm}-\psi_l$. 

To construct the convex envelope in Fig.~\ref{fig:sin_cos}, we first compute the green lines that are tangent at the equally spaced black dots. The extreme points defining the polytope that forms the convex envelope are then obtained from the intersections of neighboring tangent lines, denoted by the red squares in Fig.~\ref{fig:sin_cos}. Finally, the polytope is extended using the bounds on $V_l$ and $V_m$ in the same manner as in both the original and rotated QC relaxations~\cite{coffrin2015qc,NarimaniTPS1}. 

Formally, let $\mathcal{T}_{lm} = \{(C_{lm}^{int,1},S_{lm}^{int,1}),\, (C_{lm}^{int,2},S_{lm}^{int,2}),\,\ldots,\, \\(C_{lm}^{int,N_{seg}},S_{lm}^{int,N_{seg}})\}$ denote the coordinates of the extreme points (red squares) in Fig.~\ref{fig:sin_cos}, where $N_{seg}$ is a user-selected parameter for the number of extreme points. Extend these extreme points using the bounds on the voltage magnitudes to obtain the extreme points for a convex envelope enclosing the function $V_l V_m \cos(\theta_{lm}-\delta_{lm}-\psi_l) \sin(\theta_{lm}-\delta_{lm}-\psi_l)$, denoted as $\eta^{(k)}_{lm} \in [\underline{V_l},\overline{V_l}] \times [\underline{V_m},\overline{V_m}]\times\mathcal{T}_{lm}$, $k=1,\ldots,4N_{seg}$.
Algorithm~\ref{alg:algorithm} describes how to compute these extreme points. 
{\SetAlgoNoLine
\begin{algorithm}[t] % check
    \SetAlgoNoLine
    \small
    \caption{Compute Extreme Points}
    \label{alg:algorithm}
    % config
    \SetKwInOut{Input}{Input}
    \SetKwInOut{Output}{Output}
    \SetKwFunction{FMain}{Main}
    \SetKwFunction{FGen}{Generate}
    \DontPrintSemicolon
        \Function{Extreme\_Point($\overline{\theta}_{lm}$,$\underline{\theta}_{lm}$,$\delta_{lm}$,$\psi_l$,$N_{seg}$)}{
    $U \gets \overline{\theta}_{lm}-\delta_{lm}-\psi_l$, $L \gets \underline{\theta}_{lm}-\delta_{lm}-\psi_l$.\;
Divide the arc between $U$ and $L$ into $N_{seg}$ equal segments.\newline
\For{$i = 1,\ldots, N_{seg}$}{
$E=\left[(\tilde{C}_{lm,1},\tilde{S}_{lm,1}),\ldots,(\tilde{C}_{lm,N_{seg}},\tilde{S}_{lm,N_{seg}}) \right]$ $\gets$ Intersection of the tangent lines corresponding to both ends of the $i$-th segment.} 
$EXT=\left[E,(\tilde{C}_{lm,L},\tilde{S}_{lm,L}),(\tilde{C}_{lm,U},\tilde{S}_{lm,U})\right]$ $\gets$ Add the points on the closest tangent line to the arc at the endpoints $U$ and $L$ to $E$; see Fig.~\ref{fig:sin_cos}.

Extend the resulting points by the upper and lower bounds on voltage magnitudes, $[EXT]\times [\underline{V_l},\overline{V_l}] \times [\underline{V_m},\overline{V_m}]$. 

            \Return{$[EXT]\times [\underline{V_l},\overline{V_l}] \times [\underline{V_m},\overline{V_m}]$}.}
\end{algorithm}
}

By introducing auxiliary variables denoted as $\lambda_{lm,k}\in[0,1]$, $k=1,\ldots,4N_{seg}$, we next form a convex envelope for the function $V_l V_m \cos(\theta_{lm}-\delta_{lm}-\psi_l) \sin(\theta_{lm}-\delta_{lm}-\psi_l)$ as the convex combination of the extreme points $\eta^{(k)}$. Finally, we take projections of this convex envelope to obtain envelopes enclosing the products $V_lV_m\cos(\theta_{lm}-\delta_{lm}-\psi_l)$ and $V_lV_m\sin(\theta_{lm}-\delta_{lm}-\psi_l)$. 

Using this procedure, we obtain the following constraints that link the lifted variables $c_{lm}$ and $s_{lm}$ corresponding to the expressions $V_lV_m\cos(\theta_{lm}-\delta_{lm}-\psi_l)$ and $V_lV_m\sin(\theta_{lm}-\delta_{lm}-\psi_l)$ with the remainder of the variables in the problem (i.e., the lifted variables $\widecheck{C}_{lm}$ and $\widecheck{S}_{lm}$ for the cosine and sine terms, $\widecheck{C}_{lm} \in \left\langle\cos(\theta_{lm}\!-\!\delta_{lm}\!-\!\psi_l)\right\rangle^{C}\!\!$ and $\widecheck{S}_{lm} \in\left\langle\sin(\theta_{lm}\!-\!\delta_{lm}\!-\!\psi_l)\right\rangle^{S}$, and the variables $\theta_{lm}$, $V_l$, and $V_m$):

%\vspace{-.2cm}
\begin{align}
\label{eq:quadrilinear_term_envelopes1}
    \nonumber & \tilde{c}_{lm} = \!\!\!\!\!\!\!\! \sum_{k=1,\ldots,4N_{seg}}\!\!\!\! \lambda_{lm,k}\, \eta^{(k)}_{lm,1} \eta^{(k)}_{lm,2} \eta^{(k)}_{lm,4}\!\!, \;\;\;\;V_l =\!\!\!\!\!\!\! \sum_{k=1,\ldots,4N_{seg}} \!\!\!\!\!\!\!\lambda_{lm,k}\eta^{(k)}_{lm,1},\;\;  \\
    \nonumber & \tilde{s}_{lm} = \!\!\!\!\!\!\!\! \sum_{k=1,\ldots,4N_{seg}}\!\!\!\!\!\!\! \lambda_{lm,k}\, \eta^{(k)}_{lm,1} \eta^{(k)}_{lm,2} \eta^{(k)}_{lm,5},~
    V_m =\!\!\!\!\!\!\!\!\! \sum_{k=1,\ldots,4N_{seg}}\!\!\!\!\!\!\!\!\! \lambda_{lm,k}\eta^{(k)}_{lm,2},\\
    \nonumber &\tilde{C}_{lm} =\!\!\!\!\!\!\!\!\!\! \sum_{k=1,\ldots,4N_{seg}} \!\!\!\!\!\!\!\lambda_{lm,k}\eta^{(k)}_{lm,4},
    \quad\tilde{S}_{lm} =\!\!\!\!\!\!\!\! \sum_{k=1,\ldots,4N_{seg}}\!\!\!\!\! \lambda_{lm,k}\eta^{(k)}_{lm,5}\,
    ,\;\;\\ 
    \nonumber &\theta_{lm} =\!\!\!\!\!\!\!\!\!\! \sum_{k=1,\ldots,4N_{seg}} \!\!\!\!\!\!\!\lambda_{lm,k}\eta^{(k)}_{lm,3},\\
  &\sum_{k=1,\ldots,4N_{seg}} \!\!\!\!\lambda_{lm,k} = 1, \qquad \!\!\!\!\!\!\lambda_{lm,k} \geqslant 0, \quad \!\!\!k=1,\ldots,4N_{seg}.
\end{align}

Fig.~\ref{Fig1.c} visualizes a projection of the convex envelope obtained using this approach. Comparing Fig.~\ref{Fig1.c} with Figs.~\ref{Fig1.a} and~\ref{Fig1.b} demonstrates the superiority of the proposed approach in providing tighter envelopes compared to those in~\cite{coffrin2015qc, NarimaniTPS1}. 
Note that~\eqref{eq:quadrilinear_term_envelopes1} precludes the need for the linking constraint~\eqref{eq:qc_linking} that relates the common term $V_l V_m$ in the products $V_l V_m \sin(\theta_{lm}-\delta_{lm}-\psi_l)$ and $V_l V_m \cos(\theta_{lm}-\delta_{lm}-\psi_l)$.

\section{Tighter Trigonometric Envelopes}
\label{Envelopes_trig_terms}

Having addressed the product terms, we next turn our attention to the trigonometric functions $\cos(\theta_{lm}-\delta_{lm}-\psi_l)$ and $\sin(\theta_{lm}-\delta_{lm}-\psi_l)$. This section leverages certain characteristics of the sine and cosine functions along with the changes in their curvature to provide tighter convex envelopes derived using multiple tangent lines to these functions. Figs.~\ref{fig:trig_function.a} and~\ref{fig:trig_function.b} illustrate these tangent lines for the sine and cosine functions, respectively. The remainder of this section focuses solely on the cosine envelopes since the sine envelopes can be constructed as rotated versions of the cosine envelopes. In this section, we present an overview of the key ideas without delving into extensive mathematical details. The complete mathematical derivations related to the concepts discussed in this section can be found in~\ifarxiv the appendix\else \cite[Appendix]{narimani2023tightening}\fi.

\begin{figure}[t]
\centering
\subfloat[Envelope for $\sin(\theta_{lm}-\delta_{lm}-\psi_l), \theta_{lm} \in {[}{-}75^\circ, 75^\circ {]}$ \text{and} $\psi_l+\delta_{lm}= 5^\circ$]{\label{fig:trig_function.a}
\begin{tikzpicture}
    \node[] (image) at (0,0) {\includegraphics[width=0.6\columnwidth]{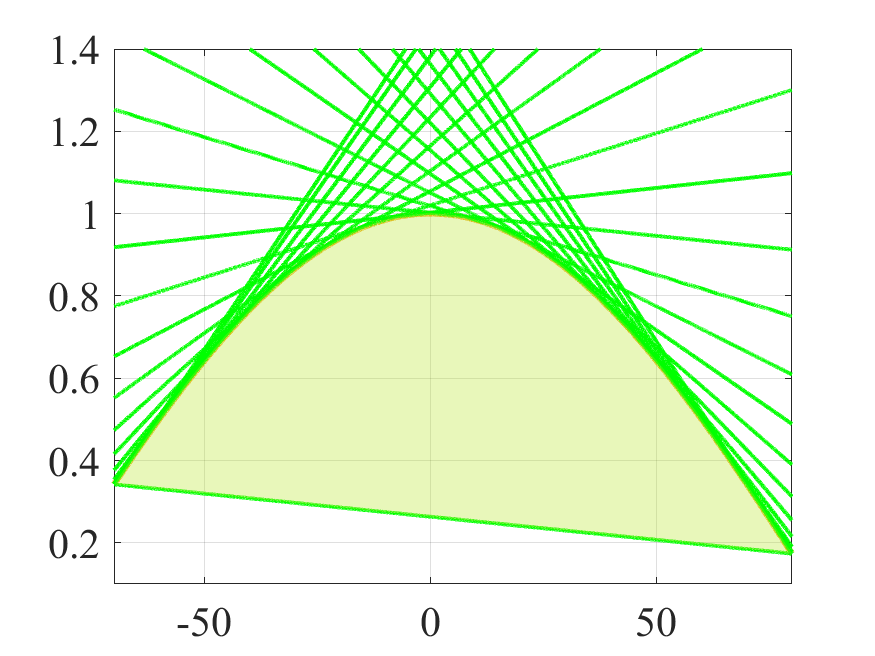}};
       \node [anchor=west, font=\normalsize \rotatebox{90}{$\sin(\theta_{lm}-\delta_{lm}-\psi_l)$}] at (-3.3,0.2) {};
       \node [anchor=west, font=\normalsize] at (-1.3,-2.1) {$\theta_{lm}-\delta_{lm}-\psi_l$};
\end{tikzpicture}
}\\%\vspace{-.3cm} 
\subfloat[Envelope for $\cos(\theta_{lm}-\delta_{lm}-\psi_l)$, $\theta_{lm} \in {[}{-}75^\circ, 75^\circ{]}$ and $\psi_l+\delta_{lm} = 5^\circ$]{\label{fig:trig_function.b}
\centering
\centering
\begin{tikzpicture}
    \node[] (image) at (0,0) {\includegraphics[width=0.6\columnwidth]{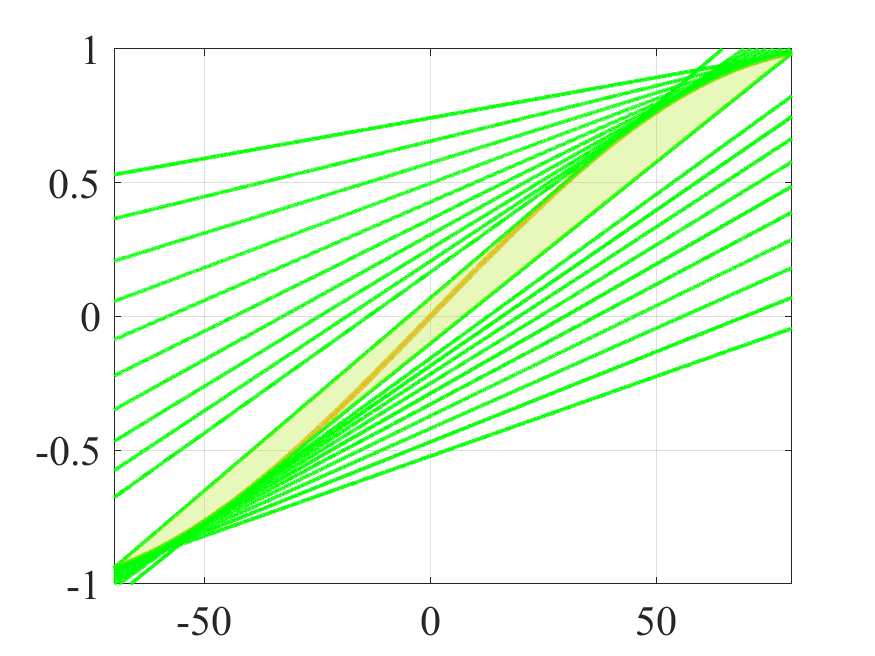}};
       \node [anchor=west, font=\normalsize \rotatebox{90}{$\cos(\theta_{lm}-\delta_{lm}-\psi_l)$}] at (-3.3,0.2) {};
       \node [anchor=west, font=\normalsize] at (-1.3,-2.1) {$\theta_{lm}-\delta_{lm}-\psi_l$};
\end{tikzpicture}
}
\caption{Convex regions (light green) constructed using tangents (solid green lines) to the sine and cosine functions,  $\sin(\theta_{lm}-\delta_{lm}-\psi_l)$ and $\cos(\theta_{lm}-\delta_{lm}-\psi_l)$, as described by \ifarxiv Algorithm~
\ref{alg:algorithm1}
% \ref{eq:roots} 
in the appendix.
\else Algorithm~2 in~\cite[Appendix]{narimani2023tightening}.
\fi}
\end{figure}

We note that the method proposed in this section is a specific form of an approach recently developed in~\cite{sundar2022sequence} that uses a sequence of linear programming relaxations which converge towards the convex hull of a univariate function. Our proposed method is an explicit form for a sequence of polyhedral relaxations that convexify the trigonometric terms in the power flow equations. In contrast to the approach in~\cite{sundar2022sequence}, which requires solving a series of linear programs to identify the convex hull of a univariate function, our proposed method does not necessitate solving any optimization problems to construct convex envelopes for the trigonometric function.

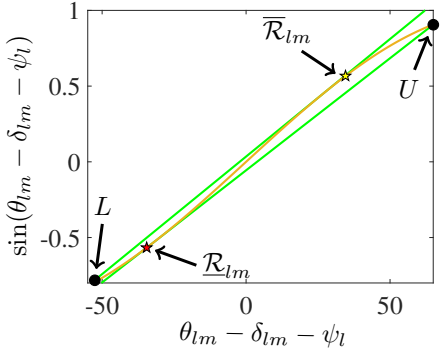
\begin{figure}
    \centering
\begin{tikzpicture}
\node at (0,0) {\includegraphics[scale=0.4]{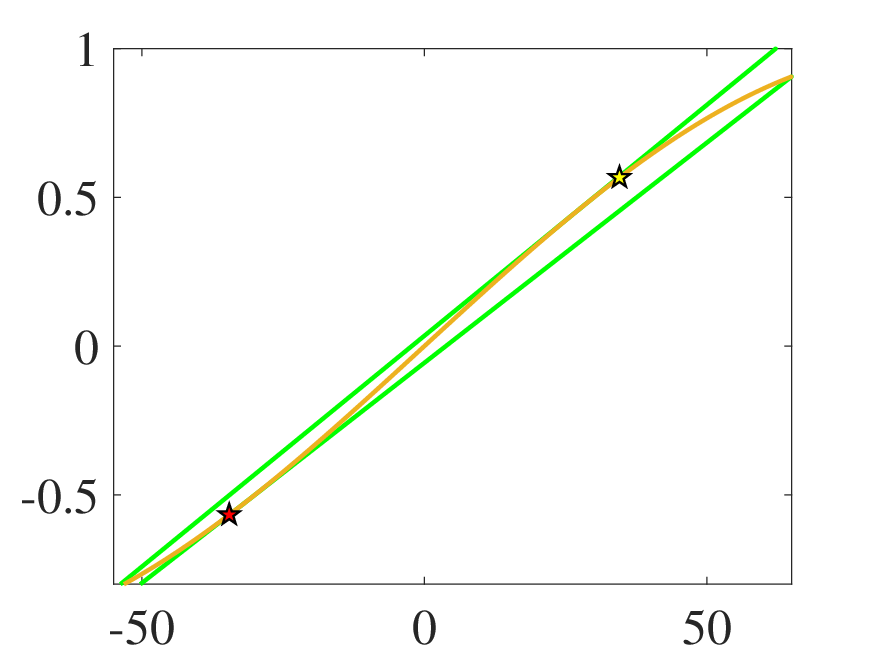}};
\node [anchor=west, rotate=-1, font=\normalsize] at (-1.1,-2.4) {$\theta_{lm}-\delta_{lm}-\psi_l$};
\node [anchor=west, rotate=90, font=\normalsize] at (-3.1,-1.5) {$\sin(\theta_{lm}-\delta_{lm}-\psi_l)$};

\node [anchor=west, font=\normalsize] at (-1,-1.5) {\begin{tabular}{c}$\underline{\mathcal{R}}_{lm}$\end{tabular}};
\node [anchor=west, font=\normalsize] at (-0.2,1.6) {\begin{tabular}{c}$\overline{\mathcal{R}}_{lm}$ \end{tabular}};
     	        \draw [-> , line width=0.4mm, black] (-0.75,-1.5) -- (-1.3,-1.3);
        \draw [ -> , line width=0.4mm, black] (0.65,1.35) -- (1.1,1.05);
        %%%%%%%%%%%%%%%%%%%
\node [anchor=west, font=\normalsize] at (-2.45,-0.75) {\begin{tabular}{c}$L$\end{tabular}};
\node [anchor=west, font=\normalsize] at (1.6,0.8) {\begin{tabular}{c}$U$ \end{tabular}};
     	        \draw [-> , line width=0.4mm, black] (-2,-0.95) -- (-2.08,-1.55);
        \draw [ -> , line width=0.4mm, black] (2.15,1.1) -- (2.33,1.55);
        %%%%%%%%%%%%%%%%%%
 \filldraw [black] (-2.1,-1.7) circle (2pt);
 \filldraw [black] (2.4,1.7) circle (2pt);
\end{tikzpicture}%
%\vspace{-.4cm}
	\caption{Envelope for $\sin(\theta_{lm}-\delta_{lm}-\psi_l)$, $\theta_{lm} \in \text{[}-60^\circ, 60^\circ \text{]}$, $\delta_{lm}=-5^\circ$, and $\psi_l= 0^\circ$. The curvature of the sine function changes within the interval considered here. $\underline{\mathcal{R}}_{lm}$ and $\overline{\mathcal{R}}_{lm}$ are shown by red and yellow stars, respectively. Tangent lines to these points from the endpoints of the interval are plotted. The endpoints of the interval, i.e., $L=\underline{\theta}_{lm}-\delta_{lm}-\psi_l$ and $U=\overline{\theta}_{lm}-\delta_{lm}-\psi_l$, are shown by black circles. The trigonometric function's curvature does not change sign within the intervals $[L,\underline{\mathcal{R}}_{lm}]$ and $[\overline{\mathcal{R}}_{lm},U]$. Thus, tangent lines to points in these ranges can be selected to form an enclosing envelope for the trigonometric function.}
	\label{fig:Enclosing_points}
	   % \vspace{-0.3cm}
\end{figure}

\iffalse
\begin{figure}
%\centering
\subfloat[Tight linear envelopes obtained by the proposed approach in Algorithm~\ref{alg:algorithm1}  for $\sin(\theta_{lm}-\delta_{lm}-\psi_i)$ function.]{\label{fig:trig_function.a}
\begin{minipage}{0.5\columnwidth}%
\centering
\begin{tikzpicture}
    \node[] (image) at (0,0) {\includegraphics[width=1.4\columnwidth]{sin.eps}};
    \begin{scope}[]
       \node [anchor=west, font=\normalsize \rotatebox{90}{$\sin(\theta_{lm}-\delta_{lm}-\psi_i)$}] at (-3.6,0.2) {};
       \node [anchor=west, font=\normalsize] at (-1.5,-2.5) {$\theta_{lm}-\delta_{lm}-\psi_i$};
    \end{scope}
\end{tikzpicture}
\end{minipage}
}\\
\subfloat[Tight linear envelopes obtained by the proposed approach in Algorithm~\ref{alg:algorithm1}  for $\sin(\theta_{lm}-\delta_{lm}-\psi_i)$ function.]{\label{fig:trig_function.b}
\begin{minipage}{0.5\columnwidth}%
\centering
\begin{tikzpicture}
    \node[] (image) at (0,0)
{\includegraphics[width=1.4\columnwidth]{cosine.eps}};
    \begin{scope}[]
          \node [anchor=west, font=\normalsize \rotatebox{90}{$\cos(\theta_{lm}-\delta_{lm}-\psi_i)$}] at (-3.6,0.2) {};
       \node [anchor=west, font=\normalsize] at (-1.5,-2.5) {$\theta_{lm}-\delta_{lm}-\psi_i$};
    \end{scope}
\end{tikzpicture}
\end{minipage}
}
\caption{Tangent lines to sine and cosine functions that accurately estimates the convex hull of these functions.}
\end{figure}
\fi

%
Convex envelopes constructed using tangent lines were also previously used to convexify the cosine function in the Linear Programming AC (LPAC) approximation proposed in~\cite{Linear_AC}. However, those envelopes are specific to arguments ranging from ${-}90^\circ$ and $90^\circ$. Since the arguments for the trigonometric functions in our formulation change with the values of $\delta_{lm}$ and $\psi_l$, we must consider ranges that admit any possible argument, including ranges for which the curvature changes. This is challenging since a tangent line to the trigonometric function at one point may intersect the function in another point, with the resulting envelope failing to enclose the function. 

This section addresses this issue by finding the largest ranges of values for which tangent lines to the trigonometric function form an enclosing envelope.
These ranges are defined by the lower bound of the argument, $\underline{\theta}_{lm} - \delta_{lm} - \psi_l$, to a point denoted as $\underline{\mathcal{R}}_{lm}$ as well as from a point denoted as $\overline{\mathcal{R}}_{lm}$ to the upper bound of the argument, $\overline{\theta}_{lm} - \delta_{lm} - \psi_l$.
Fig.~\ref{fig:Enclosing_points} shows $\underline{\mathcal{R}}_{lm}$ and $\overline{\mathcal{R}}_{lm}$ via red and yellow stars, respectively. More specifically, to assist in finding $\underline{\mathcal{R}}_{lm}$ and $\overline{\mathcal{R}}_{lm}$, we define a function $F(\theta_{lm})$ which represents the difference between the trigonometric function $\cos(\theta_{lm} - \delta_{lm}-\psi_l)$ itself and the line which connects the endpoints of $\cos(\theta_{lm} - \delta_{lm}-\psi_l)$ at $\underline{\theta}_{lm}$ and $\overline{\theta}_{lm}$:
\begin{align}
\label{eq:auxi_equation1}
\nonumber \noindent &F(\theta_{lm})  = \cos(\theta_{lm} - \delta_{lm}-\psi_l) - \cos(\overline{\theta}_{lm}-\delta_{lm}-\psi_l)\\
\nonumber &\qquad\quad-\frac{\cos(\overline{\theta}_{lm}-\delta_{lm}-\psi_l)-\cos(\underline{\theta}_{lm}-\delta_{lm}-\psi_l)}{\overline{\theta}_{lm}-\underline{\theta}_{lm}}\\&\qquad\quad\times\left(\theta_{lm}-\overline{\theta}_{lm}-\psi_l\right).
\end{align}
The set of zeros of the first derivative of $F(\theta_{lm})$, i.e., the set of solutions to $\frac{dF(\theta_{lm})}{d \theta_{lm}} = 0$, is a key quantity to determine if the curvature of  $\cos(\theta_{lm}-\delta_{lm}-\psi_l)$ changes between $\underline{\theta}_{lm}-\delta_{lm}-\psi_l$ and $\overline{\theta}_{lm}-\delta_{lm}-\psi_l$. 

If $\frac{dF(\theta_{lm})}{d \theta_{lm}} = 0$ has no solutions, then the curvature of the trigonometric function does not change between $\underline{\theta}_{lm} - \delta_{lm} - \psi_l$ and $\overline{\theta}_{lm} - \delta_{lm} - \psi_l$. Accordingly, any tangent lines can be selected to form an enclosing envelope for the trigonometric function. We select equally spaced tangent lines within the range $[\underline{\theta}_{lm} - \delta_{lm} - \psi_l, \overline{\theta}_{lm} - \delta_{lm} - \psi_l]$ as illustrated in Fig.~\ref{fig:trig_function.a}.

Conversely, if $\frac{dF(\theta_{lm})}{d \theta_{lm}} = 0$ has one or more solutions, then the trigonometric function's curvature changes. This necessitates special consideration, i.e., finding $\underline{\mathcal{R}}_{lm}$ and $\overline{\mathcal{R}}_{lm}$, to select appropriate tangent lines to $\cos(\theta_{lm} - \delta_{lm}-\psi_l)$.

To compute $\underline{\mathcal{R}}_{lm}$ and $\overline{\mathcal{R}}_{lm}$ for the cosine function, we first identify the tangent line to the cosine function that also passes through the endpoint $\overline{\theta}_{lm}-\delta_{lm}-\psi_l$. We then define another auxiliary function representing the difference between this tangent line and the cosine function. The value of $\underline{\mathcal{R}}_{lm}$ is given by the root of the first derivative of this auxiliary function that is between $\overline{\theta}_{lm}-\delta_{lm}-\psi_l$ and $\underline{\theta}_{lm}-\delta_{lm}-\psi_l$. Note that the voltage angle difference restriction, i.e., ${-}90^\circ \le \theta_{lm} \le 90^\circ$, ensures that the sine and cosine functions have at most one curvature sign change in any given interval. $\overline{\mathcal{R}}_{lm}$ is computed similarly by formulating the tangent line to the cosine function that also pass through the endpoint $\underline{\theta}_{lm}-\delta_{lm}-\psi_l$ and following the steps above. A comprehensive explanation of how to compute $\underline{\mathcal{R}}_{lm}$ and $\overline{\mathcal{R}}_{lm}$ is available in 
\ifarxiv
the appendix.
\else
\cite[Appendix]{narimani2023tightening}.
\fi

By construction, the trigonometric function's curvature does not change sign within the intervals $[\underline{\theta}_{lm}-\delta_{lm}-\psi_l,\underline{\mathcal{R}}_{lm}]$ and $[\overline{\mathcal{R}}_{lm},\overline{\theta}_{lm}-\delta_{lm}-\psi_l]$. Accordingly, tangent lines to points in these ranges can be selected to form an enclosing envelope for the trigonometric function. As shown in Fig.~\ref{fig:trig_function.b}, we choose equally spaced tangent lines within each of these ranges.

Our proposed QC relaxation uses envelopes $\left\langle \sin(\theta_{lm}-\delta_{lm}-\psi_l)\right\rangle^{S^{\prime}}$ and $\left\langle \cos(\theta_{lm}-\delta_{lm}-\psi_l)\right\rangle^{C^{\prime}}$ based on the tangent lines described above. The formulations of the upper and lower bounds of these envelopes depend on the curvature's sign and the number of solutions for $\frac{dF(\theta_{lm})}{d \theta_{lm}} = 0$. For brevity, we present a summary of the envelopes for the cosine function. Future details for the cosine function along with expressions for the sine envelopes are given \ifarxiv in the appendix.
\else in~\cite[Appendix]{narimani2023tightening}.
\fi  
\begin{small}
\begin{subequations}
\label{eq:convex_envelopes_sin&cos2}
\begin{align}
%\raisetag{18pt} 
\label{eq:cos envelope11}
&\nonumber\qquad\qquad\quad\text{If $\frac{dF(\theta_{lm})}{d \theta_{lm}} = 0$ has one or more solutions: }\\
 &\left\langle \cos(\theta_{lm}-\delta_{lm}-\psi_l)\right\rangle^{C^{\prime}} =
\begin{cases}
\!\widecheck{C}^{\prime}\!:\begin{cases}
\widecheck{C}^{\prime}\! \le \overline{L}_{\cos,i},~i=1,\ldots,N_{tan}\\
\widecheck{C}^{\prime}\! \ge \underline{L}_{\cos,i},~i=1,\ldots,N_{tan}\end{cases}
\end{cases}\\
%\raisetag{18pt} 
\label{eq:cosine envelope12}
&  \nonumber\quad\text{If $\frac{dF(\theta_{lm})}{d \theta_{lm}} = 0$ has no solutions \& curvature sign is negative: }\\
&\left\langle\cos(\theta_{lm}-\delta_{lm}-\psi_l)\right\rangle^{C^{\prime}} =
\begin{cases}
\!\widecheck{C}^{\prime}\!:\begin{cases}\widecheck{C}^{\prime}\! \le 
\overline{L}_{\cos,i},~i=1,\ldots,N_{tan}\\
%\bigcap\limits_{i=1}^{N_{upper}}\overline{L}_{\cos,i}\\
\widecheck{C}^{\prime}\! \ge \underline{L}_{\cos,i},~i=1\end{cases}
\end{cases}\\
%\raisetag{18pt} 
\label{eq:cosine envelope13}
& \nonumber\quad\text{If $\frac{dF(\theta_{lm})}{d \theta_{lm}} = 0$ has no solutions \& curvature sign is positive: } \\
&\left\langle\cos(\theta_{lm}-\delta_{lm}-\psi_l)\right\rangle^{C^{\prime}} =
\begin{cases}
\!\widecheck{C}^{\prime}\!:\begin{cases}\widecheck{C}^{\prime}\! \le 
\overline{L}_{\cos,i},~i=1\\
%\bigcap\limits_{i=1}^{N_{upper}}\overline{L}_{\cos,i}\\
\widecheck{C}^{\prime}\! \ge \underline{L}_{\cos,i},~i=1,\ldots,N_{tan}\end{cases}
\end{cases}%
\end{align}%
\end{subequations}%
\end{small}%
where $\overline{L}_{\sin,i}$, $\overline{L}_{\cos,i}$, $\underline{L}_{\sin,i}$, and $\underline{L}_{\cos,i}$ are the $i^{th}$ tangent lines which upper and lower bound the sine and cosine functions, respectively. When the sign of the trigonometric function's curvature does not change within an interval, either the upper or lower boundary of the envelope (depending on the sign of the curvature) is defined via the line connecting the endpoints of the trigonometric function, as defined in~\eqref{eq:cosine envelope12} and~\eqref{eq:cosine envelope13}.%

The envelopes in~\eqref{eq:convex_envelopes_sin&cos2} are valid for any argument $\theta_{lm}-\delta_{lm}-\psi_l$. 
The lifted variables $\widecheck{S}^{\prime}$ and $\widecheck{C}^{\prime}$ are associated with the envelopes for the functions $\sin(\theta_{lm}-\delta_{lm}-\psi_l)$ and $\cos(\theta_{lm}-\delta_{lm}-\psi_l)$. Detailed expressions for $\overline{L}_{\sin,i}$, $\overline{L}_{\cos,i}$, $\underline{L}_{\sin,i}$, and $\underline{L}_{\cos,i}$ are available in 
\ifarxiv
the appendix.
\else
\cite[Appendix]{narimani2023tightening}.
\fi

\section{The Linear Rotated QC Relaxation}
\label{TLRQC_relaxation}

This section brings together each improvement from this paper (multiple angle rotations associated with each bus in Section~\ref{Rotated_OPF}, new envelopes for the product terms in Section~\ref{Higher_dimansional_space}, and tighter trigonometric envelopes in Section~\ref{Envelopes_trig_terms}) to formulate our proposed QC relaxation in~\eqref{eq:LRQC OPF} below. 
We subsequently call this formulation the ``Linear Rotated QC'' (LRQC) relaxation since the polytopes for the convex envelopes are constructed with linear inequalities and the power flow equations are rotated versions of the original expressions. 
\begin{subequations}
\label{eq:LRQC OPF}
\begin{align}
& \min\quad \eqref{eq:RQC obj}\\
&\nonumber \text{subject to} \quad \left(\forall i\in\mathcal{N}, \forall   \left(l,m\right) \in\mathcal{L}\right)\\
&\text{Equations~}\eqref{eq:RQC active power injection}\text{--}\eqref{eq:shifted OPF power flow limit lm},\\
\label{eq:RQC active power flow lm2}
& \nonumber \tilde{P}_{lm} \!=\! \left(Y_{lm}\cos(\delta_{lm} \!+\! \psi_l) - b_{c,lm}/2 \sin(\psi_l)\right)w_{ll}\\  &\qquad - Y_{lm}\tilde{c}^{\prime}_{lm},\\[-0.35em]
  \label{eq:RQC reactive power flow lm2}
& \nonumber \tilde{Q}_{lm}\!=\! -\left(Y_{lm}\sin(\delta_{lm}\!+\!\psi_l)+b_{c,lm}/2\cos(\psi_l)\right)w_{ll}\\  
  &\qquad -Y_{lm}\tilde{s}^{\prime}_{lm},\\[-0.35em]
  \label{eq:RQC active power flow ml2}
& \tilde{P}_{ml}\!=\!
  -Y_{lm}\tilde{c}^{\prime}_{lm} \!+\! \left( Y_{lm} \cos(\delta_{lm}\!+\! \psi_l) - b_{c,lm}/2 \sin(\psi_l)\right)w_{mm},\\[-0.35em]
  \label{eq:RQC reactive power flow ml2}
& \tilde{Q}_{ml} \!=\!
  Y_{lm}\tilde{s}^{\prime}_{lm} - \left(Y_{lm}  \sin(\delta_{lm}\!+\!\psi_l) + b_{c,lm}/2 \cos(\psi_l)\right)w_{mm},\\
 \label{eq:RQC squared current2}
& \tilde{P}_{lm}^2 + \tilde{Q}_{lm}^2 \leq w_{ll}\, \tilde{\ell}^{\prime}_{lm},\\
 \label{eq:RQC squared current expression2}
&\nonumber \tilde{\ell}^{\prime}_{lm} =\bigg(b_{c,lm}^2/4+Y_{lm}^2- Y_{lm}b_{c,lm}\cos(\delta_{lm}+\psi_l)\sin(\psi_l)\\[-0.35em]
&\nonumber\quad\quad+Y_{lm}b_{c,lm}\sin(\delta_{lm}+\psi_l)\cos(\psi_l)\bigg)V_l^2+Y_{lm}^2 V_m^2\\&\quad\quad\nonumber+\left(-2Y_{lm}^2\cos(\delta_{lm}+\psi_l)+Y_{lm}b_{c,lm}\sin(\psi_l)\right)\tilde{c}^{\prime}_{lm}\\&\quad\quad
+\left(2Y_{lm}^2\sin(\delta_{lm}+\psi_l)+Y_{lm}b_{c,lm}\cos(\psi_l)\right)\tilde{s}^{\prime}_{lm},\\
\nonumber & \tilde{c}^{\prime}_{lm} = \!\!\!\!\!\!\!\! \sum_{k=1,\ldots,4N_{seg}} \!\!\!\!\!\lambda_{lm,k}\, \xi^{(k)}_{lm,1} \xi^{(k)}_{lm,2} \xi^{(k)}_{lm,4},~~~ V_l =\!\!\!\!\!\!\!\!\! \sum_{k=1,\ldots,4N_{seg}} \!\!\!\!\!\!\!\!\lambda_{lm,k}\xi^{(k)}_{lm,1}, \\
\nonumber & \tilde{s}^{\prime}_{lm} = \!\!\!\!\!\!\!\! \sum_{k=1,\ldots,4N_{seg}} \!\!\!\!\!\lambda_{lm,k}\, \xi^{(k)}_{lm,1} \xi^{(k)}_{lm,2} \xi^{(k)}_{lm,5},\\
\nonumber &\theta_{lm} =\!\!\!\!\!\!\!\!\!\! \sum_{k=1,\ldots,N_{seg}} \!\!\!\!\!\!\!\lambda_{lm,k}\eta^{(k)}_{lm,3},~~~~V_m =\!\!\!\!\! \sum_{k=1,\ldots,4N_{seg}}\!\!\!\!\! \lambda_{lm,k}\xi^{(k)}_{lm,2},\\
    \nonumber & \tilde{S}_{lm}^{\prime} =\!\!\!\!\!\!\!\! \sum_{k=1,\ldots,4N_{seg}}\!\!\!\!\! \lambda_{lm,k}\xi^{(k)}_{lm,5},~~~~\tilde{C}_{lm}^{\prime} =\!\!\!\!\!\!\!\!\!\! \sum_{k=1,\ldots,4N_{seg}}\lambda_{lm,k}\xi^{(k)}_{lm,4},\\
     \nonumber&
     \!\!\! \sum_{k=1,\ldots,4N_{seg}} \!\!\!\!\lambda_{lm,k} = 1, \qquad \!\!\!\!\!\!\lambda_{lm,k} \geqslant 0, \quad \!\!\!k=1,\ldots,4N_{seg},\\
     \label{eq:convex hull trilinear terms2} & \tilde{C}_{lm}^{\prime} \in \left\langle\cos(\theta_{lm}\!-\!\delta_{lm}\!-\!\psi_l)\right\rangle^{C^{\prime}}\!\!, \; \tilde{S}_{lm}^{\prime} \in\left\langle\sin(\theta_{lm}\!-\!\delta_{lm}\!-\!\psi_l)\right\rangle^{S^{\prime}}\!\!. 
\end{align}
\end{subequations}

The lifted variables $\tilde{c}^{\prime}_{lm}$ and $\tilde{s}^{\prime}_{lm}$ represent relaxations of the product terms $V_lV_m\cos(\theta_{lm}-\delta_{lm}-\psi_l)$ and $V_lV_m\sin(\theta_{lm}-\delta_{lm}-\psi_l)$, respectively, with~\eqref{eq:convex hull trilinear terms2} formulating an ``extreme point'' representation of the convex hulls for the product terms $V_l V_m \widecheck{C}^{\prime}_{lm}\widecheck{S}^{\prime}_{lm}$.
The extreme points of $V_lV_m\widecheck{C}^{\prime}_{lm}\widecheck{S}^{\prime}_{lm}$ are $\xi^{(k)} \in [\underline{V_l},\overline{V_l}] \times [\underline{V_m},\overline{V_m}]\times\mathcal{T}_{lm}$, $k=1,\ldots,4N_{seg}$.  
$\mathcal{T}_{lm}$ denotes the coordinates of the extreme points (red squares) in Fig.~\ref{fig:sin_cos}

\ifarxiv 
 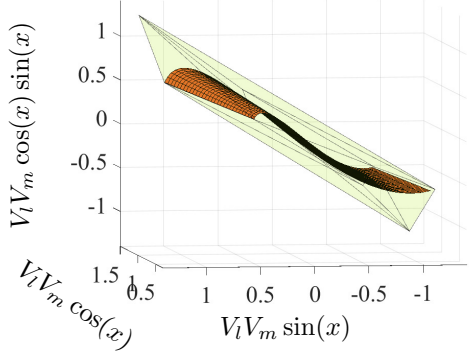
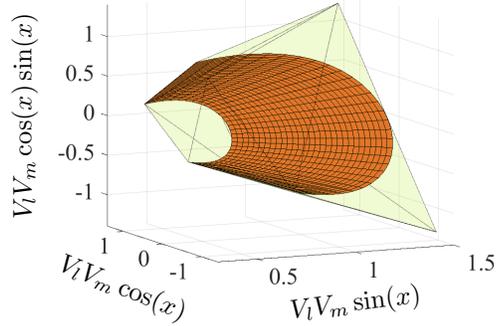
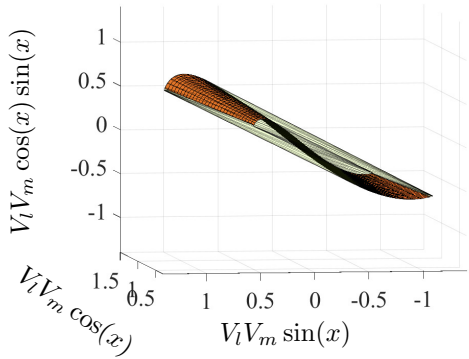
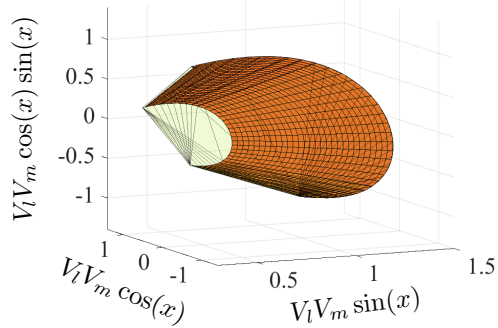
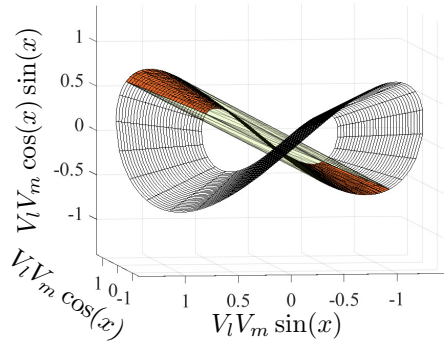
\begin{figure*}
\centering
\def\tabularxcolumn#1{m{#1}}
\begin{tabularx}{\linewidth}{@{}cXX@{}}
\begin{tabular}{cc}
\subfloat[Envelope from the original QC relaxation in~\cite{Sundar}.]{\label{fig4.a}
\centering
\begin{tikzpicture}
    \node[] (image) at (0,0) {\includegraphics[width=0.67\columnwidth]{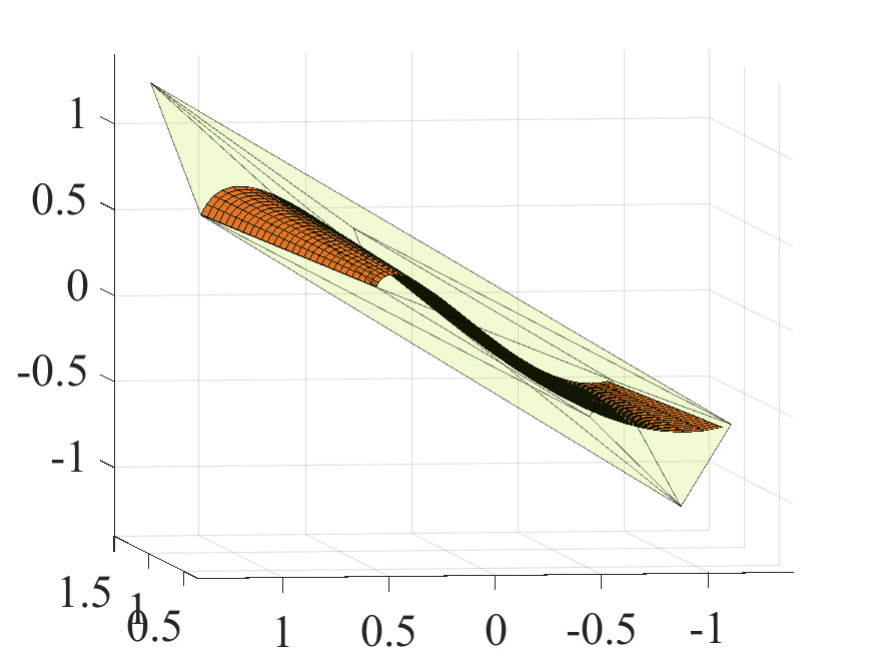}};
       \node [ font=\normalsize] at (0.0,-2.5) {$V_lV_m\sin(x)$};
       \node [ font=\normalsize \rotatebox{-35}{$V_lV_m\cos(x)$}] at (-2.8,-2.2) {};
       \node [ font=\normalsize \rotatebox{90}{$V_lV_m\cos(x)\sin(x)$}] at (-3.5,0.2) {};
\end{tikzpicture}
}%
&\subfloat[Envelope from the original QC relaxation in~\cite{Sundar} (alternate view).]{\label{fig4.b}
\centering
\begin{tikzpicture}
    \node[]  at (-0.8,0)  {\includegraphics[width=0.67\columnwidth]{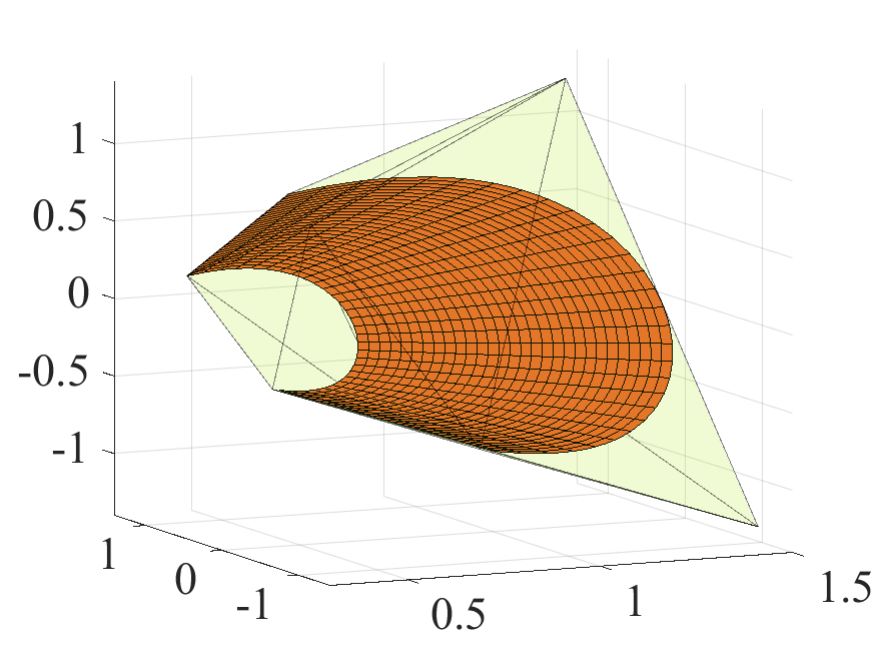}};
       \node [ font=\normalsize \rotatebox{7} {$V_lV_m\sin(x)$}] at (0.3,-2.3) {};
       \node [ font=\normalsize \rotatebox{-22}{$V_lV_m\cos(x)$}] at (-2.8,-2.1) {};
       \node [ font=\normalsize \rotatebox{90}{$V_lV_m\cos(x)\sin(x)$}] at (-4.1,0.2) {};
\end{tikzpicture}
}
%\vspace{-.4cm}
\\ 
\subfloat[Envelope from the proposed LRQC relaxation~\eqref{eq:LRQC OPF}.]{\label{fig4.c}
\centering
\centering
\begin{tikzpicture}
    \node[] (image) at (0,0) {\includegraphics[width=0.67\columnwidth]{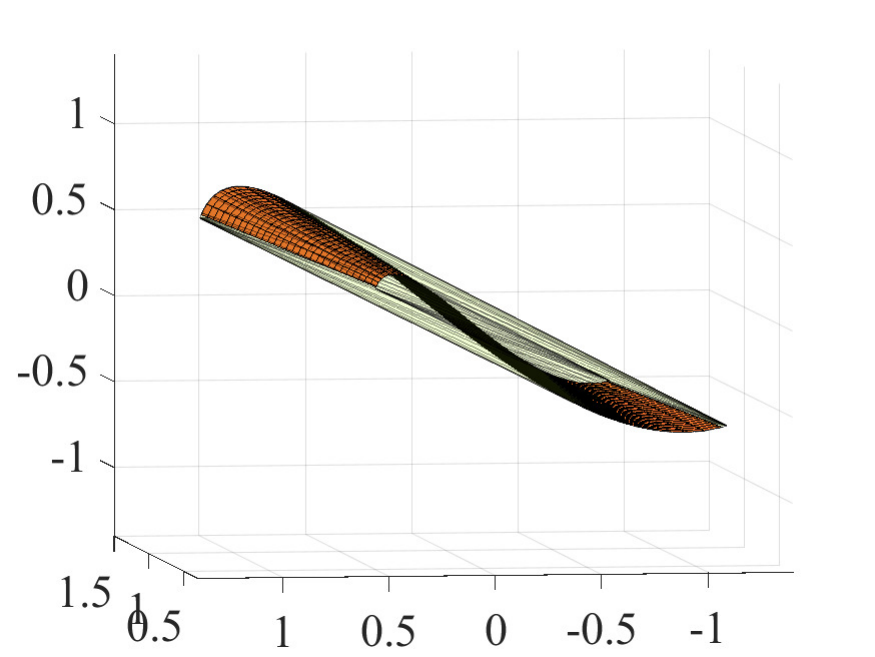}};
       \node [ font=\normalsize] at (0.0,-2.5) {$V_lV_m\sin(x)$};
       \node [ font=\normalsize \rotatebox{-35}{$V_lV_m\cos(x)$}] at (-2.8,-2.2) {};
       \node [ font=\normalsize \rotatebox{90}{$V_lV_m\cos(x)\sin(x)$}] at (-3.5,0.2) {};
\end{tikzpicture}
}
& 
\subfloat[Envelope from the proposed LRQC relaxation~\eqref{eq:LRQC OPF} (alternate view).]{\label{fig4.d}
\centering
\centering
\begin{tikzpicture}
    \node[] (image) at (-0.8,0) {\includegraphics[width=0.67\columnwidth]{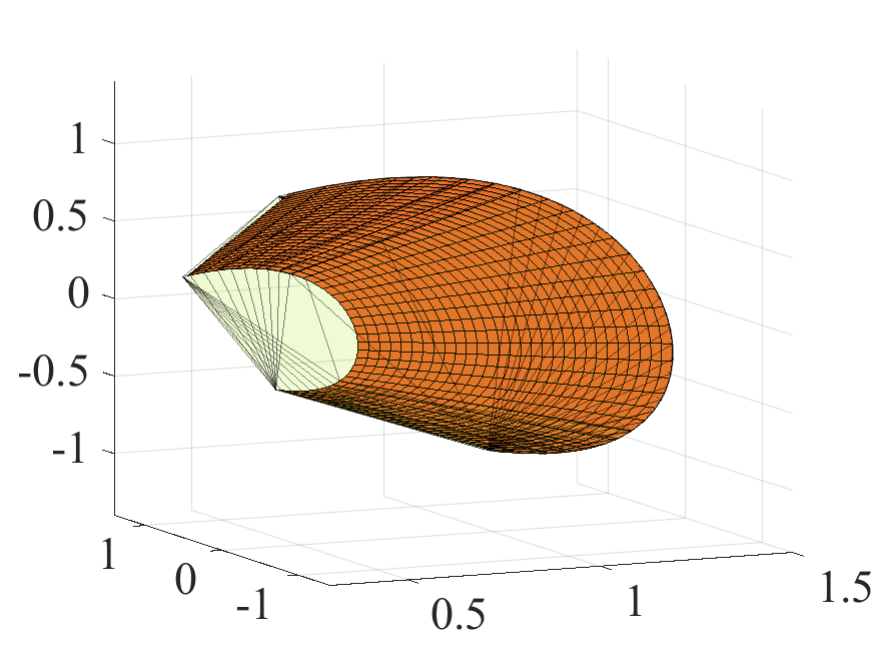}};
    \begin{scope}[]
       \node [ font=\normalsize \rotatebox{7} {$V_lV_m\sin(x)$}] at (0.3,-2.3) {};
       \node [ font=\normalsize \rotatebox{-22}{$V_lV_m\cos(x)$}] at (-2.8,-2.1) {};
       \node [ font=\normalsize \rotatebox{90}{$V_lV_m\cos(x)\sin(x)$}] at (-4.1,0.2) {};
    \end{scope}
\end{tikzpicture}
}
\\ 
\captionsetup{width=1.0\linewidth, margin={0cm,-9cm}}
\subfloat[Envelope from the proposed LRQC relaxation~\eqref{eq:LRQC OPF}. The black and white region shows the entire function $V_lV_m\cos(x)\sin(x)$ while the orange region is the portion of this function within the voltage magnitude and phase angle difference bounds as in Figs.~\ref{fig4.a}--\ref{fig4.d}.]{\label{fig4.e}
\begin{tikzpicture}
  \hspace{4.5cm}
    \node[] (image) at (0,0)
    {\includegraphics[width=0.67\columnwidth]{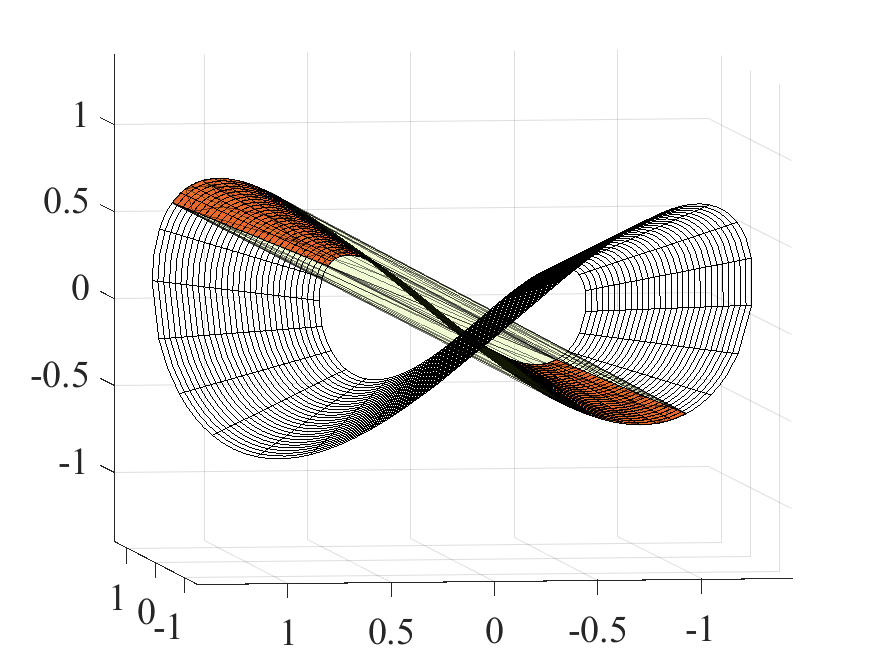}};
    \hspace{1.0cm}
    \begin{scope}[]
       \node [ font=\normalsize] at (-.8,-2.4) {$V_l V_m\sin(x)$};
       \node [ font=\normalsize \rotatebox{-35}{$V_l V_m\cos(x)$}] at (-3.6,-2.0) {};
       \node [ font=\normalsize \rotatebox{90}{$V_l V_m\cos(x)\sin(x)$}] at (-4.1,0.2) {};
    \end{scope}
\end{tikzpicture}
}
 \end{tabular}
\end{tabularx}
\caption{Projections of the function $V_lV_m\cos(x)\sin(x)$ in terms of $V_l V_m \cos(x)$ and $V_l V_m \sin(x)$. The argument $x$ indicates the angle difference $\theta_{lm}$ for the original QC relaxation in Figs.~\ref{fig4.a} and~\ref{fig4.b} and the rotated argument from the polar admittance representation, $\theta_{lm}-\delta_{lm}-\psi_l$, for the LRQC relaxation from~\eqref{eq:LRQC OPF} in Figs.~\ref{fig4.c}--\ref{fig4.e}. The orange region common to Figs.~\ref{fig4.a}--\ref{fig4.e} is the function $V_lV_m\cos(x)\sin(x)$ that we seek to enclose in a convex envelope. The light green regions correspond to the surfaces of the convex envelopes proposed in the original QC relaxation~\cite{coffrin2015qc} for Figs.~\ref{fig4.a} and \ref{fig4.b} and the proposed LRQC relaxation~\eqref{eq:LRQC OPF} from Section~\ref{TLRQC_relaxation} for Figs.~\ref{fig4.c}--\ref{fig4.e}. Note that Figs.~\ref{fig4.b} and \ref{fig4.d} on the right side show rotated views of the same projections as Figs.~\ref{fig4.a} and~\ref{fig4.c}, respectively, on the left side.}
\end{figure*}
\fi

\ifarxiv 
To illustrate the tightness of the envelopes in the LRQC relaxation, Figs.~\ref{fig4.a}--\ref{fig4.e} show a projection of the function $V_lV_m\sin(\theta_{lm}-\delta_{lm}-\psi_l)\cos(\theta_{lm}-\delta_{lm}-\psi_l)$ along with the convex envelopes from both the approach in~\cite{Sundar} and our proposed method. 
\else
To illustrate the tightness of the envelopes in the LRQC relaxation, Figs.~\ref{fig4.a}--\ref{fig4.e} show a projection of the function $V_lV_m\sin(\theta_{lm}-\delta_{lm}-\psi_l)\cos(\theta_{lm}-\delta_{lm}-\psi_l)$ along with the convex envelopes from both the approach in~\cite{Sundar} and our proposed method.
\fi
The orange region common to each figure corresponds to different views of the function itself and the light green polytopes are the convex envelopes. 
\ifarxiv 
Figs.~\ref{fig4.a}--\ref{fig4.b} show envelopes from the original QC relaxation and Figs.~\ref{fig4.c}--\ref{fig4.d} show our proposed envelopes. Observe that our proposed envelopes can be significantly tighter than those in the original QC relaxation. Fig.~\ref{fig4.e} shows these same envelopes with the full function $V_lV_m\sin(\theta_{lm}-\delta_{lm}-\psi_l)\cos(\theta_{lm}-\delta_{lm}-\psi_l)$ where regions outside of the voltage magnitude and angle difference bounds are transparent rather than orange.
\else
Fig.~\ref{fig4.a} shows the envelope from the original QC relaxation and Fig.~\ref{fig4.c} shows our proposed envelopes. Observe that our proposed envelopes can be significantly tighter than those in the original QC relaxation. Fig.~\ref{fig4.e} shows these same envelopes with the full function $V_lV_m\sin(\theta_{lm}-\delta_{lm}-\psi_l)\cos(\theta_{lm}-\delta_{lm}-\psi_l)$,  where regions outside of the voltage magnitude and angle difference bounds are transparent rather than orange.
\fi

\section{Choosing the Rotation Angles}
\label{Best_Rotation_Angle}

\begin{figure}[t!]
%\vspace{-.7cm}
    \centering
\begin{tikzpicture}
\node at (0,0) {\includegraphics[scale=0.41]{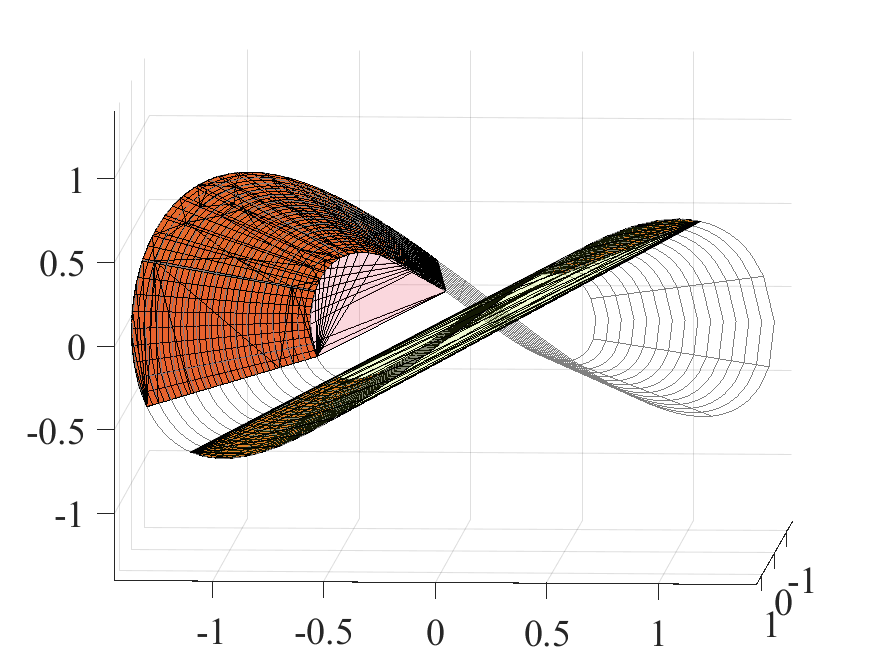}};
\node [anchor=west, rotate=-1, font=\normalsize] at (-1.0,-2.8) {$V_lV_m\sin(x)$};
\node [anchor=west, rotate=65, font=\normalsize] at (2.7,-3.0) {$V_lV_m\cos(x)$};
\node [anchor=west, rotate=90, font=\normalsize] at (-3.5,-1.3) {$V_lV_m\sin(x)\cos(x)$};
\node [anchor=west, font=\normalsize] at (-1.85,1.9) {\begin{tabular}{c}$\psi_l=-67.8^\circ$ \\ $\delta_{lm}=83^\circ$\end{tabular}};
\node [anchor=west, font=\normalsize] at (0.2,-1.3) {\begin{tabular}{c}$\psi_l=-174^\circ$ \\ $\delta_{lm}=83^\circ$\end{tabular}};
     	        \draw [-> , line width=0.4mm, black] (0.4,-1.1) -- (0.2,-0.05);
        \draw [ -> , line width=0.4mm, black] (-0.8,1.5) -- (-0.5,0.2);
\end{tikzpicture}%
%\vspace{-.4cm}
	\caption{Projection of $V_lV_m\cos(x)\sin(x)$ (orange regions), where $x$ indicates the $\theta_{lm}-\delta_{lm}-\psi_l$, for different values of $\psi_l$. The associated envelopes in light orange and light green show how choosing $\psi_l$ affects the convex envelopes.}
	\label{fig_volume}
	   %\vspace{-0.5cm}
\end{figure}

The rotation angles~$\psi_l$ play an important role in the performance of the proposed LRQC relaxation~\eqref{eq:LRQC OPF}. Since the admittance angles $\delta_{lm}$ vary between branches, different rotation angles $\psi_l$ may yield tighter envelopes for the trigonometric terms $\cos(\theta_{lm}-\delta_{lm}-\psi_l)$ and $\sin(\theta_{lm}-\delta_{lm}-\psi_l)$. To illustrate the impact of the rotation angle $\psi_l$ on the convex envelopes, Fig.~\ref{fig_volume} shows two envelopes associated with different choices of $\psi_l$. This section proposes and analyzes a heuristic approach for choosing the rotation angle $\psi_l$ for each bus. This heuristic is based on minimizing the convex envelopes' volumes using the intuition that smaller volumes correspond to tighter envelopes. The results in Section~\ref{Numerical_results} show this heuristic's merits via improved optimality gaps for various test cases.

For each bus, we determine the best $\psi_l$ by calculating the summation of volumes associated with the convex envelopes that enclose $V_lV_m\cos(\theta_{lm}-\delta_{lm}-\psi_l)\sin(\theta_{lm}-\delta_{lm}-\psi_l)$ terms for all the lines connected to bus $l$. To this end, we begin by sweeping the value of $\psi_l$ from $-90^\circ$ to $90^\circ$ in $1^\circ$ increments. We then compute the volume of the polytope depicted in Fig.~\ref{fig_volume}  for all lines connected to bus $l$. Finally, we choose the rotation angle for each bus based on the minimum sum of the volumes.

Finding the volume-minimizing rotation angle for each bus is a time consuming process especially for larger test systems due to the need to perform many volume computations.
Since this volume-minimization heuristic only requires knowledge of the line admittances connected to each bus, the volume computations can be performed once offline and can be reused for multiple OPF problems with the same system so long as the topology remains unchanged. Furthermore, the evaluation of this heuristic can be performed in parallel for each line. If the topology does change, only the values of $\psi_l$ associated with buses~$l$ that are directly associated with the modified topology need to be updated. Thus, while potentially time consuming in its first evaluation, we anticipate this heuristic would nevertheless be practically relevant. However, if one wishes to avoid time-consuming offline computations, we observed that most of the resulting rotation angles $\psi_l$ for the PGLib-OPF test cases are in the intervals $[-90^\circ,-85^\circ]$ and $[85^\circ,90^\circ]$. The numerical results indicate that selecting a value of $\psi^{\dagger}_{l}=-85^\circ$ for all buses~$l$ yields small optimality gaps for most test cases, suggesting that this value could be used directly with limited impacts on the relaxation's tightness.

\section{Choosing The Number of Extreme Points\\ for the $\sin(x)
\cos(x)$ Envelopes}
\label{intersection_volume}
\mrn{Since our convex envelopes are polytopes, the associated constraints in the LRQC formulation are linear. This contrasts with prior QC relaxations, where more computationally complex convex quadratic constraints are commonly used. Our formulation also enables tailoring the tightness of these envelopes by adjusting the number of segments in the polytopes to balance tractability and tightness. This section presents an analytical assessment regarding this trade-off. An empirical assessment is provided by the numeric results in Section~\ref{sec:results_Accurracy_Time}.}
%, specifically, in Tables~\ref{tab:optimality_gap_Nseg} and~\ref{tab:execution_time_Nseg}.}

\begin{figure}[t]
    \centering
\begin{tikzpicture}
\node at (0,0) {\includegraphics[scale=0.38]{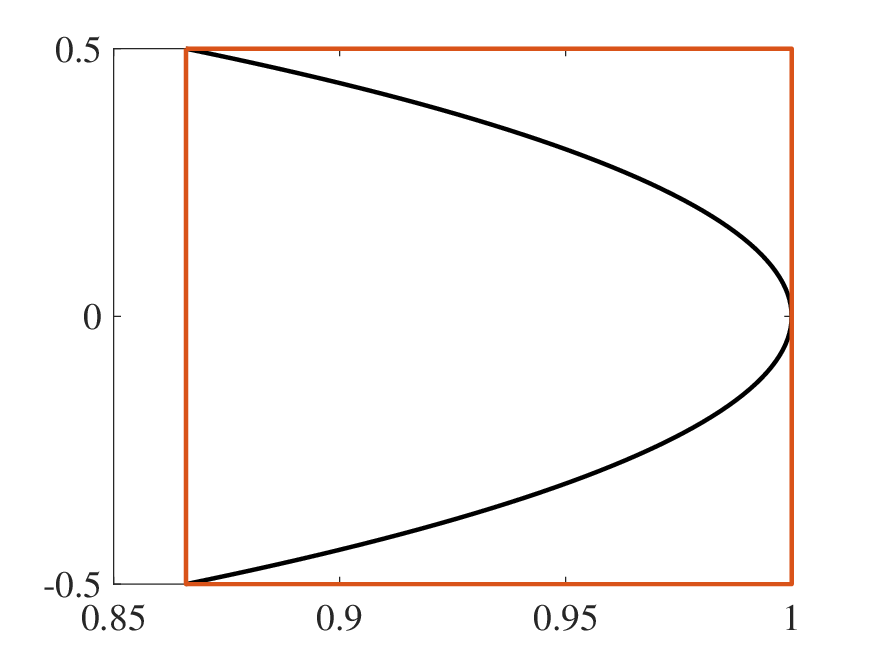}};
\node [anchor=west, rotate=-1, font=\normalsize] at (-0.5,-2.3) {$\cos(x)$};
\node [anchor=west, rotate=90, font=\normalsize] at (-2.7,-0.7) {$\sin(x)$};
\end{tikzpicture}%
%\vspace{-.4cm}
	\caption{\mrn{Envelope for the function $\cos(x)\sin(x)$ in terms of $\cos(x)$ and $\sin(x)$. The argument $x$ indicates  $\theta_{lm}$ for the original QC relaxations. The black curve is the function $\cos(x)\sin(x)$ that is enclosed in a convex envelope shown by the red lines.}}
	\label{fig_envelopes_emprical}
	% \vspace{-0.5cm}
\end{figure}

\begin{figure}[t]
\centering
%\begin{figure}
\def\tabularxcolumn#1{m{#1}}
\begin{tabularx}{\linewidth}{@{}cXX@{}}
\begin{tabular}{cc}
\subfloat[$N_{seg} = 3$]{\includegraphics[width=4.0cm]{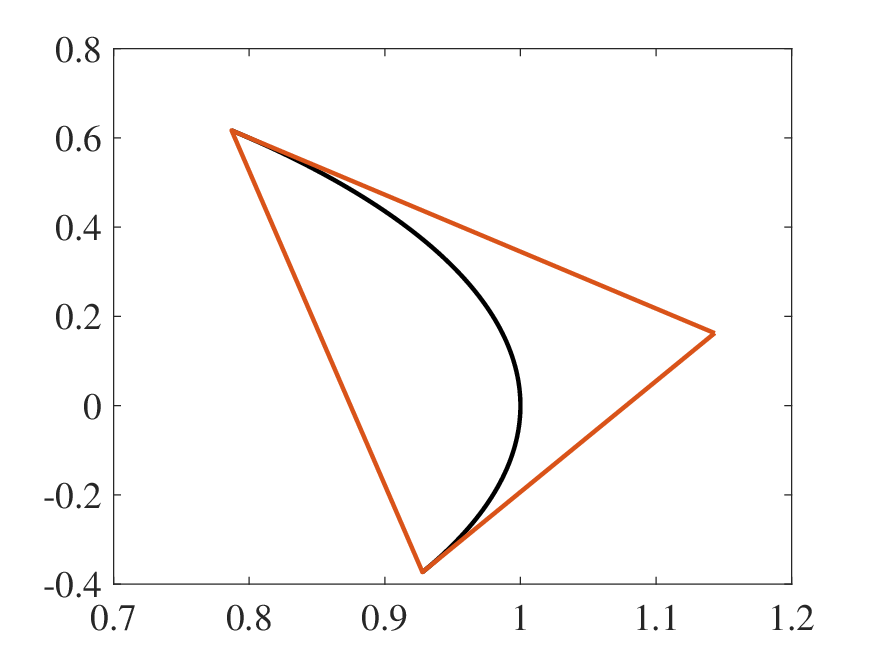}\label{fig11}} 
   &\subfloat[$N_{seg} = 6$]{\includegraphics[width=4.0cm]{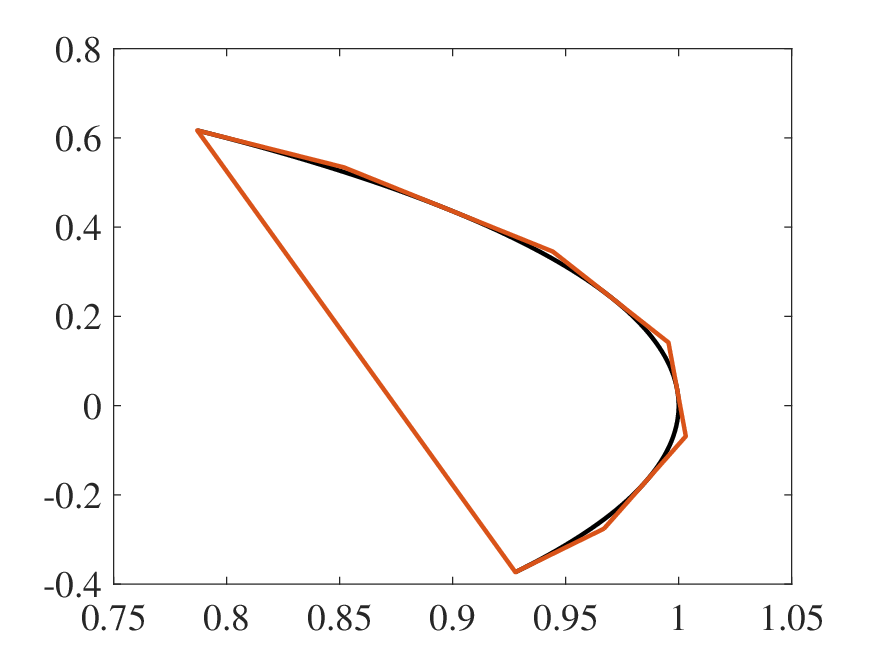}\label{fig12}}  \\
  % \vspace{-.1cm}
   \subfloat[$N_{seg} = 12$]{\includegraphics[width=4.0cm]{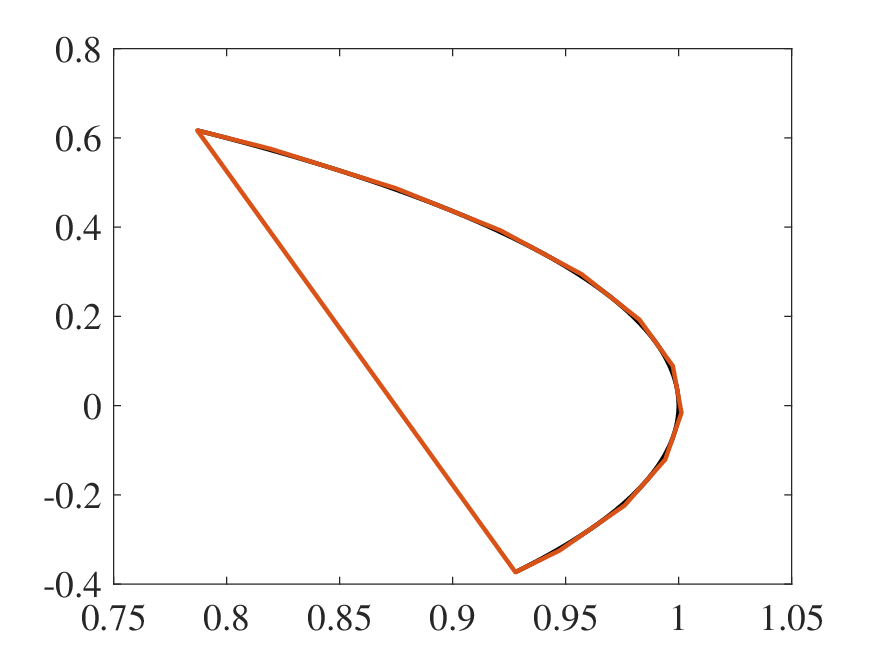}\label{fig13}} 
   &\subfloat[$N_{seg} = 22$]{\includegraphics[width=4.0cm]{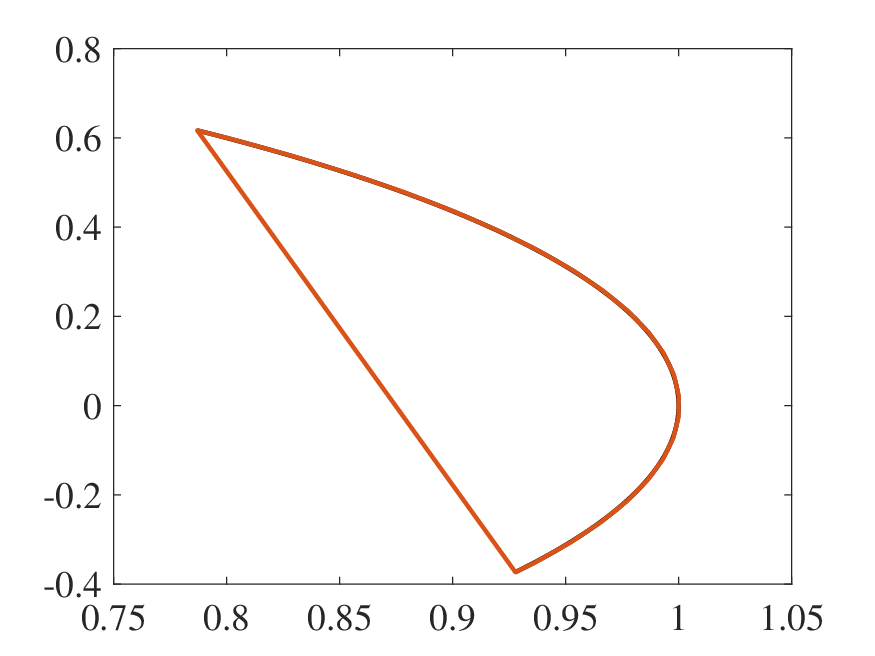}\label{fig14}}
   \end{tabular}
\end{tabularx}
\centering
\caption{\mrn{Various envelopes for the function $\cos(x)\sin(x)$ in terms of $\cos(x)$ and $\sin(x)$. The argument $x$ indicates the rotated argument from the polar admittance representation, $\theta_{lm}-\delta_{lm}-\psi_l$, for the rotated QC relaxations in Figs.~\ref{fig11}-\ref{fig14}. The black curve common to Figs.~\ref{fig11}-\ref{fig14} is the function $\cos(x)\sin(x)$ that we enclose in a convex envelope. The red envelopes in these figures are the convex envelopes in our proposed formulation in Algorithm~\ref{alg:algorithm} for (a):~$N_{seg}=3$, (b):~$N_{seg}=6$, (c):~$N_{seg}=12$, and (d):~$N_{seg}=22$.}}
\label{fig:envelops}
%\end{figure}
\end{figure}

\mrn{We next analytically characterize the tightness of the proposed envelopes for the product terms $V_lV_m\sin(x)\cos(x)$ as the number of extreme points varies. We specifically compare the volume associated with a projection of the envelope for the $V_lV_m\sin(x)\cos(x)$ terms with respect to the expression $\sin(x)\cos(x)$. The normalized volume associated with the envelopes for the $\sin(x)\cos(x)$ expression using the formulation in the original QC relaxation, as shown in Fig.~\ref{fig_envelopes_emprical}, is 0.134. The normalized volume of the relevant envelopes in the proposed LRQC relaxation depends on the number of extreme points, which is itself determined by the $N_{\text{seg}}$ parameter in Algorithm~\ref{alg:algorithm}. Figs.~\ref{fig11}-\ref{fig14} show the $\sin(x)\cos(x)$ function and its convex envelopes for $N_{\text{seg}}=3$, $N_{\text{seg}}=6$, $N_{\text{seg}}=12$, and $N_{\text{seg}}=22$, respectively. The normalized volumes enclosed by the convex envelope in these figures are 0.144, 0.093, 0.0911, and 0.0907. This indicates that increasing the number of segments from 3 to 6 significantly decreases this volume, whereas increases from 6 to 12 has a much smaller impact, suggesting diminishing returns to increasing this parameter. This is consistent with the numerical results in Section~\ref{sec:results_Accurracy_Time}.}

\mrn{As we will show empirically in Table~\ref{tab:execution_time_Nseg} and discuss in Section~\ref{sec:results_Accurracy_Time}, increasing the number of the segments from 6 to 12 can significantly increase computational times. We therefore recommend selecting $N_{\text{seg}}=5$.}

%\vspace{.5cm}
\section{Numerical Results}
\label{Numerical_results}
This section demonstrates the proposed improvements using selected test cases from the PGLib-OPF v18.08 benchmark library~\cite{pglib}. With large optimality gaps between the objective values from the best known local optima and the lower bounds from various relaxations, these test cases challenge a variety of solution algorithms and are therefore suitable for our purposes. Our implementations use Julia 0.6.4, JuMP v0.18~\cite{JuMP}, PowerModels.jl~\cite{powermodels}, and Gurobi 8.0 as modeling tools and the solver. \mrn{For comparison purposes, we also use the second-order cone programming (SOCP) relaxation from~\cite{jabr2006radial} as implemented in PowerModels.jl~\cite{powermodels} as well as a Matlab implementation of the semidefinite programming (SDP) relaxation from~\cite{molzahn2013sdp} solved with Mosek 10.1.} The results are computed using a laptop with an Intel i7 1.80~GHz processor and 16 GB of RAM. 

\subsection{Optimality Gaps and Solution Times}

Table~\ref{tab:various_QC_results} summarizes the results from applying the QC~\eqref{eq:QC relaxation}, RQC~\eqref{eq:RQC OPF}, \mrn{SOCP~\cite{jabr2006radial}, SDP~\cite{molzahn2013sdp}, and} the proposed LRQC~\eqref{eq:LRQC OPF} relaxations to selected test cases. 
To get illustrative results for the LRQC relaxation, we set $N_{seg} = 5$. The first column lists the test cases. The next group of columns represents optimality gaps as defined in~\eqref{eq:optimality gap}. The optimality gaps are computed using the local solutions to the non-convex problem~\eqref{OPF formulation} from \mbox{PowerModels.jl}:
\begin{align}
\label{eq:optimality gap}
\text{\emph{Optimality~Gap}}=\left(\dfrac{\text{\emph{Local~Solution}} - \text{\emph{QC Bound}}}{\text{\emph{Local~Solution}}}\right).
\end{align}

Upon comparing the \mrn{fifth} and \mrn{sixth} columns of Table~\ref{tab:various_QC_results}, it is evident that the RQC relaxation from our previous work in~\cite{NarimaniTPS1} outperforms the original QC relaxation for all test cases by converging to tighter lower bounds. The best rotation angle $\psi^*$ for the RQC relaxation in the \mrn{seventh} column of Table~\ref{tab:various_QC_results} is obtained by sweeping $\psi$ from $-90^\circ$ to $90^\circ$ in steps of $0.5^\circ$. The RQC relaxation in~\cite{NarimaniTPS1} with $\psi^*$ (the best value of $\psi$ for each case) provides optimality gaps that are at least as good as those obtained by the original QC relaxation~\eqref{eq:QC relaxation} for all test cases, resulting in an improvement of 1.36\% on average compared to the original QC relaxation. 

By comparing the \mrn{sixth and ninth} columns of Table~\ref{tab:various_QC_results}, it can be seen that the proposed LRQC relaxation is superior to the RQC relaxation as it converges to tighter lower bounds for all test cases.
The \mrn{eighth} column in Table~\ref{tab:various_QC_results} lists the results for the LRQC relaxation with rotation angles computed by minimizing the volume of the envelope enclosing the function $V_lV_m\cos(\theta_{lm}-\delta_{lm}-\psi_l)\sin(\theta_{lm}-\delta_{lm}-\psi_l)$. Comparing the \mrn{eighth} column with the \mrn{third, fifth, and sixth} columns demonstrates that the LRQC relaxation improves the optimality gaps for all the test cases in Table~\ref{tab:various_QC_results} compared to the SOCP, QC and RQC relaxations, with some cases exhibiting substantial improvements. For instance, the RQC and LRQC relaxations have 0.63\% and 0.27\% optimality gaps for ``case3$\_$lmbd'' test case, respectively. 

\mrn{Moreover, the proposed LRQC relaxation finds better lower bounds for some test cases compared to the SDP relaxation. For instance, the SDP relaxation has 0.38\%, 4.99\%, and 2.54\% optimality gaps for the ``case3$\_$lmbd'', ``case3$\_$lmbd$\_\_$api'', and ``case24$\_$ieee$\_$rts$\_\_$sad'' test cases, respectively, while the proposed LRQC relaxation's optimality gaps are 0.26\%, 3.65\%, and 1.82\%. This indicates that the proposed LRQC relaxation can find better lower bounds for some problems while also being much faster than the SDP relaxation, as shown in Table~\ref{tab:various_QC_results_time}.}

We also observe that the proposed LRQC relaxation improves the optimality gaps for both small and large systems. For instance, the LRQC relaxation closes the optimality gap for the ``case2868$\_$rte$\_\_$api'' test case, where the previous RQC relaxation had an optimality gap of 0.16\%.
Moreover, the LRQCs relaxation with both the volume-minimizing $\psi_l$ and suggested $\psi^{\dagger}_{l}=85^\circ$ outperform the QC and RQC relaxations for all test cases. As expected from the analysis in Section~\ref{Best_Rotation_Angle}, applying the suggested $\psi^\dagger_{l}=85^\circ$ results in good performance across a variety of test cases.

% Our proposed LRQC approach demonstrates a significant improvement in OPF solution across multiple test cases. This improvement highlights the algorithm's substantial impact on cost savings within the power industry.} %Furthermore, it is crucial to note that various optimization problems in the power system, such as planning, transmission switching, unit commitment, and reliability, hinge on the effectiveness of the OPF solution where the proposed LRQC can play an important role.

\mrn{To provide further context for these results, we conducted a comparison between the optimality gap improvements achieved by the strongest previously known QC relaxation over the proceeding SOCP relaxation from~\cite{jabr2006radial}. For the PGLib-OPF test cases in Table~\ref{tab:various_QC_results}, this comparison shows optimality gap reductions from 0.0\% to 6.84\%, with an average across these test cases of 1.20\%. We note that these improvements were achieved via a number of advancements detailed in a series of papers including~\cite{coffrin2015qc,coffrin2016strengthen_tps,coffrin2016quadtrig}. In comparison, the LRQC relaxation's improvements over the previous state-of-the-art QC relaxation range from 0.0\% to 9.58\%, with an average of 1.31\% across the PGLib-OPF test cases in Table~\ref{tab:various_QC_results}.} \mrn{Thus, the tighter optimality gaps here are comparable in size to prior advances in QC relaxation formulations. We also note that optimality gap improvements of this size are meaningful given the large-scale nature of power systems. As an analogy, one might compare the considerable effort expended to close optimality gaps for mixed-integer linear programming solvers to within 0.5\% or 0.1\% in a variety of power systems applications like unit commitment.}

% \mrn{
% \begin{figure}
% \centering
% \subfloat[\mrn{Correlation of the optimality gap with the averaged distance to local optimality for the QC and LRQC relaxations.}]{\label{fig:average_distance}
% \begin{tikzpicture}
%     \node[] (image) at (0,0) {\includegraphics[width=0.75\columnwidth]{average_distance.eps}};
%        \node [anchor=west, font=\normalsize \rotatebox{90}{Average distance (\%)}] at (-3.9,0.2) {};
%        \node [anchor=west, font=\normalsize] at (-1.3,-2.8) {Optimality gap (\%)};
% \end{tikzpicture}
% }\\%\vspace{-.3cm} 
% \subfloat[\mrn{Correlation of the optimality gap with the cumulative constraint violation for the QC and LRQC relaxations.}]{\label{fig:cumulative_distance}
% \centering
% \centering
% \begin{tikzpicture}
%     \node[] (image) at (0,0) {\includegraphics[width=0.75\columnwidth]{cumulative_distance.eps}};
%        \node [anchor=west, font=\normalsize \rotatebox{90}{cumulative violation (\%)}] at (-3.9,0.2) {};
%        \node [anchor=west, font=\normalsize] at (-1.3,-2.8) {Optimality gap (\%)};
% \end{tikzpicture}
% }
% \caption{\mrn{Check the correlation of the optimality gap with the cumulative constraint violation and the averaged distance to local optimality for the QC and LRQR relaxations. Note that the axes are logarithmic. The red squares correspond to the LRQC relaxation, and the black circles correspond to the QC relaxation.}}
% \end{figure}}

\begin{table*}
\caption{Results from Applying Various Relaxations to Selected PGLib Test Cases}
%\vspace{-25pts}
\label{tab:various_QC_results}
\begin{center}\mrn{
\begin{tabular}{|c|c|c|c|c|c|c|c|c|}
\hline 
\multirow{2}{*}{Test Case} & \multirow{2}{*}{AC} & \multirow{2}{*}{SOCP Gap} & \multirow{2}{*}{SDP Gap} & \multirow{2}{*}{QC Gap} & \multicolumn{2}{c|}{{\thead{RQC $(\psi^*)$ }}}  & \multirow{2}{*}{{\thead{LRQC}}} & \multirow{2}{*}{{\thead{LRQC Gap (\%) \\$(\psi^\dagger_{l}=85^\circ)$ }}}
\tabularnewline
\cline{6-7} \cline{7-7} 
 & \thead{(\$/hr)} &\thead{Gap (\%)}  & \thead{Gap (\%)} &\thead{Gap (\%)}  & \thead{Gap (\%)} & \thead{$\psi^\ast$} & \thead{Gap (\%)} & \tabularnewline
\hline 
case3\_lmbd & 5812.64 & 1.32 & \mrn{0.39} & 0.97 & 0.63 & 11 & 0.26 & 0.27\tabularnewline
\hline 
case14\_ieee & 2178.08 & 0.11 & \mrn{0.00} & 0.11 & 0.10 & -23 & 0.09 & 0.10 \tabularnewline
\hline 
case30\_ieee & 8208.52 & 18.84 & \mrn{0.00} & 18.67 & 11.82 & -25 & 9.08 & 12.06\tabularnewline
\hline 
case39\_epri & 138415.56 & 0.55 & \mrn{0.01} & 0.54 & 0.51 & 0 & 0.50 & 0.51\tabularnewline
\hline 
case89\_pegase & 107285.67 & 0.75 & \mrn{0.30} & 0.75 & 0.74 & 77 & 0.73 & 0.74\tabularnewline
\hline 
case118\_ieee & 97213.61 & 0.90 & \mrn{0.07} & 0.77 & 0.62 & 70 & 0.55 & 0.56\tabularnewline
\hline 
case240\_pserc & 3329670.06 & 2.77 & \mrn{1.43} & 2.72 & 2.54 & 8 & 2.39 & 2.41\tabularnewline
\hline 
case300\_ieee & 565219.97 & 2.62 & \mrn{0.12} & 2.56 & 2.24 & -13 & 2.18 & 2.16\tabularnewline
\hline 
case1951\_rte & 2085581.84 & 0.13 & \mrn{0.01} & 0.13 & 0.11 & -10 & 0.11 & 0.11\tabularnewline
\hline 
case2316\_sdet & 1775325.55 & 1.79 & \mrn{0.66} & 1.79 & 1.78 & -9 & 1.76 & 1.77\tabularnewline
\hline 
case2848\_rte & 1286608.19 & 0.12 & \mrn{0.05} & 0.12 & 0.12 & -48 & 0.11 & 0.11\tabularnewline
\hline 
case2869\_pegase & 2462790.43 & 1.01 & \mrn{0.08}& 1.00 & 0.98 & -10 & 0.98 & 0.98\tabularnewline
\hline 
case6515\_rte & 2825499.64 & 6.40 & \mrn{5.57} & 6.39 & 6.38 & 82 & 6.37 & 6.37\tabularnewline
\hline 
case9241\_pegase & 6243090.38
 & 2.54 & \mrn{2.10} & 1.71 & 1.69 & -10 & 1.66 & 1.67\tabularnewline
\hline 
case3\_lmbd\_\_api & 11242.12 & 9.32 & \mrn{7.34} & 4.57 & 3.93 & -71 & 3.65 & 3.68\tabularnewline
\hline 
case14\_ieee\_\_api & 5999.36 & 5.13 & \mrn{0.00} & 5.13 & 5.13 & 63 & 5.13 & 5.13\tabularnewline
\hline 
case24\_ieee\_rts\_\_api & 134948.17 & 17.87 & \mrn{2.06} & 11.02 & 6.98 & -11 & 4.47 & 6.15\tabularnewline
\hline 
case30\_fsr\_\_api & 701.15 & 2.76 & \mrn{0.28} & 2.75 & 2.69 & 78 & 2.58 & 2.63\tabularnewline
\hline 
case30\_ieee\_\_api & 18043.92 & 5.45 & \mrn{0.00} & 5.45 & 5.29 & -23 & 4.50 & 5.25\tabularnewline
\hline 
case73\_ieee\_rts\_\_api & 422726.14 & 12.88 & \mrn{2.91} & 9.54 & 7.24 & -10 & 6.21 & 6.92\tabularnewline
\hline 
case118\_ieee\_\_api & 242054.0 & 28.81 & \mrn{11.16}
 & 28.67 & 26.38 & -8 & 24.17 & 26.00\tabularnewline
\hline 
case162\_ieee\_dtc\_\_api & 120996.12 & 4.36 & \mrn{1.42} & 4.32 & 4.27 & -9 & 4.24 & 4.26\tabularnewline
\hline 
case179\_goc\_\_api & 1932120.33 & 9.88 & \mrn{0.55} & 5.86 & 4.06 & -78 & 3.16 & 3.16\tabularnewline
\hline 
case300\_ieee\_\_api & 650147.21 & 0.89 & \mrn{0.08} & 0.83 & 0.70 & -15 & 0.64 & 0.64\tabularnewline
\hline 
case2848\_rte\_\_api & 1496368.95 & 0.22 & \mrn{0.06} & 0.22 & 0.21 & 79 & 0.18 & 0.20\tabularnewline
\hline 
case2869\_pegase\_\_api & 2934160.71 & 1.33 & \mrn{0.45} & 1.32 & 1.30 & -10 & 1.11 & 1.29\tabularnewline
\hline 
case6515\_rte\_\_api & 3162434.34
 & 1.95 & \mrn{1.18} & 1.91 & 1.91 & -8 & 1.90 & 1.91\tabularnewline
\hline 
case3\_lmbd\_\_sad & 5959.33 & 3.74 & \mrn{1.86} & 1.38 & 1.02 & 68 & 0.92 & 0.92\tabularnewline
\hline 
case14\_ieee\_\_sad & 2777.35 & 21.54 & \mrn{0.09} & 19.16 & 15.39 & -12 & 12.70 & 14.20\tabularnewline
\hline 
case24\_ieee\_rts\_\_sad & 76943.24 & 9.55 & \mrn{4.36} & 2.74 & 2.12 & -12 & 1.82 & 1.84\tabularnewline
\hline 
case30\_ieee\_\_sad & 8208.52 & 9.69 & \mrn{0.00} & 5.66 & 4.45 & 66 & 3.94 & 4.11\tabularnewline
\hline 
case39\_epri\_\_sad & 148354.41 & 0.66 & \mrn{0.02} & 0.20 & 0.17 & 82 & 0.16 & 0.16\tabularnewline
\hline 
case57\_ieee\_\_sad & 38663.88 & 0.70 & \mrn{0.05} & 0.32 & 0.31 & 84 & 0.29 & 0.29\tabularnewline
\hline 
case73\_ieee\_rts\_\_sad & 227745.73 & 6.74 & \mrn{2.75} & 2.37 & 1.82 & 78 & 1.55 & 1.56\tabularnewline
\hline 
case118\_ieee\_\_sad & 103292.3082 & 8.21 & \mrn{1.43} & 6.67 & 5.07 & 69 & 4.00 & 4.48\tabularnewline
\hline 
case162\_ieee\_dtc\_\_sad & 108695.95 & 6.48 & \mrn{2.08} & 6.22 & 5.54 & 76 & 4.51 & 5.02\tabularnewline
\hline 
case300\_ieee\_\_sad & 565712.83 & 2.60 & \mrn{0.14} & 2.34 & 1.59 & 83 & 1.32 & 1.37\tabularnewline
\hline 
case1951\_rte\_\_sad & 2092788.97 & 0.48 & \mrn{0.30} & 0.43 & 0.42 & -10 & 0.40 & 0.40\tabularnewline
\hline 
case2746wop\_k\_\_sad & 1234338.04 & 2.36 & \mrn{0.71} & 1.99 & 1.84 & 80 & 1.77 & 1.76\tabularnewline
\hline 
case6515\_rte\_\_sad & 2882577.97
 & 8.26 & \mrn{6.21} & 8.22 & 8.21 & 75 & 8.15 & 8.17
\tabularnewline
\hline 
case9241\_pegase\_\_sad & 6319549.55
 & 2.48 & \mrn{2.77} & 2.42 & 2.40 & 82 & 2.39 & 2.39\tabularnewline
\hline 
\end{tabular}}
{\\ \raggedright \emph{QC Gap}: Optimality gap for the QC relaxation from~\eqref{eq:QC relaxation},  \emph{RQC Gap}: Optimality gap for the relaxation from~\cite[Eq.~(16)]{NarimaniTPS1}, 
\mrn{\emph{SOCP Gap}: Optimality gap for the SOCP relaxation from~\cite{jabr2006radial},
\emph{SDP Gap}: \mrn{Optimality gap for the SDP relaxation from~\cite{molzahn2013sdp},}}
\emph{LRQC Gap}: Optimality gap for the relaxation from~\eqref{eq:LRQC OPF}, $\psi^*$:~Use of the volume-minimizing $\psi_l$ for this case. \par}
\end{center}
\end{table*}

\begin{table*}
 \caption{Execution Time from Applying Various Relaxations to Selected PGLib Test Cases}
%\vspace{-25pts}
\label{tab:various_QC_results_time}
\begin{center}\mrn{
\begin{tabular}{|c|c|c|c|c|c|c|}
\hline Test Cases
& AC Time & SOCP Time & SDP Time & QC Time & RQC Time  & LRQC Time\tabularnewline
\hline 
\hline 
case3\_lmbd &0.72 &  0.09 & \mrn{0.06} & 0.26 & 0.01 & 0.02\tabularnewline
\hline 
case14\_ieee &0.01  & 0.04 & \mrn{0.06} & 0.37 & 0.02 & 0.05\tabularnewline
\hline 
case30\_ieee &0.03  & 0.12 & \mrn{0.13} & 0.33 & 0.03 & 0.12\tabularnewline
\hline 
case39\_epri &0.05  & 0.16 & \mrn{0.14} & 0.38 & 0.06 & 0.12\tabularnewline
\hline 
case89\_pegase &0.23  & 0.45 & \mrn{1.09} & 0.87 & 0.44 & 1.10\tabularnewline
\hline 
case118\_ieee &0.16  & 0.32 & \mrn{0.59} & 0.55 & 0.23 & 0.37\tabularnewline
\hline 
case240\_pserc &4.37  & 1.06 & \mrn{1.07} & 1.13 & 0.92 & 1.42\tabularnewline
\hline 
case300\_ieee &1.22   & 0.94 & \mrn{1.57} & 1.54 & 3.15 & 1.33\tabularnewline
\hline 
case1951\_rte &8.07  & 6.56 & \mrn{17.84} & 7.70 & 8.97 & 12.65\tabularnewline
\hline 
case2316\_sdet &5.2  & 5.68 & \mrn{59.21} & 6.11 & 7.37 & 7.36\tabularnewline
\hline 
case2848\_rte &11.77  & 8.04 & \mrn{32.26} & 10.56 & 10.63 & 18.04\tabularnewline
\hline 
%case2868\_rte &13.63  & 8.29 & - & 13.79 & 13.15 & 23.47\tabularnewline
%\hline 
case2869\_pegase &11.29  & 9.16 & \mrn{34.96} & 15.54 & 13.17 & 19.88\tabularnewline
\hline 
case6515\_rte & 71.92
  & 23.72 & \mrn{224.47} & 41.13 & 17.62 & 49.74\tabularnewline
\hline 
case9241\_pegase &82.83  & 34.82 & \mrn{394.62} & 101.47 & 120.31 & 138.52\tabularnewline
\hline 
case3\_lmbd\_\_api &0.02  & 0.09 & \mrn{0.01} & 0.51 & 0.01 & 0.02\tabularnewline
\hline 
case14\_ieee\_\_api &0.02  & 0.04 & \mrn{0.06} & 0.37 & 0.03 & 0.05\tabularnewline
\hline 
case24\_ieee\_rts\_\_api &0.11  & 0.19 & \mrn{0.14} & 0.71 & 0.04 & 0.04\tabularnewline
\hline 
case30\_ieee\_\_api &0.03  & 0.12 & \mrn{0.14} & 0.31 & 0.06 & 0.08\tabularnewline
\hline 
case30\_fsr\_\_api &0.04  & 0.08 & \mrn{0.14} & 0.37 & 0.14 & 0.12\tabularnewline
\hline 
case73\_ieee\_rts\_\_api &0.24  & 0.28 & \mrn{0.44} & 1.00 & 0.37 & 0.43\tabularnewline
\hline 
case118\_ieee\_\_api &0.24  & 0.34 & \mrn{0.79} & 0.53 & 0.97 & 0.41\tabularnewline
\hline 
pcase162\_ieee\_dtc\_\_api &0.28  & 0.46 & \mrn{2.42} & 0.68 & 0.44 & 0.83\tabularnewline
\hline 
pcase179\_goc\_\_api &1.08  & 1.53 & \mrn{0.83} & 0.82 & 0.64 & 1.51\tabularnewline
\hline 
case300\_ieee\_\_api &0.90  & 0.91 & \mrn{2.24} & 1.22 & 2.36 & 1.56\tabularnewline
\hline 
case2848\_rte\_\_api &22.61  & 7.51 & \mrn{32.84} & 12.60 & 12.56 & 18.52\tabularnewline
\hline 
%case2868\_rte\_\_api &15.74  & 7.73 & - & 12.52 & 13.29 & 18.63\tabularnewline
%\hline 
case2869\_pegase\_\_api &11.48  & 8.56 & \mrn{39.98} & 10.72 & 11.00 & 17.26\tabularnewline
\hline 
case6515\_rte\_\_api &81.42  & 22.49 & \mrn{254.93} & 45.08 & 43.36 & 6.72\tabularnewline
\hline 
case3\_lmbd\_\_sad &0.01  & 0.12 & \mrn{0.01} & 0.44&  0.01 & 0.02\tabularnewline
\hline 
case14\_ieee\_\_sad &0.02
 & 0.04 & \mrn{0.04} & 0.35  & 0.03 & 0.05\tabularnewline
\hline 
case24\_ieee\_rts\_\_sad &0.12  & 0.09 & \mrn{0.11} & 0.40 & 0.06 & 0.06\tabularnewline
\hline 
case30\_ieee\_\_sad &0.03  & 0.07 & \mrn{0.08} & 0.32 & 0.06 & 0.13\tabularnewline
\hline 
case39\_epri\_\_sad &0.04  & 0.14 & \mrn{0.12} & 0.36 & 0.11 & 0.12\tabularnewline
\hline 
case57\_ieee\_\_sad &0.06  & 0.14 & \mrn{0.21} & 0.38 & 0.11 & 0.16\tabularnewline
\hline 
case73\_ieee\_rts\_\_sad &0.13  & 0.21 & \mrn{0.38} & 0.41 & 0.41 & 0.44\tabularnewline
\hline 
case118\_ieee\_\_sad &0.18  & 0.28 & \mrn{0.62} & 0.58 & 0.39 & 0.68\tabularnewline
\hline 
case162\_ieee\_dtc\_\_sad &0.43  & 0.46 & \mrn{2.09} & 0.86 & 0.84 & 0.91\tabularnewline
\hline 
case300\_ieee\_\_sad &0.52  & 1.16 & \mrn{1.55} & 1.94 & 2.06 & 1.94\tabularnewline
\hline 
case1951\_rte\_\_sad &7.74  & 6.03 & \mrn{20.58} & 7.89 & 8.45 & 17.25\tabularnewline
\hline 
case2746wop\_k\_\_sad &5.98  & 4.38 & \mrn{67.03}& 6.90 & 7.91 & 10.47\tabularnewline
\hline 
case6515\_rte\_\_sad &71.38  & 22.79 & \mrn{233.99} & 46.16 & 52.88 & 103.51	\tabularnewline
\hline
case9241\_pegase\_\_sad &129.78   & 34.82 & \mrn{478.62} & 89.37 & 85.27 & 196.28	\tabularnewline
\hline 
\end{tabular}}
\end{center}
\end{table*}

As shown in \mrn{Table~\ref{tab:various_QC_results_time}}, the LRQC relaxation's improved tightness comes at the cost of slower (but still tractable) computational times for some test cases. 
When analyzing the last two columns of \mrn{Table~\ref{tab:various_QC_results_time}}, it becomes evident that the impact of adding the proposed envelopes on execution time is quite diverse. For instance, in cases like ``case300$\_$ieee,'' implementing these envelopes leads to a reduction in execution time by over \mrn{$57.7$\%}. On the contrary, for other test cases, there is a considerable increase in execution time, reaching up to \mrn{$200$\%} in some instances. However, on average across all the test cases, enforcing the proposed envelopes results in a moderate increase of less than $38$\% in the time required to solve the RQC relaxation from~\cite{NarimaniTPS1}.

\subsection{Assessing Decision Variable Quality}
\label{sec:results_quality}

% \mrn{Note that our proposed LRQC relaxation stands out in terms of tractability when compared to existing SDP relaxation methods. The complexity of SDP relaxation methods tends to exponentially increase with the size of the network, rendering them impractical for large test systems with thousands of nodes. In contrast, our proposed algorithm demonstrates remarkable efficiency by easily providing a lower bound for very large OPF test cases with more than ten thousand nodes in a matter of minutes.}

\mrn{The optimality gap is a key measure of relaxation tightness that is both most commonly used to benchmark the performance of various relaxations and, as discussed in the introduction, is highly relevant for many applications. For other applications, one may also be interested in the quality of the decision variable values for voltage phasors, line flows, generator outputs, etc. in a relaxation's solution. To assess this, reference~\cite{venzke2020inexact} proposes two metrics for gauging the proximity to local optimality and to AC feasibility. A main finding of~\cite{venzke2020inexact} is that many power flow relaxations have a nonlinear relationship between the optimality gap and these metrics. Small optimality gaps typically correspond to small values of these metrics, but moderate to large optimality gaps may have either small or large values of these metrics. Our proposed LRQC relaxation exhibits similar behavior.}

\mrn{For the proximity to local optimality metric (``Average Normalized Distance to a Local Solution'') in~\cite{venzke2020inexact}, we note that our proposed LRQC method outperforms the QC relaxation in 67\% of PGLib-OPF test cases with an average improvement of 12\%. Similarly, for the AC feasibility metric (``Cumulative Normalized Constraint Violation'') in~\cite{venzke2020inexact}, the LRQC relaxation outperforms the QC relaxation in 58\% of PGLib-OPF test cases with an average improvement of 7\%. Thus, our proposed LRQC relaxation often outperforms the original QC relaxation on both metrics in~\cite{venzke2020inexact}.}

\mrn{We also note that recent work in~\cite{taheri2024} proposes a new solution restoration method that significantly improves upon the simplistic AC power flow method used in~\cite{venzke2020inexact}, often resulting in several orders-of-magnitude improvements in the accuracy of the restored AC power flow feasible solutions. We therefore focus on the optimality gaps as our primary metric for comparing relaxations.}

\subsection{Impacts of Bound Tightening on Optimality Gaps}
\label{sec:results_OBBT}
\mrn{As key parameters in forming the convex envelopes for the trigonometric functions, the accuracy of QC relaxations strongly depend on the tightness of the bounds on voltage magnitudes and phase angle differences. To characterize this, we applied the bound tightening method described in~\cite{dmitry2019} to several selected test cases and then executed both the original QC and the proposed LRQC relaxations on the tightened test cases. As expected, the results indicate that bound tightening has a substantial impact on the optimality gaps for both relaxations. Comparing the optimality gaps for both QC and LRQC relaxations in Table~\ref{tab:results_for_tightened_cases} with their corresponding values in Table~\ref{tab:various_QC_results} reveals that applying bound tightening reduces the optimality gaps. \mrn{For instance, applying the bound tightening approach reduces the gaps for the ``case39$\_$epri'' and ``case118 ieee$\_\_$api'' cases by 0.31\% and 11.53\%, respectively.} %This reinforces the fact that the accuracy of the QC and LRQC relaxations strongly depend on the tightness of the bounds. Note that the proposed LRQC relaxation finds better lower bounds for tightened test cases compared to the original QC relaxation, which proves its superiority over the original QC relaxation.
We emphasize that the proposed LRQC relaxation still finds better lower bounds for the bound-tightened test cases compared to the original QC relaxation, again demonstrating its superiority over the original QC relaxation.}

\subsection{Balancing Execution time and LRQC Relaxation Tightness}
\label{sec:results_Accurracy_Time}
\mrn{The parameter $N_{seg}$ plays an important role in determining the tightness of the proposed LRQC relaxation. To demonstrate its impacts on both solution time and tightness, we applied the LRQC relaxation with various $N_{seg}$ values to selected test cases. This parameter's influences on the optimality gap and execution time of the LRQC relaxation are presented in Tables~\ref{tab:optimality_gap_Nseg} and~\ref{tab:execution_time_Nseg}, respectively. The results in these tables indicate that increasing $N_{seg}$ beyond five typically has limited impacts on the optimality gaps but can cause significant increases to the solution times. Supporting the analytical assessment in Section~\ref{intersection_volume}, these empirical results suggest that $N_{seg} = 5$ provides a good balance between tightness and tractability of the LRQC relaxation. Note that with $N_{seg}=5$, the proposed LRQC relaxation finds tighter lower bounds for all test cases compared to the original QC, RQC, and SOCP relaxations.}

\begin{table}
  \caption{\centering{Optimality Gaps (\%) of the QC and LRQC Relaxations for Selected Tightened PGLib Test Cases}}
\label{tab:results_for_tightened_cases}
    \centering
\mrn{
\begin{tabular}{|c|c|c|c|}
\hline 
Test Cases & AC (\$/hr) & QC & LRQC\tabularnewline
\hline 
\hline 
case3\_lmbd & 5812.64 & 0.8 & 0.26\tabularnewline
\hline 
case39\_epri & 138415.5627 & 0.3 & 0.20\tabularnewline
\hline 
case118\_ieee & 97213.61 & 0.4 & 0.37\tabularnewline
\hline 
case240\_pserc & 3329670.06 & 2.5 & 2.13\tabularnewline
\hline 
case300\_ieee & 565219.97 & 1.5 & 1.15\tabularnewline
\hline 
case3\_lmbd\_\_api & 11242.12 & 3.9 & 3.32\tabularnewline
\hline 
case30\_fsr\_\_api & 701.15 & 2.5 & 2.40\tabularnewline
\hline 
case73\_ieee\_rts\_\_api & 422726.14 & 5.5 & 4.14\tabularnewline
\hline 
case118\_ieee\_\_api & 242054.01 & 20.8 & 14.47\tabularnewline
\hline 
case162\_ieee\_dtc\_\_api & 120996.09 & 4.1 & 3.66\tabularnewline
\hline 
case300\_ieee\_\_api & 650147.21 & 0.6 & 0.42\tabularnewline
\hline 
case3\_lmbd\_\_sad & 5959.33 & 1.4 & 0.94\tabularnewline
\hline 
case24\_ieee\_rts\_\_sad & 76943.24 & 1.4 & 1.12\tabularnewline
\hline 
case39\_epri\_\_sad & 148354.41 & 0.1 & 0.12\tabularnewline
\hline 
case57\_ieee\_\_sad & 38663.88 & 0.3 & 0.25\tabularnewline
\hline 
case73\_ieee\_rts\_\_sad & 227745.73 & 1.7 & 1.25\tabularnewline
\hline 
case118\_ieee\_\_sad & 105216.67 & 5.9 & 3.70\tabularnewline
\hline 
case162\_ieee\_dtc\_\_sad & 108695.95 & 5.4 & 3.94\tabularnewline
\hline 
case300\_ieee\_\_sad & 565712.83 & 1.4 & 1.14\tabularnewline
\hline 
\end{tabular}
}
\end{table}

\begin{table}
\caption{Optimality Gaps (\%) of the LRQC Relaxation for Selected PGLib Test Cases for Differing Numbers of Segments}
     \label{tab:optimality_gap_Nseg}
     \centering
 \begin{scriptsize}
    \mrn{\begin{tabular}{|c|c|c|c|c|}
\hline 
\multirow{2}{*}{Test Cases} & \multicolumn{4}{c|}{Number of segments ($N_{seg}$) }\tabularnewline
\cline{2-5} \cline{3-5} \cline{4-5} \cline{5-5} 
 &   $N_{seg}=3$  & $N_{seg}=5$ & $N_{seg}=10$ & $N_{seg}=20$\tabularnewline
\hline 
case3\_lmbd & 0.40 & 0.27 & 0.23 & 0.21\tabularnewline
\hline 
case14\_ieee & 0.10 & 0.10 & 0.10 & 0.07\tabularnewline
\hline 
case24\_ieee\_rts & 0.01 & 0.01 & 0.01 & 0.01\tabularnewline
\hline 
case30\_ieee & 16.48 & 12.06 & 8.46 & 7.89\tabularnewline
\hline 
case118\_ieee & 0.62 & 0.56 & 0.48 & 0.45\tabularnewline
\hline 
case300\_ieee & 2.22 & 2.16 & 1.67 & 1.60\tabularnewline
\hline 
case3\_lmbd\_\_api & 3.64 & 3.68 & 3.59 & 3.58\tabularnewline
\hline 
case14\_ieee\_\_api & 5.13 & 5.13 & 5.13 & 5.11\tabularnewline
\hline 
case24\_ieee\_rts\_\_api & 9.25 & 6.15 & 5.48 & 5.07\tabularnewline
\hline 
case30\_ieee\_\_api & 5.45 & 5.25 & 4.26 & 3.87\tabularnewline
\hline 
case118\_ieee\_\_api & 27.28 & 26.00 & 23.64 & 22.99\tabularnewline
\hline 
case300\_ieee\_\_api & 0.69 & 0.64 & 0.52 & 0.49\tabularnewline
\hline 
case3\_lmbd\_\_sad & 0.91 & 0.92 & 0.91 & 0.91\tabularnewline
\hline 
case14\_ieee\_\_sad & 15.13 & 14.20 & 13.66 & 13.57\tabularnewline
\hline 
case24\_ieee\_rts\_\_sad & 1.91 & 1.84 & 1.83 & 1.82\tabularnewline
\hline 
case30\_ieee\_\_sad & 4.36 & 4.11 & 4.11 & 4.09\tabularnewline
\hline 
case118\_ieee\_\_sad & 4.80 & 4.48 & 4.27 & 4.23\tabularnewline
\hline 
case300\_ieee\_\_sad & 1.40 & 1.37 & 1.21 & 1.19\tabularnewline
\hline 
\end{tabular}}
\end{scriptsize}
\end{table}

\begin{table}
  \caption{Execution Times of the LRQC Relaxation for Selected PGLib Test Cases for Differing Numbers of Segments}
 \label{tab:execution_time_Nseg}
    \centering
     \begin{scriptsize}\mrn{
   % Preview source code for paragraph 0
\begin{tabular}{|c|c|c|c|c|}
\hline 
\multirow{2}{*}{Test Cases} & \multicolumn{4}{c|}{Number of segments $N_{seg}$}\tabularnewline
\cline{2-5} \cline{3-5} \cline{4-5} \cline{5-5} 
 & $N_{seg}=3$ & $N_{seg}=5$ & $N_{seg}=10$ & $N_{seg}=20$\tabularnewline
\hline 
case3\_lmbd & 0.02 & 0.02 & 0.05 & 0.06\tabularnewline
\hline 
case14\_ieee & 0.05 & 0.05 & 0.27 & 0.60\tabularnewline
\hline 
case24\_ieee\_rts & 0.12 & 0.14 & 0.35 & 0.72\tabularnewline
\hline 
case30\_ieee & 0.12 & 0.12 & 0.65 & 1.62\tabularnewline
\hline 
case118\_ieee & 0.32 & 0.37 & 1.65 & 3.09\tabularnewline
\hline 
case300\_ieee & 1.26 & 1.33 & 13.48 & 24.43\tabularnewline
\hline 
case3\_lmbd\_\_api & 0.02 & 0.02 & 0.07 & 0.13\tabularnewline
\hline 
case14\_ieee\_\_api & 0.05 & 0.05 & 0.18 & 0.53\tabularnewline
\hline 
case24\_ieee\_rts\_\_api & 0.04 & 0.04 & 0.43 & 1.63\tabularnewline
\hline 
case30\_ieee\_\_api & 0.07 & 0.08 & 0.43 & 0.84\tabularnewline
\hline 
case118\_ieee\_\_api & 0.37 & 0.41 & 2.36 & 4.30\tabularnewline
\hline 
case300\_ieee\_\_api & 1.50 & 1.56 & 17.77 & 44.40\tabularnewline
\hline 
case3\_lmbd\_\_sad & 0.02 & 0.02 & 0.03 & 0.08\tabularnewline
\hline 
case14\_ieee\_\_sad & 0.05 & 0.05 & 0.20 & 1.02\tabularnewline
\hline 
case24\_ieee\_rts\_\_sad & 0.05 & 0.06 & 0.32 & 0.61\tabularnewline
\hline 
case30\_ieee\_\_sad & 0.12 & 0.13 & 0.47 & 1.09\tabularnewline
\hline 
case118\_ieee\_\_sad & 0.60 & 0.68 & 7.28 & 3.52\tabularnewline
\hline 
case300\_ieee\_\_sad & 1.90 & 1.94 & 7.18 & 17.79\tabularnewline
\hline 
\end{tabular}}
\end{scriptsize}
\end{table}

\section{Conclusion}
\label{conclusion}
This paper has proposed tighter envelopes for the product and trigonometric terms in the power flow equations to improve the tightness of the QC relaxation. These envelopes are developed by considering a particular nonlinear function whose projections are the expressions appearing in the power flow equations. 
Additionally, we exploit characteristics of the sine and cosine expressions along with the changes in their curvature to tighten convex envelopes associated with the trigonometric terms. Comparison to a state-of-the-art RQC relaxation implementation demonstrates the value of these improvements via reduced optimality gaps on challenging test cases while maintaining computational tractability.

\ifarxiv
\appendix

\label{apx:enclosing_points}

\begin{figure}
\def\tabularxcolumn#1{m{#1}}
\subfloat[$-165^\circ \leq \theta_{lm}-\delta_{lm}-\psi_l \leq 15^\circ$\label{fig:orig1}]{
\begin{tabular}{cc}
\small $\cos(\theta_{lm}-\delta_{lm}-\psi_l)$  & \small $F(\theta_{lm},\delta_{lm},\psi_l)$ \\
\includegraphics[width=4.0cm]{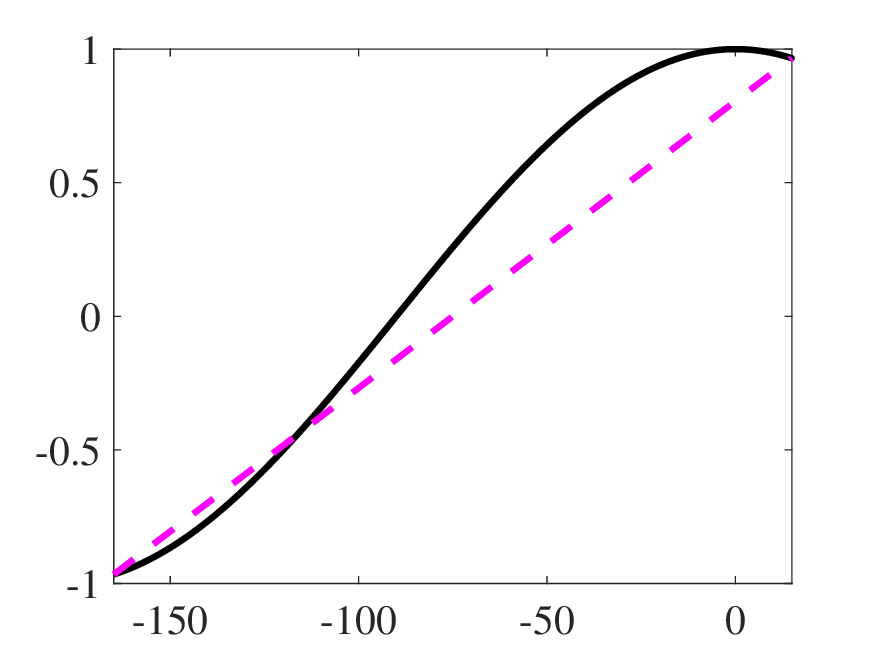} & \includegraphics[width=4.0cm]{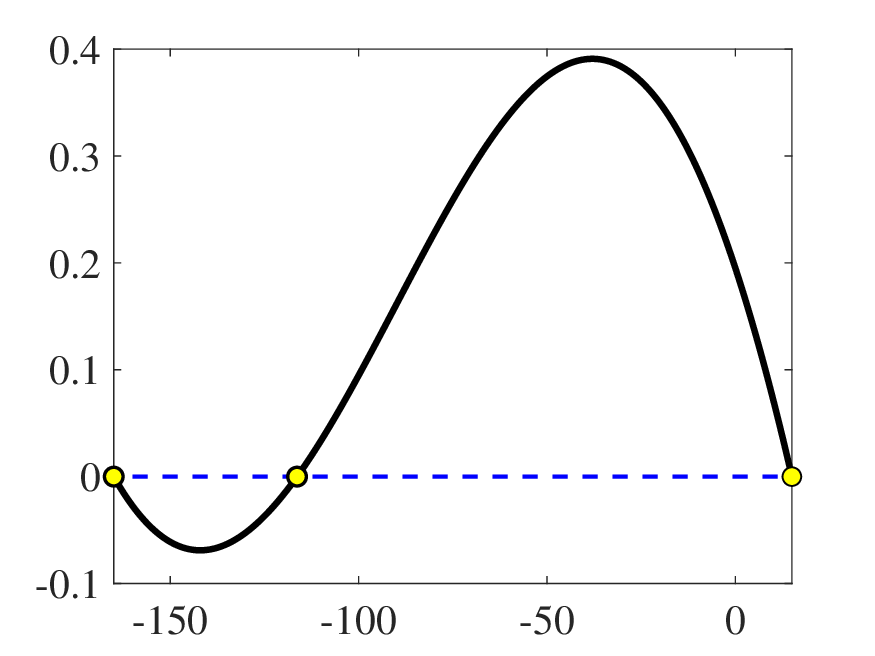}
\end{tabular}%
}\\
\subfloat[$-105^\circ \leq \theta_{lm}-\delta_{lm}-\psi_l \leq 75^\circ$\label{fig:orig3}]{
\begin{tabular}{cc}
\small $\cos(\theta_{lm} - \delta_{lm}-\psi_l)$  & \small $F(\theta_{lm},\delta_{lm},\psi_l)$  \\
\includegraphics[width=4.0cm]{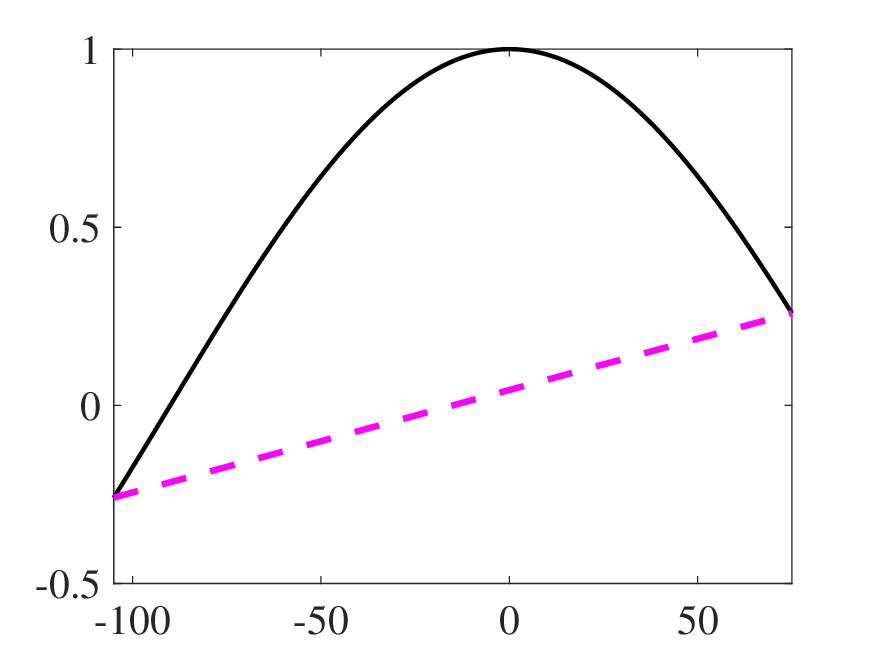} & \includegraphics[width=4.0cm]{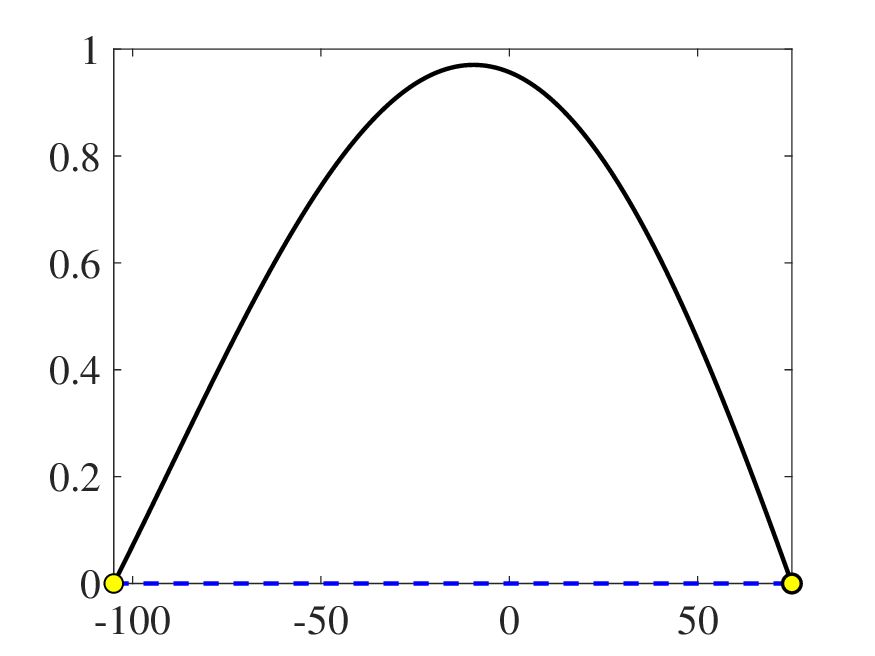}
\end{tabular}%
}
\centering
\caption{The figures in the left column show visualizations of the function $\cos(\theta_{lm} - \delta_{lm} - \psi_l)$ (black curve) and the line connecting the endpoints of this function at $\underline{\theta}_{lm}$ and $\overline{\theta}_{lm}$ (dashed magenta line) for different values of $\delta_{lm}$, $\underline{\theta}_{lm}$, and $\overline{\theta}_{lm}$. The figures in the right column show the corresponding function $F(\theta_{lm})$ and its roots between $\underline{\theta}_{lm}-\delta_{lm}-\psi_l$ and $\overline{\theta}_{lm}-\delta_{lm}-\psi_l$.}
\label{fig:proof}  %\vspace{-0.5cm}
\end{figure}

This appendix details the process of selecting appropriate tangent lines to the sine and cosine functions for arguments that take a general range of values in order to construct valid convex envelopes. While finding tangent lines for sine and cosine functions at any point is relatively straightforward --- one can simply take the derivative to obtain the slope and align the line with the tangent point --- the real challenge lies in identifying tangent lines that do not intersect the sine or cosine function at any other points.
This occurs due to the possible variation in curvature of the sine and cosine functions for arguments within a general range of values. If the curvature's sign does not change, tangent lines will not intersect the sine and cosine function at any other points.

The case that poses a challenge is when the curvature's sign changes, leading to the possibility of a tangent line intersecting the function at another point. Hence, the initial step in selecting appropriate tangent lines is to examine whether the curvature's sign changes within the specified range for the trigonometric function's argument.
To tackle this issue, we introduce an auxiliary function, denoted as $F(\theta_{lm})$ in~\eqref{eq:F}, which helps us ascertain whether the curvature of the trigonometric function changes within the range  $[\underline{\theta}_{lm}-\delta_{lm}-\psi_l$, $\overline{\theta}_{lm}-\delta_{lm}-\psi_l]$. 
The function $F(\theta_{lm})$ captures the difference between the trigonometric function $\cos(\theta_{lm} - \delta_{lm}-\psi_l)$ and the straight line connecting its endpoints at $\underline{\theta}_{lm}$ and $\overline{\theta}_{lm}$:
\begin{align}
\label{eq:F}
\nonumber \noindent &F(\theta_{lm})  = \cos(\theta_{lm} - \delta_{lm}-\psi_l) - \cos(\overline{\theta}_{lm}-\delta_{lm}-\psi_l)\\
\nonumber &\qquad\quad-\frac{\cos(\overline{\theta}_{lm}-\delta_{lm}-\psi_l)-\cos(\underline{\theta}_{lm}-\delta_{lm}-\psi_l)}{\overline{\theta}_{lm}-\underline{\theta}_{lm}}\\&\qquad\quad\times\left(\theta_{lm}-\overline{\theta}_{lm}-\psi_l\right).
\end{align}

The set of zeros of the derivative of $F(\theta_{lm})$, i.e., the set of solutions to $\frac{dF(\theta_{lm})}{d \theta_{lm}} = 0$, is a key quantity when determining if the curvature of $\cos(\theta_{lm}-\delta_{lm}-\psi_l)$ changes between $\underline{\theta}_{lm}-\delta_{lm}-\psi_l$ and $\overline{\theta}_{lm}-\delta_{lm}-\psi_l$.
In Fig.~\ref{fig:proof}, on the left side, we present illustrative examples of the function $\cos(\theta_{lm} - \delta_{lm}-\psi_l)$ (depicted by the black curve) along with the dashed magenta line representing the connection between its endpoints at $\underline{\theta}_{lm}-\delta_{lm}-\psi_l$ and $\overline{\theta}_{lm}-\delta_{lm}-\psi_l$. The right side visualizes the function $F(\theta_{lm})$ itself, with its roots indicated by yellow circles.

The derivative of $F(\theta_{lm})$ is
\begin{align}
\nonumber & \frac{dF(\theta_{lm})}{d \theta_{lm}} = -\sin(\theta_{lm}-\delta_{lm}-\psi_l)
 \\ 
& \quad - \frac{\cos(\overline{\theta}_{lm} -\delta_{lm}-\psi_l)-\cos(\underline{\theta}_{lm}-\delta_{lm}-\psi_l)}{\overline{\theta}_{lm}-\underline{\theta}_{lm}}.
\end{align}
  We denote the set of zeroes for $\frac{dF(\theta_{lm})}{d \theta_{lm}} = 0$ by $\mathcal{Z}_{\underline{\theta}_{lm}, \overline{\theta}_{lm}, \delta_{lm},\psi_l}$, where the subscripts indicate that the set is parameterized by $\underline{\theta}_{lm}$, $\overline{\theta}_{lm}$, $\delta_{lm}$, and $\psi_l$:
\begin{align}
\label{eq:roots}
\nonumber & \mathcal{Z}_{\underline{\theta}_{lm}, \overline{\theta}_{lm}, \delta_{lm}, \psi_l} = \\ 
\nonumber & \Big\lbrace (-1)^\kappa \arcsin\left(\frac{\cos(\underline{\theta}_{lm}-\delta_{lm}) - \cos(\overline{\theta}_{lm}-\delta_{lm})}{\left(\overline{\theta}_{lm}-\underline{\theta}_{lm}\right)}\right) + \pi\kappa,\\
& \qquad \kappa = \ldots, -3, -2, -1, 0, 1, 2, 3, \ldots \Big\rbrace.
\end{align}
Let $\left|\,\cdot\,\right|$ denote the cardinality of a set. The cardinality of $\mathcal{Z}_{\underline{\theta}_{lm}, \overline{\theta}_{lm}, \delta_{lm},\psi_l}$ determines whether the tangent lines to the trigonometric function at different points within $\underline{\theta}_{lm}-\delta_{lm}-\psi_l$ and $\overline{\theta}_{lm}-\delta_{lm}-\psi_l$ have multiple intersections with the trigonometric function.
If $\mathcal{Z}_{\underline{\theta}_{lm}, \overline{\theta}_{lm}, \delta_{lm},\psi_l}$ is empty, it indicates that the curvature of the trigonometric function does not change.
Thus, no tangent lines will have another intersection with the function. In this case, the slopes of the $k^{th}$ tangent line (referred to as $m_k$) at an arbitrary point within the interval,  $(\phi_{lm},\cos(\phi_{lm}-\delta_{lm}-\psi_l)$, is the derivative of the trigonometric function at this point $\left(m_k=\frac{d(\cos(\theta_{lm}-\delta_{lm}-\psi_l))}{d\theta_{lm}}|_{\phi_{lm}}=-\sin(\phi_{lm}-\delta_{lm}-\psi_l)\right)$. 
The coordinate of the point itself gives us the offset, represented as $b_{int}$, of the tangent line which is equal to $\cos(\phi_{lm}-\delta_{lm}-\psi_l)$. 

Conversely, when $\left|\mathcal{Z}_{\underline{\theta}_{lm}, \overline{\theta}_{lm}, \delta_{lm},\psi_l}\right|\ge 1$, the tangent line to the cosine function at some point may intersect the cosine function again at a different point. To address this concern, we identify points, denoted as $\underline{\mathcal{R}}_{lm}$ and $\overline{\mathcal{R}}_{lm}$, at the boundaries of ranges for which tangent lines do not intersect the cosine function. $\underline{\mathcal{R}}_{lm}$ and $\overline{\mathcal{R}}_{lm}$ are shown by red and yellow stars, respectively, in Fig.~\ref{fig:Enclosing_points}.
Note that the voltage angle difference restriction, i.e., ${-}90^\circ \le \theta_{lm} \le 90^\circ$, guarantees that the sine and cosine functions experience at most one curvature sign change. Consequently, selecting tangent lines to the cosine function in the ranges $[\underline{\theta}_{lm}-\delta_{lm}-\psi_l,\underline{\mathcal{R}}_{lm}]$ and  $[\overline{\mathcal{R}}_{lm}, \overline{\theta}_{lm}-\delta_{lm}-\psi_l$] is a straightforward process since the curvature of the cosine function remains consistent, with the same sign, within these ranges. This method yields linear envelopes that provide a close outer approximation to the cosine function. The technique is applicable to any angle difference range, irrespective of the trigonometric function's curvature. Next, we will elaborate on the computation process for $\underline{\mathcal{R}}_{lm}$ and $\overline{\mathcal{R}}_{lm}$.

To determine $\underline{\mathcal{R}}_{lm}$ and $\overline{\mathcal{R}}_{lm}$ for the cosine function, we start by formulating the tangent lines to the cosine function that pass through the endpoints $\overline{\theta}_{lm}-\delta_{lm}-\psi_l$ and $\underline{\theta}_{lm}-\delta_{lm}-\psi_l$, respectively. The tangent line to the cosine function that pass through the point $\overline{\theta}_{lm}-\delta_{lm}-\psi_l$ is given by $\underline{F}_{tang}(\theta_{lm})$:
\begin{align}
\label{eq:tangent_line}
\nonumber & \underline{F}_{tang}(\theta_{lm})=-\sin(\theta_{lm}-\delta_{lm}-\psi_l)(\theta_{lm}-\overline{\theta}_{lm})\\&\qquad\qquad\qquad+\cos(\overline{\theta}_{lm}-\delta_{lm}-\psi_l).
\end{align}
We subsequently define an auxiliary function, denoted as $\underline{G}(\theta_{lm})$, which represents the difference between $\underline{F}_{tang}(\theta_{lm})$ and $\cos(\theta_{lm}-\delta_{lm}-\psi_l)$:
\begin{align}
\label{eq:tangent_line_aux}
\nonumber & \underline{G}(\theta_{lm})=\sin(\theta_{lm}-\delta_{lm}-\psi_l)(\theta_{lm}-\overline{\theta}_{lm})\\&\qquad -\cos(\overline{\theta}_{lm}-\delta_{lm}-\psi_l)+\cos(\theta_{lm}-\delta_{lm}-\psi_l).
\end{align}
The root of the derivative of $\underline{G}(\theta_{lm})$, i.e., the solution to $\frac{d\underline{G}(\theta_{lm})}{d \theta_{lm}} = 0$, within the interval between $\underline{\theta}_{lm}-\delta_{lm}-\psi_l$ and $\overline{\theta}_{lm}-\delta_{lm}-\psi_l$ corresponds to the value of $\underline{\mathcal{R}}_{lm}$.

Similarly, to find $\overline{\mathcal{R}}_{lm}$, we first formulate the tangent line to the cosine function that pass through the point $\underline{\theta}_{lm}-\delta_{lm}-\psi_l$:
\begin{align}
\label{eq:tangent_line2}
\nonumber & \overline{F}_{tang}(\theta_{lm})=-\sin(\theta_{lm}-\delta_{lm}-\psi_l)(\theta_{lm}-\underline{\theta}_{lm})\\&\qquad\qquad\qquad+\cos(\underline{\theta}_{lm}-\delta_{lm}-\psi_l).
\end{align}
Next, we define an auxiliary function, denoted as $\overline{G}(\theta_{lm})$, which represents the difference between $\overline{F}_{tang}(\theta_{lm})$ and $\cos(\theta_{lm}-\delta_{lm}-\psi_l)$:
\begin{align}
\label{eq:tangent_line_aux2}
\nonumber & \overline{G}(\theta_{lm})=\sin(\theta_{lm}-\delta_{lm}-\psi_l)(\theta_{lm}-\underline{\theta}_{lm})\\&\qquad -\cos(\underline{\theta}_{lm}-\delta_{lm}-\psi_l)+\cos(\theta_{lm}-\delta_{lm}-\psi_l).
\end{align}
Analogously to the discussion above, the root of the derivative of $\overline{G}(\theta_{lm})$, i.e., the solution to $\frac{d\overline{G}(\theta_{lm})}{d \theta_{lm}} = 0$, within the interval between $\underline{\theta}_{lm}-\delta_{lm}-\psi_l$ and $\overline{\theta}_{lm}-\delta_{lm}-\psi_l$, corresponds to the value of $\overline{\mathcal{R}}_{lm}$.

To locate the roots $\frac{d\underline{G}(\theta_{lm})}{d \theta_{lm}} = 0$ and $\frac{d\overline{G}(\theta_{lm})}{d \theta_{lm}} = 0$ within the interval $\left[\underline{\theta}_{lm}-\delta_{lm}-\psi_l, \overline{\theta}_{lm}-\delta_{lm}-\psi_l\right]$, we start by applying a bisection method to obtain a close initialization for a locally convergent Newton method to determine the precise values of the roots. The reason for employing the bisection method initially is that the Newton method may converge to solutions beyond the interval of interest for periodic functions like $\frac{d\underline{G}(\theta_{lm})}{d \theta_{lm}}$ and $\frac{d\overline{G}(\theta_{lm})}{d \theta_{lm}}$. The bisection method finds approximate solutions within the interval of interest that are refined with a Newton method. The slope of tangent line to the cosine function at $\underline{\mathcal{R}}_{lm}$ (referred to as $m_{k,\underline{\mathcal{R}}_{lm}}$), is the derivative of the cosine functions at the corresponding argument for $\underline{\mathcal{R}}_{lm}$, i.e.,  $\left(m_{k,\underline{\mathcal{R}}_{lm}}=-\sin(\phi_{lm,\underline{\mathcal{R}}_{lm}}-\delta_{lm}-\psi_l)\right)$. The coordinate of the point itself gives us the offset, represented as $b_{int,\underline{\mathcal{R}}_{lm}}$, of the tangent line which is equal to $\cos(\phi_{lm,\underline{\mathcal{R}}_{lm}}-\delta_{lm}-\psi_l)$. The slope and offset for $\overline{\mathcal{R}}_{lm}$ is computed similarly.

  		{\SetAlgoNoLine
	\begin{algorithm}[!t] % check
		\SetAlgoNoLine
		\small
		\caption{Finding tangent lines}
		\label{alg:algorithm1}
		% config
		\SetKwInOut{Input}{Input}
		\SetKwInOut{Output}{Output}
		\SetKwFunction{FMain}{Main}
		\SetKwFunction{FGen}{Generate}
		\DontPrintSemicolon
		\Function{Find\_Tangent($\overline{\theta}_{lm}$,$\underline{\theta}_{lm}$,$\delta_{lm}$,$\psi_l$,$N_{tan}$)}{
				$U \gets \overline{\theta}_{lm}-\delta_{lm} - \psi_l$,\quad $L \gets \underline{\theta}_{lm}-\delta_{lm} - \psi_l$.\; 
			Define $F(\theta_{lm})$
			and $\frac{dF(\theta_{lm})}{d \theta_{lm}}$.\;	
		$\mathcal{Z}_{\underline{\theta}_{lm}, \overline{\theta}_{lm}, \delta_{lm},\psi_l} \gets$ Find the roots of $\frac{dF(\theta_{lm})}{d \theta_{lm}}$ within $[L,U]$.
		
			\uIf{$\left|\mathcal{Z}_{\underline{\theta}_{lm}, \overline{\theta}_{lm}, \delta_{lm},\psi_l}\right| < 1$}{
			%At equally spaced points within $[L,U]$, compute tangent lines, i.e., $m_k$ and $b_{int}$.
   At $N_{tan}$ equally spaced points within $[L,U]$:\\
   \For{$i=1:N_{tan}$}{
            $m_{k,i}=-\sin(\phi_{lm,i}-\delta_{lm}-\psi_l)$
            $b_{int,i}=\cos(\phi_{lm,i}-\delta_{lm}-\psi_l)$
            }}
   \ElseIf{$\left|\mathcal{Z}_{\underline{\theta}_{lm}, \overline{\theta}_{lm}, \delta_{lm},\psi_l}\right| \ge 1$}{
                     $\underline{\mathcal{R}}_{lm}\gets\frac{dG(\theta_{lm},U)}{d \theta_{lm}}=0$\\
                     $\overline{\mathcal{R}}_{lm}\gets\frac{dG(\theta_{lm},L)}{d \theta_{lm}}=0$\\
Find the tangent lines to $\cos(\theta_{lm}-\delta_{lm}-\psi_l)$ at $\underline{\mathcal{R}}_{lm}$ and $\overline{\mathcal{R}}_{lm}$:

$m_{k,\underline{\mathcal{R}}_{lm}}=-\sin(\phi_{lm,\underline{\mathcal{R}}_{lm}}-\delta_{lm}-\psi_l)$\\
            $b_{int,\underline{\mathcal{R}}_{lm}}=\cos(\phi_{lm,\underline{\mathcal{R}}_{lm}}-\delta_{lm}-\psi_l)$\\

            $m_{k,\overline{\mathcal{R}}_{lm}}=-\sin(\phi_{lm,\overline{\mathcal{R}}_{lm}}-\delta_{lm}-\psi_l)$\\
            $b_{int,\overline{\mathcal{R}}_{lm}}=\cos(\phi_{lm,\overline{\mathcal{R}}_{lm}}-\delta_{lm}-\psi_l)$\\

					Equally divide the ranges $[\overline{\mathcal{R}}_{lm},U]$ and $[L,\underline{\mathcal{R}}_{lm}]$ into $N_{tan}$ segments

     \For{$i=1:N_{tan}$}{
            $m_{k,i}=-\sin(\phi_{lm,i}-\delta_{lm}-\psi_l)$
            $b_{int,i}=\cos(\phi_{lm,i}-\delta_{lm}-\psi_l)$
            }
     
			}	
				\Return {$m_k$} and $b_{int}$}%Tangent lines.}%
	\end{algorithm}
	}

%Once the tangent lines from $\underline{\theta}_{lm}-\delta_{lm}-\psi_l$ and $\overline{\theta}_{lm}-\delta_{lm}-\psi_l$ to $\overline{\mathcal{R}}$ and $\underline{\mathcal{R}}$ on the cosine function are derived, the tangent lines to the cosine function within $\left[\underline{\theta}_{lm}-\delta_{lm}-\psi_l,\underline{\mathcal{R}}\right]$ and $\left[\overline{\mathcal{R}},\overline{\theta}_{lm}-\delta_{lm}-\psi_l\right]$ intervals can be easily found as the curvature of the cosine function within these intervals does not change. 
After computing $\underline{\mathcal{R}}_{lm}$ and $\overline{\mathcal{R}}_{lm}$, we select equally spaced tangent lines to the cosine function within the intervals $\left[\underline{\theta}_{lm}-\delta_{lm}-\psi_l,\underline{\mathcal{R}}_{lm}\right]$ and $\left[\overline{\mathcal{R}}_{lm},\overline{\theta}_{lm}-\delta_{lm}-\psi_l\right]$.

Algorithm~\ref{alg:algorithm1} outlines the procedure for computing these convex envelopes using carefully selected tangent lines.
For notational convenience, define $L=\underline{\theta}_{lm}-\delta_{lm}-\psi_l$ and $U=\overline{\theta}_{lm}-\delta_{lm}-\psi_l$. Equations~\eqref{eq:tangent_lines_cosine_lower} and~\eqref{eq:tangent_lines_cosine_upper} define the tangent lines which form the lower and upper bounds, respectively, of the cosine envelope:
\begin{small}
\begin{subequations}
\label{eq:tangent_lines_cosine_lower}
\begin{align}
&\nonumber \text{If curvature changes within $[L,U]$ from positive to negative:}\\
&\underline{L}_{\cos,i}\!\!=\!\!\begin{cases}m_{k,i}(x-\!x_{0,i})+b_{int,i},~i=1,\ldots,N_{tan},~\text{if~} x\in [L,\underline{\mathcal{R}}_{lm}]\\m_{k,\underline{\mathcal{R}}_{lm}}(x-x_{0,\underline{\mathcal{R}}_{lm}})+b_{int,\underline{\mathcal{R}}_{lm}}~\text{if~} x\in [\underline{\mathcal{R}}_{lm},U]\end{cases}\\
&\nonumber \text{If curvature changes within $[L,U]$ from negative to positive:}\\
&\underline{L}_{\cos,i}\!\!=\!\!\begin{cases}m_{k,i}(x-\!x_{0,i})+b_{int,i},~i=1,\ldots,N_{tan},~\text{if~} x\in [\overline{\mathcal{R}}_{lm},U]\\m_{k,\underline{\mathcal{R}}_{lm}}(x-x_{0,\underline{\mathcal{R}}_{lm}})+b_{int,\underline{\mathcal{R}}_{lm}}~\text{if~} x\in [L,\overline{\mathcal{R}}_{lm}]\end{cases}\\
&\nonumber \text{If curvature does not change within $[L,U]$ and it is negative:}\\
&\underline{L}_{\cos,i}\!\!=\!\!\begin{cases}\frac{\cos(U)-\cos(L)}{U-L}(x-L)+\cos(L)~~~~~ \forall x\in [L,U], \end{cases}\\
&\nonumber \text{If curvature does not change within $[L,U]$ and is positive:}\\
&\underline{L}_{\cos,i}\!\!=\!\!\begin{cases}\frac{\cos(U)-\cos(L)}{U-L}(x-L)+\cos(L)~~~~~\forall x\in [L,U]\end{cases}
\end{align}
\end{subequations}
\end{small}

Similarly, the upper bound for the cosine envelope is:
\begin{small}
\begin{subequations}
\label{eq:tangent_lines_cosine_upper}
\begin{align}
&\nonumber \text{If curvature changes within $[L,U]$ from positive to negative:}\\
&\overline{L}_{\cos,i}\!\!=\!\!\begin{cases}m_{k,i}(x-\!x_{0,i})+b_{int,i},~i\!=\!1,\ldots,N_{tan},~\text{if~} x\in [\overline{\mathcal{R}}_{lm},U]\\m_{k,\overline{\mathcal{R}}_{lm}}(x-x_{0,\overline{\mathcal{R}}_{lm}})+b_{int,\overline{\mathcal{R}}_{lm}}~\text{if~} x\in [L,\overline{\mathcal{R}}_{lm}]\end{cases}\\
&\nonumber \text{If curvature changes within $[L,U]$ from negative to positive:}\\
&\overline{L}_{\cos,i}\!\!=\!\!\begin{cases}m_{k,i}(x-\!x_{0,i})+b_{int,i},~i\!=\!1,\ldots,N_{tan},~\text{if~} x\in [L,\underline{\mathcal{R}}_{lm}]\\m_{k,\underline{\mathcal{R}}_{lm}}(x-x_{0,\underline{\mathcal{R}}_{lm}})+b_{int,\underline{\mathcal{R}}_{lm}}~\text{if~} x\in [\underline{\mathcal{R}}_{lm},U]\end{cases}\\
&\nonumber \text{If curvature does not change within $[L,U]$ and is negative:}\\
&\overline{L}_{\cos,i}\!\!=\!\!\begin{cases}\frac{\cos(U)-\cos(L)}{U-L}(x-L)+\cos(L)~~~~~ \forall x\in [L,U] \end{cases}\\
&\nonumber \text{If curvature does not change within $[L,U]$ and is positive:}\\
&\overline{L}_{\cos,i}\!\!=\!\!\begin{cases}\frac{\cos(U)-\cos(L)}{U-L}(x-L)+\cos(L)~~~~~\forall x\in [L,U]\end{cases}
\end{align}%
\end{subequations}%
\end{small}%
where $\overline{L}_{\cos,i}$ and $\underline{L}_{\cos,i}$ are the $i^{th}$ tangent lines which upper and lower bound, respectively, the cosine function; $x$ equals $\theta_{lm} - \delta_{lm} - \psi_l$; $x_{0,i}$, $x_{0,\underline{\mathcal{R}}_{lm}}$, and $x_{0,\overline{\mathcal{R}}_{lm}}$ represent the horizontal coordinates of the corresponding points within their respective intervals; $b_{int,i}$, $b_{int,\underline{\mathcal{R}}_{lm}}$, and $b_{int,\overline{\mathcal{R}}_{lm}}$ denote the vertical coordinates of these same points, which determines the offset for the tangent lines; and 
$m_{k,i}$, $m_{k,\underline{\mathcal{R}}_{lm}}$, and $m_{k,\overline{\mathcal{R}}_{lm}}$ signify the slopes of the tangent lines at these points. Note that when the curvature of the trigonometric function does not change within an interval, either the lower or upper envelope can be defined as a line connecting both ends of the trigonometric function.
For instance, consider a specific angle $\left(\phi_{lm} - \delta_{lm} - \psi_l\right)$ within the interval $[L,\underline{\mathcal{R}}_{lm}]$. In this case, $m_{k,i}$ is equal to ${-}\sin\left(\phi_{lm} - \delta_{lm} - \psi_l\right)$, which represents the first derivative of the cosine function at $x_{0,i} = \phi_{lm} - \delta_{lm} - \psi_l$. The corresponding offset for this point is $b_{int,i} = \cos\left(\phi_{lm} - \delta_{lm} - \psi_l\right)$. Tangent lines for the lower and upper bounds of the sine function are formulated similarly.

Equation~\eqref{eq:convex_envelopes_sin&cos2} formulates the cosine function. To complete the full exposition, we similarly formulate the lower and upper envelopes for the sine function by defining a function $H(\theta_{lm})$. Here, the function $H(\theta_{lm})$ represents the difference between the trigonometric function $\sin(\theta_{lm} - \delta_{lm}-\psi_l)$ itself and the line which connects the endpoints of $\sin(\theta_{lm} - \delta_{lm}-\psi_l)$ at $\underline{\theta}_{lm}$ and $\overline{\theta}_{lm}$:
\begin{align}
\label{eq:auxi_equation2}
\nonumber \noindent &H(\theta_{lm})  = \sin(\theta_{lm} - \delta_{lm}-\psi_l) - \sin(\overline{\theta}_{lm}-\delta_{lm}-\psi_l)\\
\nonumber &\qquad\quad-\frac{\sin(\overline{\theta}_{lm}-\delta_{lm}-\psi_l)-\sin(\underline{\theta}_{lm}-\delta_{lm}-\psi_l)}{\overline{\theta}_{lm}-\underline{\theta}_{lm}}\\&\qquad\quad\times\left(\theta_{lm}-\overline{\theta}_{lm}-\psi_l\right).
\end{align}
The number of zeros of the first derivative of $H(\theta_{lm})$, i.e., the number of solutions for $\frac{dH(\theta_{lm})}{d \theta_{lm}} = 0$, is a key quantity to determine if the curvature of $\sin(\theta_{lm}-\delta_{lm}-\psi_l)$ changes between $\underline{\theta}_{lm}-\delta_{lm}-\psi_l$ and $\overline{\theta}_{lm}-\delta_{lm}-\psi_l$. 
If $\frac{dH(\theta_{lm})}{d \theta_{lm}} = 0$ has no solutions, then the curvature of the $\sin(\theta_{lm}-\delta_{lm}-\psi_l)$ does not change between $\underline{\theta}_{lm} - \delta_{lm} - \psi_l$ and $\overline{\theta}_{lm} - \delta_{lm} - \psi_l$. Our proposed QC relaxation uses envelopes $\left\langle \sin(\theta_{lm}-\delta_{lm}-\psi_l)\right\rangle^{S^{\prime}}$ formed by combining the tangent lines. Depending upon the sign of the sine function's curvature and the number of solutions to $\frac{dH(\theta_{lm})}{d \theta_{lm}} = 0$, there are different upper and lower bounds on the sine function's envelopes:
\begin{small}
\begin{subequations}
\label{eq:convex_envelopes_sin&cos3}
\begin{align}
%\raisetag{18pt} 
\label{eq:sin envelope11}
&\nonumber\qquad\qquad\quad\text{If $\frac{dH(\theta_{lm})}{d \theta_{lm}} = 0$ has one or more solutions: }\\
 &\left\langle \sin(\theta_{lm}-\delta_{lm}-\psi_l)\right\rangle^{S^{\prime}} =
\begin{cases}
\!\widecheck{S}^{\prime}\!:\begin{cases}
\widecheck{S}^{\prime}\! \le \overline{L}_{\sin,i},~i=1,\ldots,N_{tan}\\
\widecheck{S}^{\prime}\! \ge \underline{L}_{\sin,i},~i=1,\ldots,N_{tan}\end{cases}
\end{cases}\\
%\raisetag{18pt} 
\label{eq:sine envelope12}
&  \nonumber\quad\text{If $\frac{dH(\theta_{lm})}{d \theta_{lm}} = 0$ has no solutions \& curvature sign is negative: }\\
&\left\langle\sin(\theta_{lm}-\delta_{lm}-\psi_l)\right\rangle^{S^{\prime}} =
\begin{cases}
\!\widecheck{S}^{\prime}\!:\begin{cases}\widecheck{S}^{\prime}\! \le 
\overline{L}_{\sin,i},~i=1,\ldots,N_{tan}\\
\widecheck{S}^{\prime}\! \ge \underline{L}_{\sin,i},~i=1\end{cases}
\end{cases}\\
%\raisetag{18pt} 
\label{eq:sine envelope13}
& \nonumber\quad\text{If $\frac{dH(\theta_{lm})}{d \theta_{lm}} = 0$ has no solutions \& curvature sign is positive: } \\
&\left\langle\sin(\theta_{lm}-\delta_{lm}-\psi_l)\right\rangle^{S^{\prime}} =
\begin{cases}
\!\widecheck{S}^{\prime}\!:\begin{cases}\widecheck{S}^{\prime}\! \le 
\overline{L}_{\sin,i},~i=1\\
\widecheck{S}^{\prime}\! \ge \underline{L}_{\sin,i},~i=1,\ldots,N_{tan}\end{cases}
\end{cases}
\end{align}
\end{subequations}%
\end{small}%
where $\overline{L}_{\sin,i}$ and $\underline{L}_{\sin,i}$ are the $i^{th}$ tangent lines which upper and lower bound, respectively, the sine function. 

Equations~\eqref{eq:tangent_lines_sine_lower} and~\eqref{eq:tangent_lines_sine_upper} mathematically represent the tangent lines for the lower and upper bounds of the sine function, respectively: 
\begin{small}
\begin{subequations}
\label{eq:tangent_lines_sine_lower}
\begin{align}
&\nonumber \text{If curvature changes within $[L,U]$ from positive to negative:}\\
&\underline{L}_{\sin,i}\!\!=\!\!\begin{cases}m_{k,i}(x-\!x_{0,i})+b_{int,i},~i=1,\ldots,N_{tan},~~\text{if~} x\in [L,\underline{\mathcal{R}}_{lm}]\\m_{k,\underline{\mathcal{R}}_{lm}}(x-x_{0,\underline{\mathcal{R}}_{lm}})+b_{int,\underline{\mathcal{R}}_{lm}}~\text{if~} x\in [\underline{\mathcal{R}}_{lm},U]\end{cases}\\
&\nonumber \text{If curvature changes within $[L,U]$ from negative to positive:}\\
&\underline{L}_{\sin,i}\!\!=\!\!\begin{cases}m_{k,i}(x-\!x_{0,i})+b_{int,i},~i=1,\ldots,N_{tan},~\text{if~} x\in [\overline{\mathcal{R}}_{lm},U]\\m_{k,\underline{\mathcal{R}}_{lm}}(x-x_{0,\underline{\mathcal{R}}_{lm}})+b_{int,\underline{\mathcal{R}}_{lm}}~\text{if~} x\in [L,\overline{\mathcal{R}}_{lm}]\end{cases}\\
&\nonumber \text{If curvature does not change within $[L,U]$ and is negative:}\\
&\underline{L}_{\sin,i}\!\!=\!\!\begin{cases}\frac{\sin(U)-\sin(L)}{U-L}(x-L)+\sin(L)~~~~~ \forall x\in [L,U], \end{cases}\\
&\nonumber \text{If curvature does not change within $[L,U]$ and is positive:}\\
&\underline{L}_{\sin,i}\!\!=\!\!\begin{cases}\frac{\sin(U)-\sin(L)}{U-L}(x-L)+\sin(L)~~~~~\forall x\in [L,U]\end{cases}
\end{align}
\end{subequations}%
\end{small}%
Similarly, the upper bound for the sine function can be represented as follows:
\begin{small}
\begin{subequations}
\label{eq:tangent_lines_sine_upper}
\begin{align}
&\nonumber \text{If curvature changes within $[L,U]$ from positive to negative:}\\
&\overline{L}_{\sin,i}\!\!=\!\!\begin{cases}m_{k,i}(x-\!x_{0,i})+b_{int,i},~i=1,\ldots,N_{tan},~\text{if~} x\in [\overline{\mathcal{R}}_{lm},U]\\m_{k,\overline{\mathcal{R}}_{lm}}(x-x_{0,\overline{\mathcal{R}}_{lm}})+b_{int,\overline{\mathcal{R}}_{lm}}~\text{if~} x\in [L,\overline{\mathcal{R}}_{lm}]\end{cases}\\
&\nonumber \text{If curvature changes within $[L,U]$ from negative to positive:}\\
&\overline{L}_{\sin,i}\!\!=\!\!\begin{cases}m_{k,i}(x-\!x_{0,i})+b_{int,i},~i=1,\ldots,N_{tan},~\text{if~} x\in [L,\underline{\mathcal{R}}_{lm}]\\m_{k,\underline{\mathcal{R}}_{lm}}(x-x_{0,\underline{\mathcal{R}}_{lm}})+b_{int,\underline{\mathcal{R}}_{lm}}~\text{if~} x\in [\underline{\mathcal{R}}_{lm},U]\end{cases}\\
&\nonumber \text{If curvature does not change within $[L,U]$ and is negative:}\\
&\overline{L}_{\sin,i}\!\!=\!\!\begin{cases}\frac{\sin(U)-\sin(L)}{U-L}(x-L)+\sin(L)~~~~~ \forall x\in [L,U] \end{cases}\\
&\nonumber \text{If curvature does not change within $[L,U]$ and is positive:}\\
&\overline{L}_{\sin,i}\!\!=\!\!\begin{cases}\frac{\sin(U)-\sin(L)}{U-L}(x-L)+\sin(L)~~~~~\forall x\in [L,U]\end{cases}
\end{align}
\end{subequations}%
\end{small}%
where $x$ equals $\theta_{lm} - \delta_{lm} - \psi_l$; $x_{0,i}$, $x_{0,\underline{\mathcal{R}}_{lm}}$, and $x_{0,\overline{\mathcal{R}}_{lm}}$ represent the horizontal coordinates of the corresponding points within their respective intervals; $b_{int,i}$, $b_{int,\underline{\mathcal{R}}_{lm}}$, and $b_{int,\overline{\mathcal{R}}_{lm}}$ denote the vertical coordinates of these points within their corresponding intervals, determining the offset for the tangent lines; and 
$m_{k,i}$, $m_{k,\underline{\mathcal{R}}_{lm}}$, and $m_{k,\overline{\mathcal{R}}_{lm}}$ signify the slopes of the tangent lines at specific points within their respective intervals.
\fi

\bibliographystyle{IEEEtran}
\ifarxiv
% \IEEEtriggeratref{20}
\else
 \IEEEtriggeratref{40}
\fi
\bibliography{refs}

\end{document}